\documentclass[]{article}

\usepackage{amsmath}
\usepackage{amssymb}
\usepackage{enumerate}
\usepackage{url}
\usepackage{amsthm}
\usepackage{graphicx}
\usepackage{subfigure}
\usepackage{multirow}

\newcommand{\s}{\mathfrak{s}}
\newcommand{\Rmap}{R}
\newcommand{\TM}{\mathcal{T}}
\newcommand{\RM}{\mathcal{R}}
\newcommand{\tRmap}{\tilde{R}}
\newcommand{\RR}{\mathbb{R}}
\newcommand{\N}{\mathbb{N}}
\newcommand{\systwreset}{\eqref{eq:general_system}-\eqref{eq:reset_condition}}
\newcommand{\D}{\mathbb{D}}

\newcommand{\bz}{\bar{z}}
\newcommand{\Stran}{\Sigma^\text{NS}}
\newcommand{\Stan}{\Sigma^\text{S}}

\newcommand{\Tt}{\mathbb{T}_T}

\newcommand{\tP}{\tilde{P}}
\newcommand{\Gc}{\mathcal{G}}
\newcommand{\td}{\tilde{\delta}}

\newcommand{\x}{{\mathbf x}}
\newcommand{\y}{{\mathbf y}}
\newcommand{\maps}{\eqref{eq:map_P1}--\eqref{eq:map_P3}}

\newtheorem{theorem}{Theorem}
\newtheorem{definition}{Definition}
\newtheorem{remark}{Remark}
\newtheorem{example}{Example}

\graphicspath{
{figures/}
{figures/period_adding_like_Delta03/A0d119/}
{figures/period_adding_like_Delta03/A0d179/}
{figures/period_adding_like_Delta03/A0d245/}
{figures/period_adding_like_Delta03/A0d512/}
{figures/period_adding_like_Delta03/periods/}
{figures/period_adding_like_Delta03/trajectorya/}
{figures/period_adding_like_Delta03/trajectoryb/}
{figures/period_adding_like_Delta03/trajectoryc/}
{figures/period_adding_like_Delta03/trajectoryd/}
{figures/period_adding_like_Delta03/satisfying_Gam/}
{figures/tangent_bif_coexistence/}
{figures/2dbif_T0d5_Delta03/}
}

\begin{document}
\newcommand{\biftanSn}{\eqref{eq:bc1tan}-\eqref{eq:bc11}}
\newcommand{\biftanSnn}{\eqref{eq:bc1tan2}, \eqref{eq:bc1}-\eqref{eq:bc11}}
\newcommand{\biftranSn}{\eqref{eq:bc1tran}-\eqref{eq:bc22}}
\newcommand{\biftranSnn}{\eqref{eq:bc1tran2}, \eqref{eq:bc2}-\eqref{eq:bc22}}
\author{Albert Granados\footnotemark[2] and Gemma
Huguet\footnotemark[3]}
\title{Gluing and grazing bifurcations in periodically forced $2$-dimensional
integrate-and-fire models\thanks{This work has been partially supported by  MINECO
MTM2015-65715-P Spanish grant and Marie Curie FP7 COFUND \O rsted
fellowship. We acknowledge the use of the UPC Dynamical Systems
group's cluster for research computing
({https://dynamicalsystems.upc.edu/en/computing/})}}
\date{}
\maketitle
\renewcommand{\thefootnote}{\fnsymbol{footnote}}
\footnotetext[2]{algr@dtu,dk, Department of Applied Mathematics and
Computer Science, Technical University of Denmark, Building 303B, 2800
Kgns. Lyngby, Denmark.}
\footnotetext[3]{Department de Matem\`atiques, Universitat Polit\`ecnica de
Catalunya, Av. Diagonal 647, 08028 Barcelona, Spain.}
\renewcommand{\thefootnote}{\arabic{footnote}}
\begin{abstract}
In this work we consider a general class of $2$-dimensional hybrid
systems. Assuming that the system possesses an attracting equilibrium
point, we show that, when periodically driven with a square-wave
pulse, the system possesses a periodic orbit which may undergo
smooth and nonsmooth grazing bifurcations.  We perform a
semi-rigorous study of the existence of periodic orbits for a
particular model consisting of a leaky integrate-and-fire model with a
dynamic threshold. We use the stroboscopic map, which in this context is a
$2$-dimensional piecewise-smooth discontinuous map. For some parameter
values we are able to show that the map is a quasi-contraction
possessing a (locally) unique maximin periodic orbit. We complement
our analysis using advanced numerical techniques to provide a complete
portrait of the dynamics as parameters are varied. 
We find that for some regions of the parameter space the model undergoes a cascade of 
gluing bifurcations, while for others the model shows multistability 
between orbits of different periods.

\end{abstract}
\section{Introduction}\label{sec:intro}
Integrate-and-fire systems are hybrid systems that combine continuous dynamics with discrete resets that occur whenever the 
variables of the system satisfy a given condition (that defines a threshold). Such systems are widely used in
neuroscience to model the dynamics of neuron's membrane potential, as the continuous dynamics models subthreshold behaviour (corresponding
to the input integration) and resets 
represent neuron spikes (characteristic rapid changes in membrane potential). They can be
seen as simplified versions of slow-fast systems, as the resets replace large amplitude oscillations that occur at 
a much faster time-scale.

Examples of such systems range from simple one-dimensional models as
the leaky integrate-and-fire \cite{Izhikevich07}, which models simple repetitive spiking, to
nonlinear 2-dimensional ones exhibiting more complicated behaviour, 
such as the Izhikevich quadratic model \cite{Izhi03}, 
or the adaptive exponential model \cite{BreGer05}. 
Of special interest for this paper are 2-dimensional models consisting of an integrate-and-fire model with a dynamic threshold -- the threshold is 
treated as a variable with nonlinear dynamics.  These type of systems have been used to 
model spike threshold variability observed in different areas of the nervous system \cite{Platkiewicz10}, in particular for 
phasic neurons (those that do not respond repetitively to steady or slowly varying inputs) in the auditory brainstem~\cite{HiggsS11, MenHugRin12,HugMenRin17}. 

A general framework to study the dynamics of hybrid systems becomes
difficult to obtain, even when the input currents are assumed to be
constant (the system remains autonomous), mainly because they are
discontinuous due to the reset condition. One of the most common
strategies in the nonsmooth literature
(see~\cite{MakLam12,BerBudChaKow08}) is to smooth the dynamics by
considering the so-called \emph{impact map} (also known in neuroscience as \emph{firing
phase map} or \emph{adaptation map}) defined on the threshold where the reset
condition is applied \cite{TouboulB09, Rubinetal17a, Coombes12}. However, the impact
map does not allow one to study itineraries or trajectories that do
not hit the threshold and has some domain restrictions. 
 
In more realistic situations one considers
periodic inputs, making the analysis more complicated. 
Indeed, in the non-autonomous case, even one-dimensional
integrate-and-fire systems exhibit very rich
dynamics \cite{KeeHopRin81,TieFellSej02, Coombes12}.  Recent works show how
theory for nonsmooth systems can be used to obtain model-independent
general results~\cite{GraKruCle14,GraKru15}. However, these are
limited to one-dimensional systems exhibiting simple subthreshold
dynamics, as they are based on the theory for circle maps
(see~\cite{GraAlsKru17} for a recent review).

In this work we study 2-dimensional hybrid systems 
subject to periodic forcing. In particular, we consider an input consisting of a
square-wave pulse of period $T$, which can be seen as a simplified model for a 
periodically varying synaptic current in neuroscience \cite{ErmentroutTerman2010} 
while is widely used in electronics (PWM), amongst others.
Our goal is to provide a description of the dynamics
to determine the firing patterns that arise in the forced system.
Assuming that the unforced system possesses an attracting equilibrium point, 
we show that, when periodically
driven, the system possesses a $T$-periodic orbit which may undergo
smooth-grazing or nonsmooth-grazing bifurcations as the amplitude of the 
forcing increases and collides with the threshold. 
To study these bifurcations, we use the stroboscopic map (or time-$T$ map),
which becomes a $2$-dimensional piecewise-smooth discontinuous map. Indeed, the map 
is defined differently depending on the number of times the
corresponding trajectory of
the time-continuous system hits the threshold, thus
showing discontinuities along the so-called \emph{switching manifolds}.
An orbit of the time-continuous system undergoing a smooth or
a nonsmooth grazing bifurcation with the threshold corresponds to a fixed point
of the stroboscopic map undergoing a border collision bifurcation
when hitting the switching manifold ,
which we study in detail.

Beyond the fixed points we study other periodic orbits of the
stroboscopic map ($nT$-periodic orbits of the time-continuous
system), as well as their 
itine\-raries (sequence of regions in the domain of the stroboscopic map
visited by the periodic orbits), and their bifurcations.
Notice that in this case, the theory 
for smooth maps does not apply to describe its dynamical properties. Unfortunately, there are little theoretical 
results to describe these orbits for the cases of maps of dimension
higher than 1, in which case one relies on classical results for
circle maps~\cite{AlsLliMis00}. In this paper we 
recover a result in this direction by Gambaudo \emph{et al.} in the 80's~\cite{GamTre85,GamGleTre84}. The theorem establishes conditions for the existence of periodic orbits for a 
piecewise continuous map of any dimension 
and provides properties on the sequence of regions visited. 
We apply this result to a particular model, a leaky integrate-and-fire model with a dynamic threshold.
By means of semi-rigorous numerical arguments we can prove that for certain parameter values the stroboscopic map
becomes a quasi-contraction possessing {\em maximin} (locally) unique periodic
orbits. That is, their symbolic itineraries are contained in the Farey
tree of symbolic sequences~\cite{GraAlsKru17}.
In parallel, we use advanced numerical techniques to provide a complete portrait 
of the dynamics as parameters are varied. 
Numerically, we find that for
certain parameters the model undergoes a period-adding bifurcation (an infinite cascade of gluing
bifurcations~\cite{GamGleTre88}, up to our numerical accuracy), while 
for others the model shows multistability between orbits of different periods.
Our study allow us to assess the scope of the existing theoretical results.

The paper is organized as follows. In Section~\ref{sec:general_setting} we present the general setting for 
hybrid 2-dimensional systems. In Section~\ref{sec:strobo_map} we introduce the stroboscopic map, 
which is a piecewise-smooth
discontinuous 2-dimensional map. In Section~\ref{sec:bif_of_fp} we describe the border collision bifurcations of the fixed points, a type of 
bifurcation that can only occur in piecewise-defined maps and in Section~\ref{sec:periodic_orbits} we present existing theoretical results for the existence of 
periodic orbits of a piecewise continuous map which step onto different regions. 
In Section~\ref{sec:numerical_results} we use the previous results in combination with numerical methods 
to describe the dynamics of a leaky integrate-and-fire model with a dynamic threshold. We modify parameters to illustrate different dynamical regimes exhibited by the model. 
Finally, the Appendix includes the details of the numerical methods used to perform the computations along the paper.

\section{General setting}\label{sec:general_setting}

Let us consider the system
\begin{equation}
\dot{z}=f(z)+vI(t), \quad z \in \RR^2
\label{eq:general_system}
\end{equation}
with $v\in\RR^2$ (typically it will be chosen $v=(1,0)^T$), $f:\RR^2 \rightarrow \RR^2$ a smooth enough function
and $I(t)$ a
$T$-periodic square wave given by
\begin{equation}
I(t)=\left\{
\begin{aligned}
&A&&\text{if }t\in\left(nT,nT+dT\right]\\
&0&&\text{if }t\in(nT+dT,(n+1)T].\\
\end{aligned}\right.,\,n\in\mathbb{N}
\label{eq:forcing}
\end{equation}
Let us also consider a {\em threshold} manifold $\TM$ in $\RR^2$ given by
\begin{equation*}
\TM=\left\{ (x,y)\in\RR^2\,|\,h(x,y)=0 \right\},
\end{equation*}
where $h:\RR^2\rightarrow\RR$ is a smooth function. We then
submit system~\eqref{eq:general_system} to the following reset
condition: whenever a trajectory reaches the threshold manifold $\TM$ at a
time $t=t_*$, the variables of the system are updated to a certain
value, i.e.
\begin{equation}
h(z(t_*))=0\longrightarrow
z(t_*^+)=\Rmap(z(t_*)),
\label{eq:reset_condition}
\end{equation}
where $\Rmap$ is a reset (smooth) map:
\begin{equation}
\begin{array}{cccc}
\Rmap:&\TM&\longrightarrow& \RR^2.
\end{array}
\label{eq:reset_map}
\end{equation}
We call
\begin{equation*}
\RM:=R(\TM)
\end{equation*}
the reset manifold.  The reset condition~\eqref{eq:reset_condition} is
applied whenever a trajectory collides with the threshold manifold
$\TM$. Following the terminology in neuroscience, from now on, when
this occurs we say that system~\systwreset{} exhibits a spike.
Although, these spikes introduce discontinuities to the trajectories
of system~\eqref{eq:general_system}, they are all well defined, as one
just needs to apply the map $\Rmap$ whenever the threshold is reached.
This induces a flow,
\begin{equation*}
\phi(t;t_0,z_0),\quad \phi(t_0;t_0,z_0)=z_0,
\end{equation*}
which, provided that sliding cannot not occur along the threshold
manifold, is well defined. The flow however is discontinuous whenever
$\TM$ is reached and nonsmooth when the pulse $I(t)$ is enabled or
disabled.

We are going to assume that for $A=0$, sytem~\eqref{eq:general_system}
has an equilibrium point $z^{*} \in \RR^2$ (see {\em H.1} below). Then, we can define the subthreshold domain as
\begin{equation}
\D= \Big\{ 
z\in\RR^2\,|\, h(z) \cdot h(z^*)>0
\Big\}. 
\end{equation}
The subthreshold domain contains all points in one side of $\TM$ (the same side that contains the equilibrium point for $A=0$).
In many practical applications we will restrict this domain to 
points in the region delimited by the manifolds $\TM$ and $\RM$.

Given $z\in \D$, we will say that its
trajectory is subthreshold if $\phi(t;t_0,z)\in\D$ for all
$t \ge t_0$.  In particular, an invariant set is
subthreshold if it is contained in $\D$.\\

We assume that, for $A=0$, system~\eqref{eq:general_system}
\begin{enumerate}[\itshape H.1]
\item possesses an attracting equilibrium point $z^* \in \RR^2$,
\item for any $z\in\D$, trajectories are subthreshold, i.e. they do not exhibit spikes.
\end{enumerate}

\begin{remark}\label{rem:H2_not-required}
Hypothesis H.2 could be removed at the price of increasing the
complexity of the mathematical analysis (see
Remark~\ref{rem:reorder_strobo}),
but we decided to keep it in order to make the presentation clearer.
Moreover, we want to emphasize that it is a realistic assumption.
Indeed, in Section~\ref{sec:numerical_results} we consider an 
integrate-and-fire model which does satisfy hypothesis H.2. 
\end{remark}

\section{The stroboscopic map}\label{sec:strobo_map}
\subsection{Definition}\label{sec:strobo_def}
\begin{figure}
\begin{center}
\includegraphics[width=0.8\textwidth]{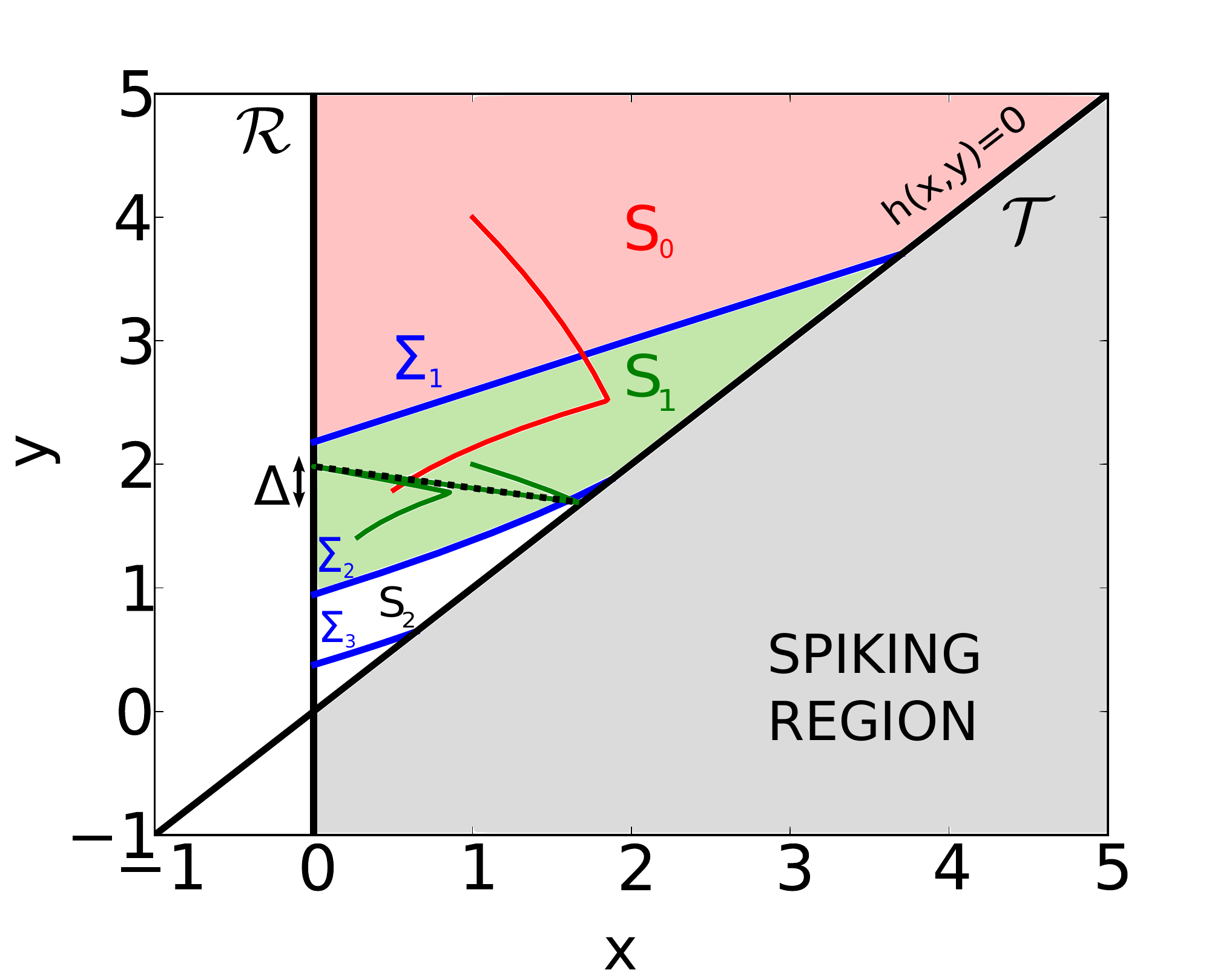}
\end{center}
\caption{Trajectories of system~\eqref{eq:system_visual}-\eqref{eq:reset_ex} with initial conditions in sets $S_0$ (red) and
$S_1$ (green), plotted on top of the pieces $S_0$, $S_1$, $S_2$ where the 
stroboscopic map is defined. Boundaries $\Sigma_1$, $\Sigma_2$ and $\Sigma_3$ are computed
using the algorithm described in Appendix~\ref{ap:sigma}. 
Solid black curves correspond to the threshold manifold $\TM$ and the reset manifold $\RM$.
Parameter values for the system are $c=0.53$,
$V_0=0.1$, $A=2$, $\tau=2$, $T=3$, $\Delta=0.3$ and $b=0.1$.}

\label{fig:sets_S0-S1}
\end{figure}

Since we consider a $T$-periodic forcing $I(t)$ (see Equation~\eqref{eq:forcing}), we follow~\cite{GraKruCle14}
and consider the stroboscopic map $\s: \D \rightarrow \D$ defined as:
\begin{equation*}
\s(z)=\phi(t_0+T;t_0,z).
\end{equation*}
Note that system~\systwreset{} is non-autonomous and, therefore, the stroboscopic
map depends on the initial time, $t_0$. However, as $t_0$ provides a family of conjugated stroboscopic
maps, we can assume $t_0=0$ from now on and abusing notation we write
$\s(z)=\phi(T;z)=\phi(T;0,z)$.\\
As detailed below, depending on the number of spikes exhibited by a solution
$\phi(t;z)$ for $t\in[0,T)$, the stroboscopic map becomes a
different combination of smooth maps given by integrating
Equation~\eqref{eq:general_system} and applying the reset
map~\eqref{eq:reset_map}. Hence, this is a piecewise smooth map. More precisely, let us define the sets (see Figure~\ref{fig:sets_S0-S1}):
\begin{equation}
S_n=\Big\{ z\in \D\, |\, \phi(t;z)\text{ exhibits $n$ spikes for } 0 < t\le dT \Big\},\;n\ge 0.
\label{eq:sets}
\end{equation}
Then, when restricted to $S_n$, the map $\s$ becomes a fixed
combination of maps alternating the integration of system \eqref{eq:general_system} and 
the application of the reset map $R$. Hence, $\s$  restricted to $S_n$ is as smooth as the map $R$ and the flow of the vector 
field $f$, as the sequence of impacts is fixed in $S_n$.
\begin{remark}
Notice that if $\phi(t;z)$ touches the threshold manifold $\TM$ at $t=dT$, the point $z$ will belong to $S_0$ or $S_1$
depending whether we apply the reset condition or not.
\end{remark}

In this work we will mainly focus on orbits involving the sets $S_0$ and $S_1$. For this reason, we first
show how to define $\s$ in these sets. We refer to
Figure~\ref{fig:sets_S0-S1} in order to illustrate what follows.\\
If $z\in S_0$, no spike occurs and $\s(z)$ becomes
\begin{equation}\label{eq:strobo_s0}
\s(z)=\s_0(z):=\varphi_0(T-dT;\varphi_A(dT;z)),
\end{equation}
where $\varphi_A$ is the flow associated with the system
$\dot{z}=f(z)+vA$.\\
In order to define the stroboscopic map in $S_1$, we consider the maps
\begin{equation}
\begin{array}[]{cccc}
P_1:& \cup_{n\ge1} S_n &\longrightarrow&\TM\times\Tt\\
&z&\longmapsto&(\varphi_A(t_*;z),t_*)
\end{array}
\label{eq:map_P1}
\end{equation}
\begin{equation}
\begin{array}[]{cccc}
\tRmap:&\TM\times\Tt&\longrightarrow&\mathcal{R}\times\Tt\\
&(z,t)&\longmapsto&(\Rmap(z),t)
\end{array}
\label{eq:map_R}
\end{equation}
\begin{equation}
\begin{array}[]{cccc}
\tP_2:&\mathcal{R}\times\Tt&\longrightarrow&\RR^2\\
&(z,t)&\longrightarrow&\varphi_A\left(dT-t;z\right)
\end{array}
\label{eq:map_P2}
\end{equation}
\begin{equation}
\begin{array}{cccc}
P_3:&\RR^2&\longrightarrow&\RR^2\\
&z&\mapsto&\varphi_0(T-dT;z)
\end{array}
\label{eq:map_P3}
\end{equation}
where $\Tt:=\RR/T\mathbb{Z}$ and the symbol $\tilde{}$ emphasizes that the
map is defined in the extended state space (time is taken as a variable).\\
The map $P_1$ sends points in $\D$ to the threshold $\TM$ by integrating the flow $\varphi_A$; it returns the hitting point on $\TM$
and the time $t_*$ needed by the trajectory to reach $\TM$.
In principle, to relate $P_1$ with system~\systwreset{}, its domain should be those points in $\D$ for
which $0 < t_*\le dT$, which is contained in $\cup_{n\ge1} S_n$. However,
the map $P_1$ can be extended to all points in $\D$ whose flow
$\varphi_A(t;z)$ reaches the threshold for some $t_*>0$,
independently on whether $t_*\le dT$ or not (see Section~\ref{sec:extension} for more
details). 
This extension becomes specially useful for numerical purposes as well as to provide insight
into the dynamics of the stroboscopic map near the switching manifolds.\\
The map $\tRmap$ is the reset map defined in~\eqref{eq:reset_map}
carrying on time. The map $\tP_2$ integrates the flow $\varphi_A$ with initial condition at the reset
manifold $\mathcal{R}$ for the remaining time until $t=dT$. Note that, similarly
as for $P_{1}$, $\tP_2$ can also be extended outside its natural domain
by letting $t<0$ (see Section~\ref{sec:extension} for more
details). Finally, the map $P_3$ is a truly stroboscopic map,
which integrates the flow $\varphi_0$ for a fixed time $T-dT$.\\
Since by hypothesis \emph{H.2}, spikes are only possible for
$A>0$ (that is, $0 < t_*\le dT$), then, for $z\in S_1$, the stroboscopic map becomes
\begin{equation}
\s(z)=\s_1(z):=P_3\circ \tP_2\circ \tRmap\circ P_1(z).
\label{eq:s_in_S1}
\end{equation}
By considering $\tP_1$ and $\tP_3$ the extended versions (to
$\RR^2\times\Tt$) of the maps $P_1$ and $P_3$ and recalling that
spikes can only occur for $0 < t\le dT$, if $z\in S_n$, $n \geq 1$ the stroboscopic map
becomes 
\begin{equation}
\s(z)= \s_n(z):=\tP_3\circ \tP_2\circ \left(\tRmap\circ \tP_1\right)^n(z).
\label{eq:strobo_general}
\end{equation}

Then, the stroboscopic map can be written as the piecewise smooth discontinuous map:
\begin{equation*}
\s(z)=\left\{
\begin{aligned}
&\s_0(z)&&\text{if }z\in S_0\\
&\s_n(z)&&\text{if }z\in S_n, n \geq 1
\end{aligned}
\right.
\end{equation*}

\begin{remark}\label{rem:reorder_strobo}
If one allows the system to exhibit spikes for $A=0$ (i.e, hypothesis
\emph{H.2} is not satisfied), then the stroboscopic map can be
similarly defined by reordering accordingly the sequence of maps
$\tP_i$ and $\tRmap$ in Equation~\eqref{eq:strobo_general}.
Moreover, the map $\s$ restricted to $S_n$ is still smooth as
long as the sequence of maps is kept constant.
\end{remark}
\begin{remark}
The stroboscopic map is discontinuous even if one identifies the
threshold and the reset manifolds: $\TM\sim R$. Although this would make
trajectories of the flow continuous, the vector field~\eqref{eq:general_system}
does not necessary coincide at the manifolds $\TM$ and $R$ and
hence the stroboscopic map would still be discontinuous.
\end{remark}
Let us now study the border, $\Sigma_1 \subset \D$, that separates the sets $S_0$
and $S_1$ and hence becomes a switching manifold of the stroboscopic
map $\s$ (see Figure~\ref{fig:sets_S0-S1}). This border is formed by the union of points whose
trajectories graze the threshold manifold $\TM$. Such a grazing
can occur in two different ways defining two different types of
points in $\Sigma_1$  
:
\begin{enumerate}[\itshape i)]
\item Smooth Grazing: points whose trajectory is tangent to $\TM$.
\item Nonsmooth Grazing: points whose trajectory is transversal to $\TM$ exactly for
$t=dT$.
\end{enumerate}
Provided that trajectories can only reach the threshold when $A>0$, (condition
\emph{H.2}), nonsmooth grazing can only occur for $t=dT$, at times when the pulse $I(t)$ is disabled. 
However, trajectories may exhibit tangent grazing
for $0 <t \leq dT$.\\
As mentioned above, the switching manifold $\Sigma_1$ can be split in
two pieces according to {\em i)} and {\em ii)}:
\begin{equation*}
\Sigma_1=\Stan_1\cup\Stran_1,
\end{equation*}
where
\begin{equation*}
\begin{aligned}
\Stan_1&= \Big\{z\in \D\,|\,
z=\varphi_A(t;z_0)\,,t\in[0,t_*], 0 < t_{*} \leq dT \, \textrm{ where } z_0 \textrm{ and } t_*\\&\quad  \textrm{ are s.t. }
h(\varphi_A(t_*;z_0))=0,
\nabla h(\varphi_A(t_*;z_0))\cdot
\frac{d}{dt}\varphi_A(t_*;z_0)=0 \Big\},
\end{aligned}
\end{equation*}
and  
\begin{equation*}
\Stran_1=\left\{ z\in \D\,|\,h(\varphi_A(dT;z))=0
\right\}.
\end{equation*}
\begin{remark}
Similarly, one can define the boundaries $\Stan_i$ and $\Stran_i$ with
$i>1$, which separate sets exhibiting more than one spike.
\end{remark}

\subsection{Virtual extension and contractiveness of the stroboscopic map}\label{sec:extension}
The maps $\s_0$ and $\s_1$ can (in some cases) be extended to their
``virtual'' domains, $S_1$ and $S_0$, respectively. In Section~\ref{sec:fixed_points_numeric} we will
show that the extended maps will be used to numerically compute feasible 
fixed points and bifurcation curves by means of a Newton method. Moreover, virtual extensions provide insight
into the dynamics of the map in the actual domain. For instance,
virtual attracting fixed points of the map suggest that the dynamics
in the actual domain pushes trajectories towards them and therefore
towards the switching manifold $\Sigma_1$.

Clearly, by ignoring the reset condition, one can always smoothly extend
$\s_0$ to $S_1$. That is, if $z\in S_1$, then we extend $\s_0$ to
$S_1$ by setting $\s_0(z)=\varphi_0(t-dT;\varphi_A(dT;z))$, which is well
defined. In words, ``keep integrating
system~\eqref{eq:general_system} with $I=A$ for a time $dT$ even if
the trajectory hits the threshold manifold $\TM$''.

Under certain conditions, one can also extend the map $\s_1$ to $S_0$.
Let $z\in S_0$ and assume that there exists $t_*>dT$ such that
$\varphi_A(t_*;z)\in\TM$. Then, although $z\notin S_1$, $\s_1$ is also
well defined at such a point by letting $t_*>dT$ in the definition of $P_1$ in \eqref{eq:map_P1} and using
$t>dT$ when applying the map $\tP_2$ defined in \eqref{eq:map_P2}, which  will consist of
integrating the flow $\varphi_A$ backwards for $t=|dT-t_*|$ .  In
words, ``keep integrating system~\eqref{eq:general_system} with $I=A$
as much time as needed until the trajectory hits the threshold
manifold $\TM$ and reset. Then, integrate the flow of system~\eqref{eq:general_system} with $I=0$ backwards in time the same amount time 
by which $dT$ was exceeded''. Note that, if $z\in S_0$ is close to
$\Stan_1$, then it may be that such $t_*$ does not exist
(the trajectory never hits the threshold manifold for $I=A$) and
hence one cannot extend $\s_1$ to $S_0$.

Regarding the contractiveness of $\s$, we first discuss the map
$\s_0$ given in Equation~\eqref{eq:strobo_s0}. Recalling that for $A=0$
system~\eqref{eq:general_system} possesses a unique attracting equilibrium
point in $\D$, this implies that $\varphi_0$ is contracting in
$\D$. For $A>0$ small, $\varphi_A$ is also contracting. By contrast, for
larger values of $A$, $\varphi_A$ may be expanding in $\D$. If this is
the case, this expansiveness can be compensated by integrating
$\varphi_0$ for large enough time, $T-dT$, which occurs if $dT$ is
smaller
enough than $T$ so that we obtain
\begin{equation*}
\left| \s_0(z)-\s_0(z') \right|<\left| z-z' \right|.
\end{equation*}
Regarding $\s_1$, the spike exhibited by trajectories of point in
$S_1$ may introduce expansiveness to $\s_1$. Arguing similarly, this
expansiveness can be compensated making $dT$ smaller enough than
$T$ so that the contracting flow $\varphi_0$ is applied for large enough
time.

\section{Border collisions of fixed points of the stroboscopic map}\label{sec:bif_of_fp}
\begin{figure}
\begin{center}
\begin{tabular}{ll}
(a) & (b) \\
\includegraphics[width=0.4\textwidth]{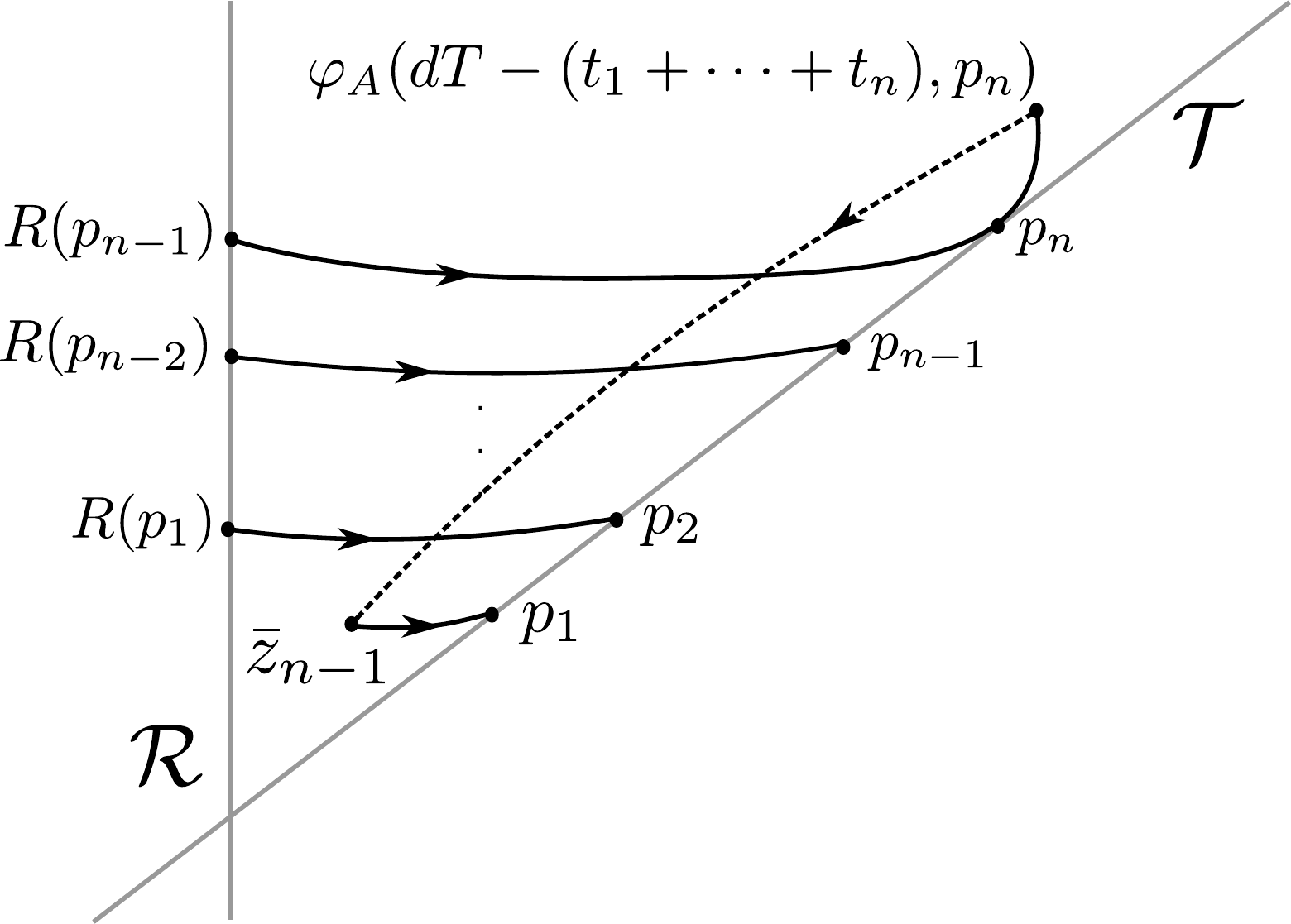}&
\includegraphics[width=0.4\textwidth]{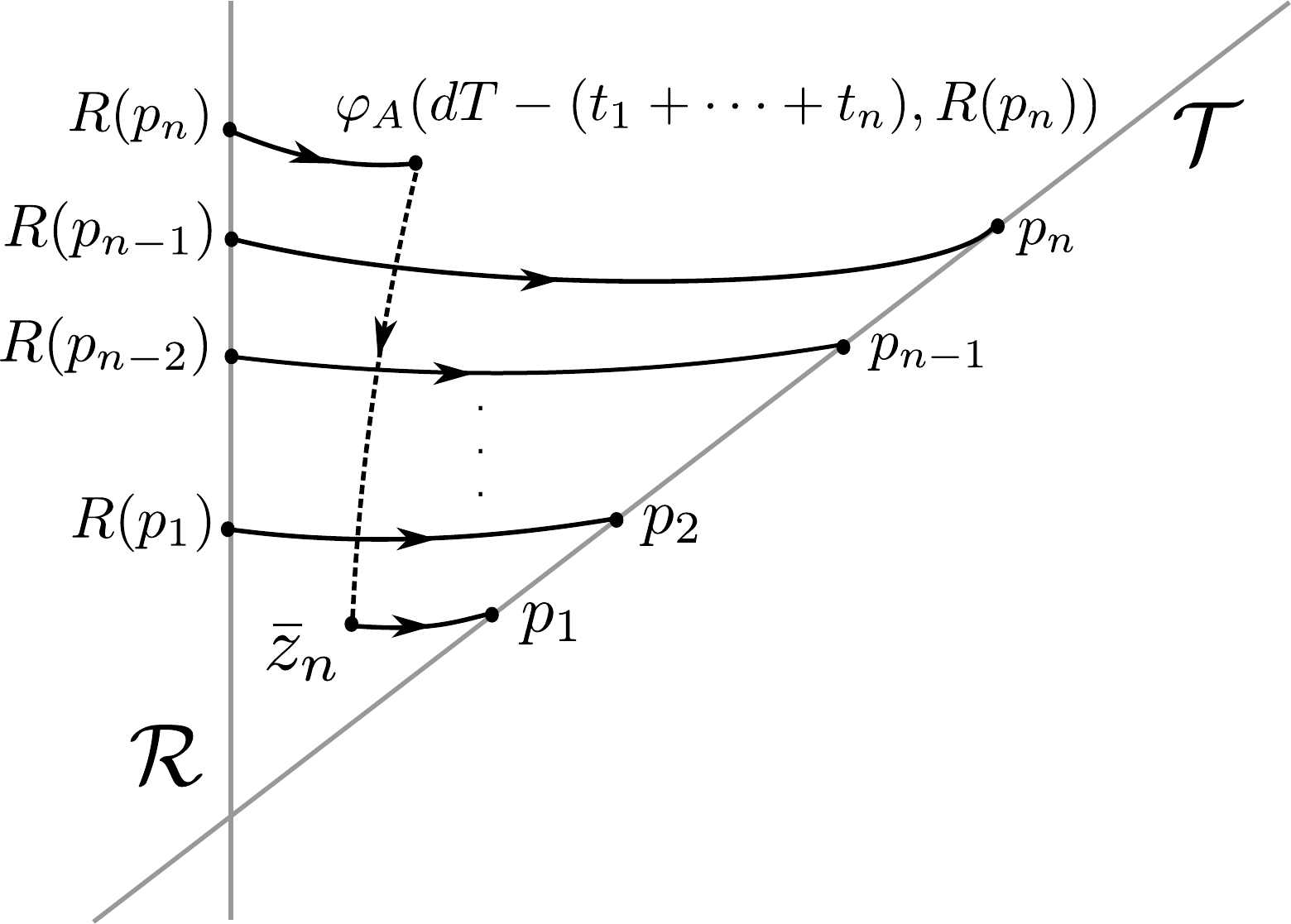}\\
(c) & (d) \\
\includegraphics[width=0.4\textwidth]{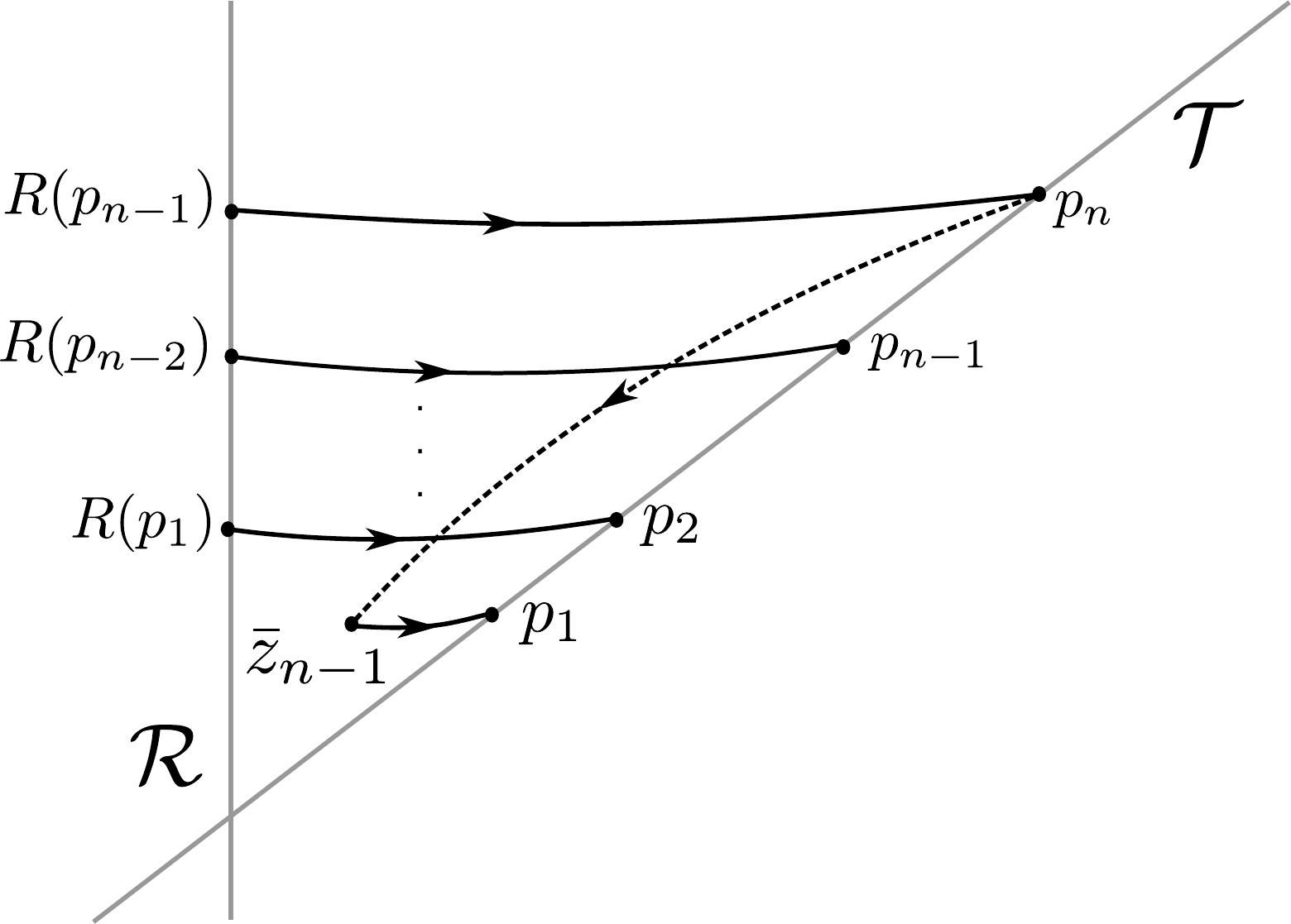}&
\includegraphics[width=0.4\textwidth]{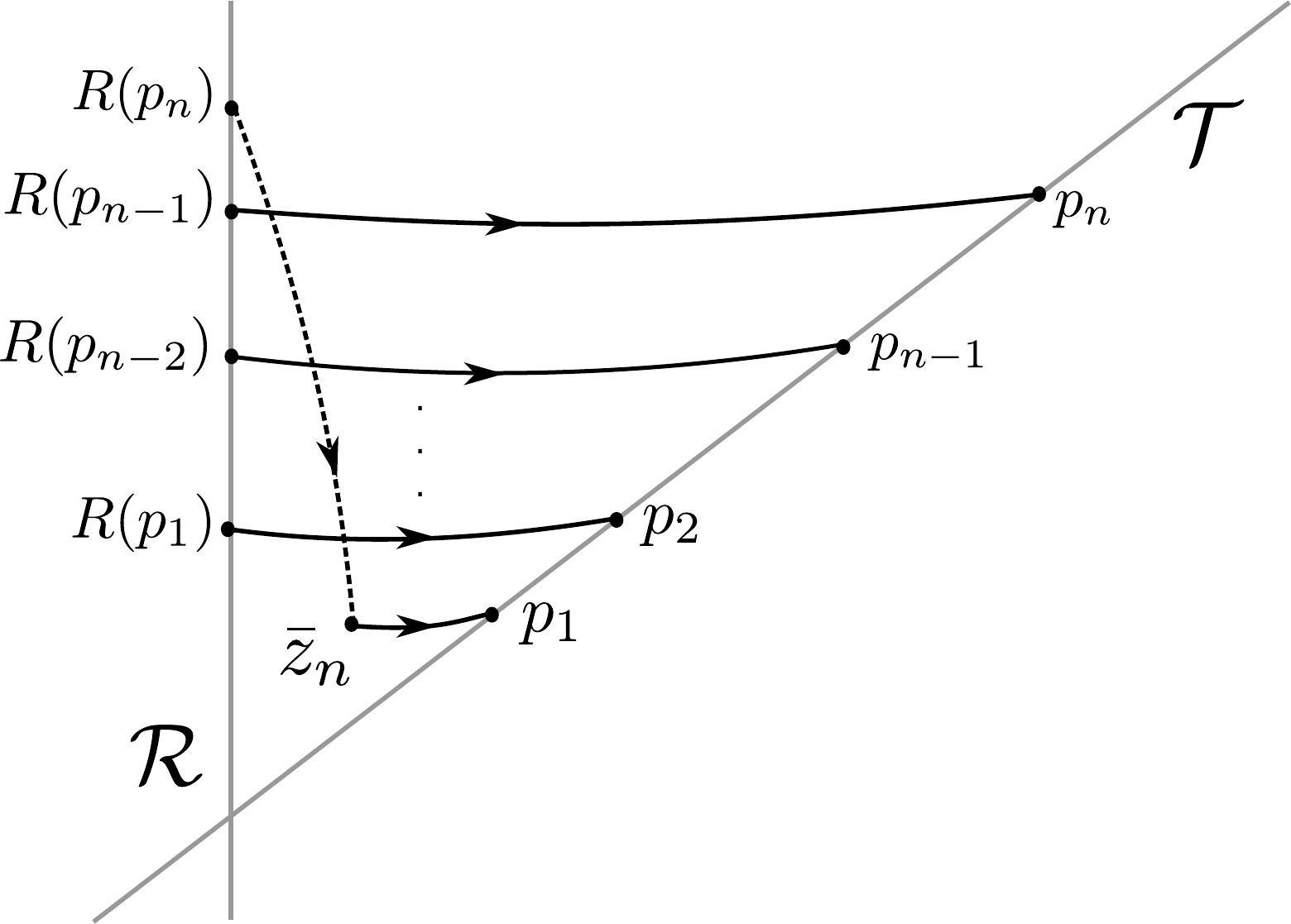}\\
\end{tabular}
\end{center}
\caption{Scheme showing a periodic orbit of the system~\systwreset{} (a,b) grazing tangentially the manifold $\mathcal{T}$ at the point $p_n$ and (c,d) 
grazing the manifold $\mathcal{T}$ at the non-differentiable point $p_n$ at time $dT$. These periodic orbits correspond to fixed points of the stroboscopic map hitting the 
manifold $\Stan_n$ (a,b) and the manifold $\Stran_n$ (c,d). Panels (a,c) correspond to a collision of a point $\bar{z}_{n-1} \in S_{n-1}$ with $\Sigma_n$, while panels (b,d)
correspond to a collision of the point $\bar{z}_{n} \in S_{n}$ with $\Sigma_n$. See Section~\ref{sec:bif_of_fp} for more details. Solid black lines correspond to trajectories of the flow $\varphi_A$, while dashed black lines 
correspond to trajectories of the flow $\varphi_0$. Grey lines correspond to the reset manifold $\mathcal{R}$ and the threshold manifold $\mathcal{T}$.}
\label{fig:scheme_bif}
\end{figure}

By assumption {\em H.1},  for $A=0$, system~\systwreset{} possesses an
attracting subthreshold equilibrium point, $z^*\in\D$.  Although the
periodic forcing $I(t)$ is not continuous in $t$, averaging
theorem~\cite{BogMit61} holds, as the system is Lipschitz in $z$. This
implies that, when $A>0$ is small enough, system~\systwreset{}
possesses a $T$-periodic orbit, which is not differentiable at
$t=0\pmod{T}$ and $t=dT\pmod{T}$ (recall that we assumed $t_0=0$).
As its amplitude increases with $A$, this periodic orbit may undergo a
{\em grazing bifurcation} \cite{Nordmark97,BerBudCha01a} if it
collides with the boundary $\TM$ when varying $A$. This
corresponds to a fixed point $\bz_0\in S_0$ of the stroboscopic map
undergoing a {\em border collision bifurcation} \cite{NusOttYor94} when
colliding with $\Sigma_1$. In general, a bifurcation occurs when a
fixed point $\bz_n \in S_n, n \geq 1$ collides with $\Sigma_{n+1}$ or
$\Sigma_{n}$.  Following {\em i)} and {\em ii)} of
Section~\ref{sec:strobo_def}, we distinguish two different types of
border collision bifurcations:
\begin{enumerate}[{Bif}.1]
\item \textbf{Smooth grazing bifurcation}: the $T$-periodic orbit of
system~\systwreset{} grazes tangentially $\TM$, and, equivalently, the
fixed point of the stroboscopic map collides with $\Stan_n$, for some $n \geq 1$.
The fixed point can collide with $\Stan_n$ from two different regions, namely, $S_{n-1}$ and $S_n$.
When the fixed point $\bar{z}_{n-1} \in S_{n-1}$
collides with $\Stan_n$, the following equations are satisfied (see
Figure~\ref{fig:scheme_bif}(a)):
\begin{equation}\label{eq:bc1tan}
\bar{z}_{n-1}=\varphi_A(-t_1,p_1) = \varphi_0(T-dT,\varphi_A(dT-(t_1+\cdots+t_n),p_n) 
\end{equation}
with
\begin{equation}\label{eq:bc1}
\begin{array}{rcl}
p_2 &=& \varphi_A (t_2,R(p_1)) \\
p_3 &=& \varphi_A (t_3,R(p_2)) \\
\vdots && \\
p_n &=& \varphi_A (t_n,R(p_{n-1})) 
\end{array}
\end{equation}
and
\begin{equation}\label{eq:bc11}
\begin{array}{c}
(f(p_n) + vA )\cdot \nabla h (p_n) =0 \\
h(p_1)=h(p_2)=\ldots=h(p_n)=0.
\end{array}
\end{equation}
While, when $\bar{z}_{n} \in S_n$ collides with $\Stan_n$ the equation is
\begin{equation}\label{eq:bc1tan2}
\bar{z}_{n}=\varphi_A(-t_1,p_1) = \varphi_0(T-dT,\varphi_A(dT-(t_1+\cdots+t_n),R(p_n)) 
\end{equation}
and the other conditions \eqref{eq:bc1}-\eqref{eq:bc11} as before (see
Figure~\ref{fig:scheme_bif}(b)).

\begin{remark}
Notice that we assume that the smooth grazing with the manifold $\TM$ occurs after the 
last spike (this will be the situation in the example considered in Section \ref{sec:numerical_results}). 
In general, this is not necessary the case, and similar equations can be written for other situations.
\end{remark}

\item \textbf{Nonsmooth grazing bifurcation}: the $T$-periodic orbit of
system~\systwreset{} grazes $\TM$ at the non-differentiable point
given by $t=dT$,  and, equivalently, the fixed point of the
stroboscopic map collides with $\Stran_n$, for some $n \geq 1$. The fixed 
point can collide with $\Stran_n$ from two different regions, namely, $S_{n-1}$ and $S_n$.
When the fixed point $\bar{z}_{n-1} \in S_{n-1}$
collides with $\Stran_n$, the following equations are satisfied (see
Figure~\ref{fig:scheme_bif}(c)):
\begin{equation}\label{eq:bc1tran}
\bar{z}_{n-1}=\varphi_A(-t_1,p_1) = \varphi_0(T-dT,p_n) 
\end{equation}
with
\begin{equation}\label{eq:bc2}
\begin{array}{rcl}
p_2 &=& \varphi_A (t_2,R(p_1)) \\
p_3 &=& \varphi_A (t_3,R(p_2)) \\
\vdots && \\
p_n &=& \varphi_A (dT-(t_1+ \cdots + t_{n-1}),R(p_{n-1})) 
\end{array}
\end{equation}
and 
\begin{equation}\label{eq:bc22}
h(p_1)=h(p_2)=\ldots=h(p_n)=0.
\end{equation}
While, when $\bar{z}_{n} \in S_n$ collides with $\Stran_n$ the equation is
\begin{equation}
\bar{z}_{n}=\varphi_A(-t_1,p_1) = \varphi_0(T-dT,R(p_n)) 
\label{eq:bc1tran2}
\end{equation}
and the other conditions \eqref{eq:bc2}-\eqref{eq:bc22} as before (see
Figure~\ref{fig:scheme_bif}(d)).

\end{enumerate}

In Section~\ref{sec:numerical_results} we compute (for a particular example) the critical 
parameter values at which the fixed points $\bar{z}_0$ and $\bar{z}_1$ undergo 
border collision bifurcations when colliding with $\Sigma_1$ and $\Sigma_2$ 
(the latter, only for $\bar{z}_1$) by means of solving numerically (using a Newton method) 
the systems of equations given above. 

After a border collision bifurcation of a fixed point in $S_n$ colliding with $\Stran_{n+1}$, 
it is expected that the map will map points in $S_n$ to points in $S_{n+1}$ and viceversa, thus causing the dynamics to 
alternate between $S_n$ and $S_{n+1}$. Therefore, it is possible that there appear 
periodic orbits of the stroboscopic map hitting both regions $S_{n}$ and $S_{n+1}$. In the next section, we 
will provide techniques to study such periodic orbits.

\section{Periodic orbits of the stroboscopic map}\label{sec:periodic_orbits}

Beyond the fixed points, we also study periodic orbits of the stroboscopic
map $\s$. Assume that the fixed point $\bar{z}_n$ of $\s$ collides with $\Stran_{n+1}$
for some parameter value undergoing a border collision bifurcation as described in 
Section~\ref{sec:bif_of_fp}. In this situation, $\s$ may possess periodic orbits 
visiting both $S_{n}$ and $S_{n+1}$. Unfortunately, there is 
little general theory that can be applied to prove the existence of such periodic orbits.
In this section, we review possibly the only result (to our knowledge)
in this direction by Gambaudo \emph{et al.}
We first introduce symbolic dynamics and some
definitions in order to characterize these periodic orbits. 

\begin{definition}\label{def:symb_iti}
Given $z\in S_{n}\cup S_{n+1}$, we define the itinerary of $z$ by $\s$ as
\begin{equation*}
I_\s(z)=\left( a(z),a\left( \s(z) \right ),a\left( \s^2(z) \right),\dots
\right),
\end{equation*}
where
\begin{equation*}
a(z)=
\left\{
\begin{aligned}
&1&&\text{if }z \in S_{n+1}\\
&0&&\text{if }z\in S_{n}.
\end{aligned}
\right.
\end{equation*}
\end{definition}

\begin{remark} 
Although we consider only periodic orbits that interact with two regions ($S_n$ and $S_{n+1}$), it is 
possible to extend the results and definitions to orbits interacting with more than two regions. However, this 
situation is out of scope of our paper.
\end{remark}

\begin{definition}\label{def:Wpq}
One calls $W_{p,q}$ the set of periodic symbolic
sequences generated by infinite concatenation of a symbolic
block of length $q$ containing $p$ symbols $1$:
\begin{equation*}
W_{p,q}=\left\{\y\in\left\{ 0,1
\right\}^\mathbb{N}\,|\,\y= \x^\infty,\,\x\in\left\{ 0,1 \right\}^q\,\text{and $\x$ contains
$p$ symbols $1$} \right\}.
\end{equation*}
\end{definition}

\begin{definition}\label{def:rotn}
One says that a sequence in $W_{p,q}$ has rotation number $p/q$.
\end{definition}
\begin{remark}\label{rem:rotn}
In the one-dimensional case, this definition of the rotation
number coincides with the classical one for one-dimensional circle
maps through a lift (see~\cite{GraAlsKru17}). However, in the planar
case, it becomes in general difficult (if possible) to define this
number by means of lifts, as the dynamics cannot always be reduced to
a 2-dimensional torus. However, following~\cite{GamTre88,GamGleTre88},
we abuse notation and call this number the ``rotation number''.
\end{remark}
\begin{definition}\label{def:order}
Symbolic sequences can be ordered using that $0<1$.  Hence,
\begin{equation*}
(\x_1\x_2\dots)< (\y_1\y_2\dots)
\end{equation*}
if and only if $\x_1=0$ and $\y_1=1$ or $\x_1=\y_1$ and there exists
some $j> 1$ such that
\begin{align*}
\x_i&=\y_i,\,\text{for all}\;i< j\\
\x_j&=0\\
\y_j&=1.
\end{align*}
\end{definition}
This order allows one to consider the following definition: 
\begin{definition}\label{def:maximin}
Let $\sigma$ be the shift operator. One says that a symbolic sequence $\x\in W_{p,q}$ is maximin if
\begin{equation*}
\min_{0\le k\le q}\left( \sigma^k(\x) \right)=\max_{\y\in
W_{p,q}}\left(
\min_{0\le k\le q}\left( \sigma^k(\y) \right)
\right).
\end{equation*}
\end{definition}

\begin{example}
Up to cyclic permutations, there exist only two periodic sequences in
$W_{2,5}$, which are represented by means of two blocks that,
when expressed in minimal form, are given by $0^3 1^2$
and $0^2101$. The maximum of the minimal blocks is $0^2 101$,
therefore the symbolic sequence generated by $0^2 101$ is maximin.
\end{example}

Intuitively, maximin symbolic sequences have ``well'' distributed
symbols $1$ along the sequence, which is related to the notion of
``well ordered'' symbolic sequences (see Definition~\ref{def:well-ordered_symbolic-sequence}).

Alternatively, maximin itineraries can be defined as those belonging
to the Farey tree of symbolic sequences. This means that they are
given by concatenation of two maximin sequences such that their
rotation numbers are Farey neighbours. See~\cite{GraAlsKru17} for a
recent review in this topic.\\
Then, we may apply the following result to study the existence of
maximin periodic orbits:
\begin{theorem}\label{teo:gambaudo}[Dynamics of quasi-contractions]
Assume that there exist sets $E_0\subset S_{n}$ and $E_1\subset
S_{n+1}$ such that
\begin{enumerate}[i)]
\item $\s(E_i)\subset E_0\cup E_1$, for $i=0,1$.
\item $\s_0$ and $\s_1$ contract in $E_0$ and $E_1$, respectively.
\item $\s^i(\Stran)\cap \Stran = \emptyset$ for all $i\ge1$.
\end{enumerate}
Then, provided that $\s$ preserves orientation, $\s$ possesses $0$ or $1$ 
periodic orbit. In the latter case,  its itinerary is maximin.
\label{theo:quasi-contraction}
\end{theorem}
The previous result was stated in~\cite{GamTre88} for
quasi-contractions in metric spaces and adapted in~\cite{GraAlsKru17}
for piecewise continuous contracting maps in $\RR^n$.

In Section \ref{sec:perorbits} we will show how Theorem
\ref{teo:gambaudo} can be applied to a particular example to prove the
existence of a periodic orbit of maximin type by checking the
hypothesis using semi-rigorous numerics.

\section{Application to a neuron model}\label{sec:numerical_results}
In this section, we apply the theoretical results presented in previous sections to a spiking
neuron model of integrate-and-fire type with a dynamic threshold.
\subsection{The model}
We consider the system proposed in~\cite{MenHugRin12}, which consists of 
a leaky inte\-grate-and-fire model with a dynamic threshold.
It is a dimensionless version of other similar models such as \cite{HiggsS11, Platkiewicz10}.
The sytem is submitted to periodic forcing $I(t)$ as in Equation~\eqref{eq:forcing}.
The equations are given by:
\begin{equation} 
\begin{aligned}
\dot{V}&=-V+V_{0}+I(t)\\
\tau_\theta\dot{\theta}&=- \theta + \theta_\infty(V) 
\end{aligned}
\label{eq:system_visual}
\end{equation}
where $(V,\theta)\in \RR^2$ are the neuron voltage and the threshold, respectively.
The function
\begin{equation}
\theta_\infty(V)=a+e^{b\left( V-c \right)}
\label{eq:theta_infty}
\end{equation}
is the steady state value of the threshold $\theta$, with $a,b,c \in \mathbb{R}$;
$\tau_{\theta}$ is the time constant for the threshold (which will be chosen only a
bit slower than the membrane time constant, i.e. $\tau_{\theta}> 1$) and $V_0$ is the voltage at the
resting state.
The spiking reset rule is given by:
\begin{equation}\label{eq:reset_ex}
\text{if }V(t_*)=\theta(t_*)\text{ then }V(t_*^+)=V_r\text{ and
}\theta(t_*^+)=\theta(t_*)+\Delta,
\end{equation}
with $V_r$ and $\Delta$ being real parameters. 
The parameters of the system along this paper are $V_0 = 0.1$, $V_r = 0$, $\Delta = 0.3$, $a = 0.08$,
$c = 0.53$ and $\tau_{\theta} = 2$. Parameter $b$ will vary along this study between 0 and 1.

Notice that system~\eqref{eq:system_visual}-\eqref{eq:reset_ex} is of the form~\systwreset{}, with a threshold manifold
\begin{equation*}
\TM=\left\{ (V,\theta) \in \RR^2 \,|\,h(V,\theta)=V-\theta=0 \right\},
\end{equation*}
reset map
\begin{equation*}
\begin{array}{cccc}
\Rmap:&\TM&\longrightarrow&\mathcal{R}\\
&(V,\theta)&\longmapsto&(V_r,\theta+\Delta),
\end{array}
\end{equation*}
reset manifold
\begin{equation*}
\mathcal{R}=\left\{ (V,\theta) \in \RR^2 \,|\, V=V_r \right\},
\end{equation*}
and subhtreshold domain
\begin{equation*}
\D=\left\{ (V,\theta),\,|\,V\ge V_r,\,V < \theta \right\}.
\end{equation*}
Notice that for biological reasons we restrict the subthreshold domain to $V\geq V_r$.

Next we show that system~\eqref{eq:system_visual} satisfies
hypothesis \emph{H.1} and \emph{H.2.}  Indeed, for $A=0$, the system has an equilibrium
point at $(V^*,\theta^*)=(V_0,\theta_\infty(V_0))$ (which satisfies $V_0<\theta_\infty(V_0)$ for the 
choice of parameters) with eigenvalues 
$\lambda_1=-1$ and $\lambda_2=-1/\tau_{\theta}$. Thus, it is an attracting node. The associated
eigenvectors are $v_1=(1, \theta'_{\infty}(V^{*})/(1-\tau_\theta))$ and $v_2=(0,1)$, respectively. Then, assuming
that $\tau_{\theta}>1$, trajectories approach the equilibrium point
tangentially to the $V$-nullcline ($V=V_0$).
Thus, given the geometry of the domain $\D$, if there exist points
whose orbits intersect the threshold manifold $\TM$, then there exists
an orbit in $\D$ which is tangent to $\TM$, and separates the orbits
of those points that intersect $\TM$ and those that do not.  Imposing
that $\nabla h$ must be perpendicular to the vector field on the
manifold $\TM$ at the tangent point $(\bar{V},\bar{V})$ we have that \[
\theta_{\infty}(\bar{V})=V_0.\] 
Choosing parameters $a,b,c$ such that $\theta_{\infty}(V) > V_0$ for all values of $V$ 
(i.e. $\theta_{\infty}(V)>a+e^{-b(V_r - c)}>V_0$), we have that all points in $\D$ belong to orbits
that do not intersect $\TM$ for $A=0$.

\begin{remark}
If one choses a function $\theta_{\infty}(V)$ for which hypothesis H.2
is not satisfied, the mathematical analysis presented in
Sections~\ref{sec:general_setting}-\ref{sec:bif_of_fp} still follows
if $z^*=(V^{*},\theta^{*})$ is enough isolated from those points not satisfying hypothesis \emph{H.2}. In
this case, one can safely remove these points from $\D$ and the
analysis in the mentioned sections holds nevertheless. 
\end{remark}

\subsection{Fixed points and their bifurcations}\label{sec:fixed_points_numeric}
In this section we analyze the bifurcations of the fixed points of
the stroboscopic map (corresponding to $T$-periodic orbits of
the time-continuous system~\eqref{eq:system_visual}-\eqref{eq:reset_ex}) when varying parameters of the
system. We focus on bifurcations exhibited by the fixed points
$\bz_0\in S_0$ and $\bz_1\in S_1$, as they are more relevant from an
applied point of view, given that they combine spiking and
subthreshold dynamics. A similar analysis can be done for fixed points
exhibiting more spikes, $\bz_n\in S_n$. We focus on bifurcations 
associated to piecewise smooth maps (border collisions), although other bifurcations of smooth maps such as the
Neimark-Sacker bifurcation may occur. As explained in Section~\ref{sec:bif_of_fp},
such border collision bifurcations correspond to a periodic orbit
of the time-continuous system grazing the threshold, which can occur through a tangency (smooth
grazing or, equivalently, border collision with $\Stan_1$ or
$\Stan_2$) or when disabling the pulse (nonsmooth grazing or,
equivalently, border collision with $\Stran_1$ or $\Stran_2$).

\begin{figure}
\begin{center}
\begin{tabular}{cc}
(a) $T=5$, $b=0.1$ and $A=2$& \\
\includegraphics[width=0.36\textwidth]{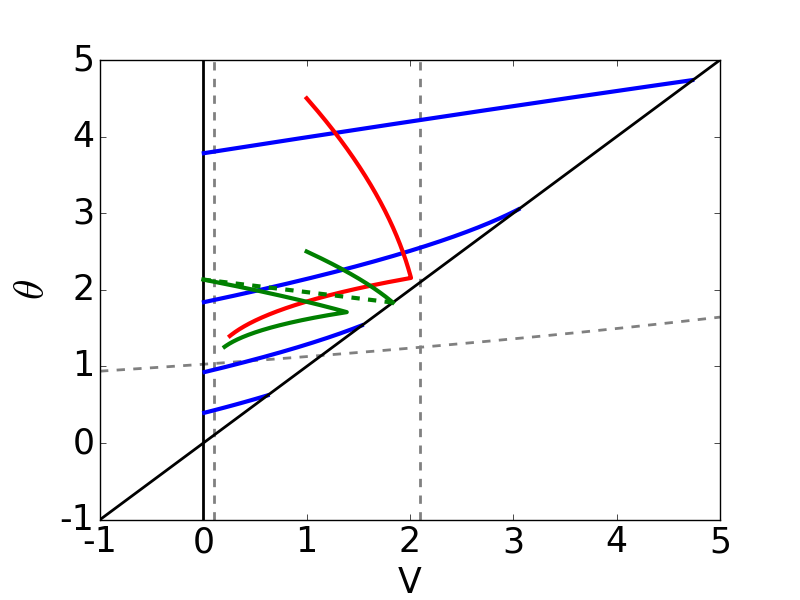}&
\includegraphics[width=0.36\textwidth]{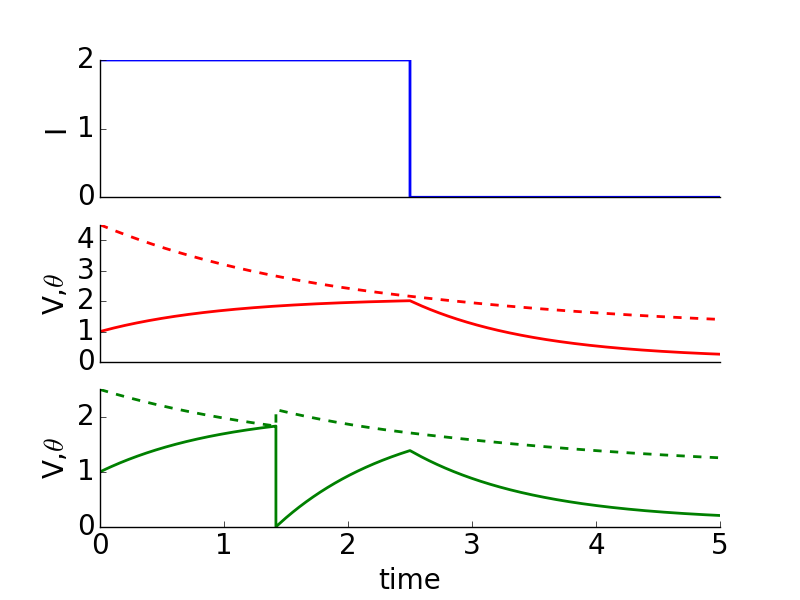}\\
(b) $T=0.5$, $b=0.1$ and $A=2$ & \\
\includegraphics[width=0.36\textwidth]{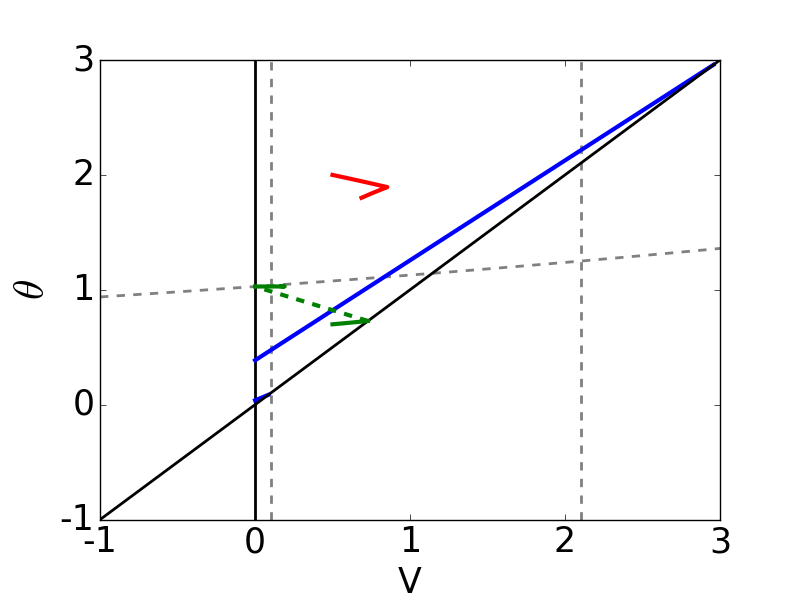}&
\includegraphics[width=0.36\textwidth]{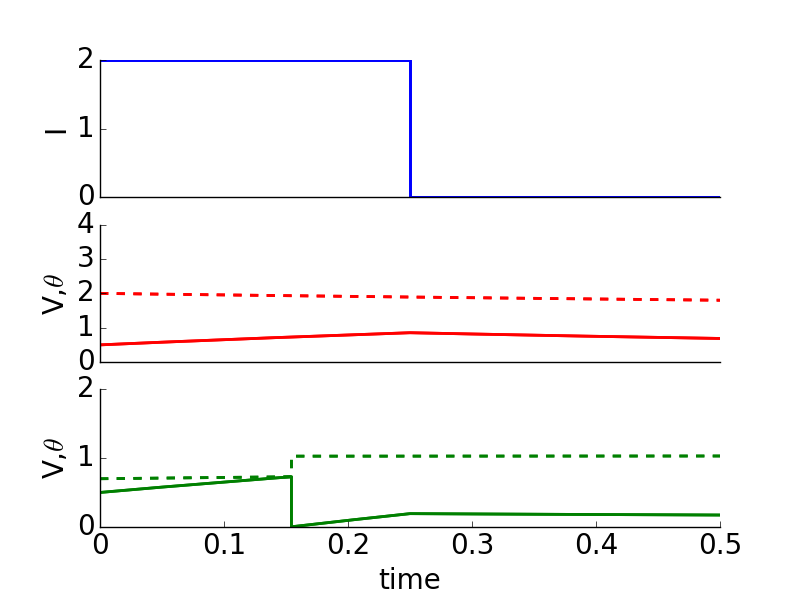}\\
(c) $T=5$, $b=0.5$ and $A=2$ & \\
\includegraphics[width=0.36\textwidth]{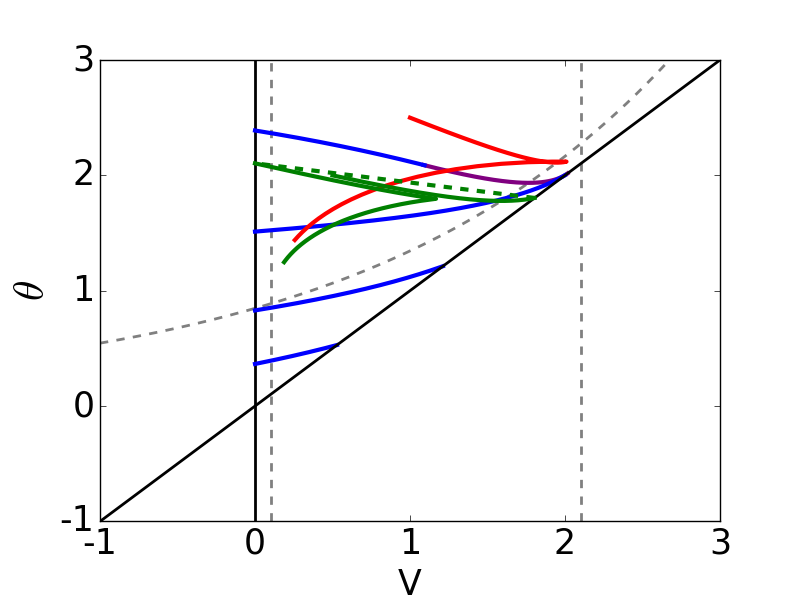}&
\includegraphics[width=0.36\textwidth]{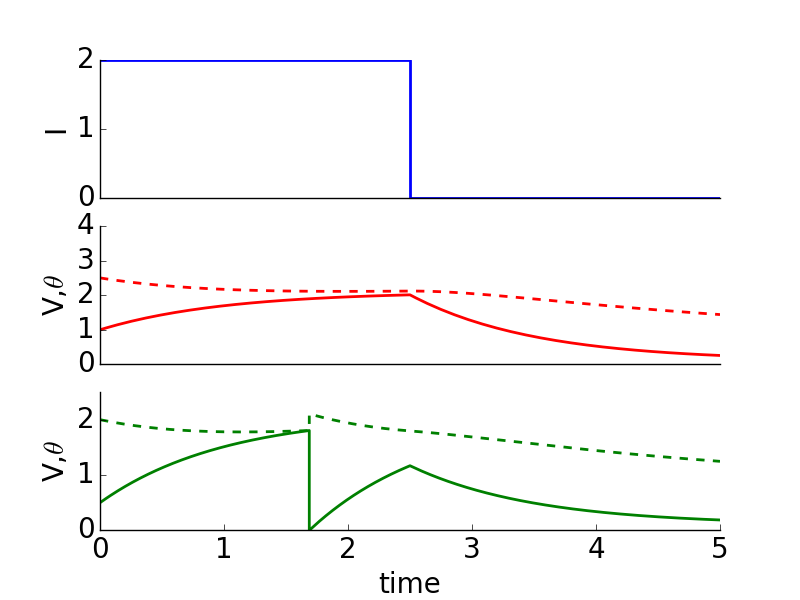} \\
(d) $T=0.5$, $b=0.5$ and $A=2$ & \\
\includegraphics[width=0.36\textwidth]{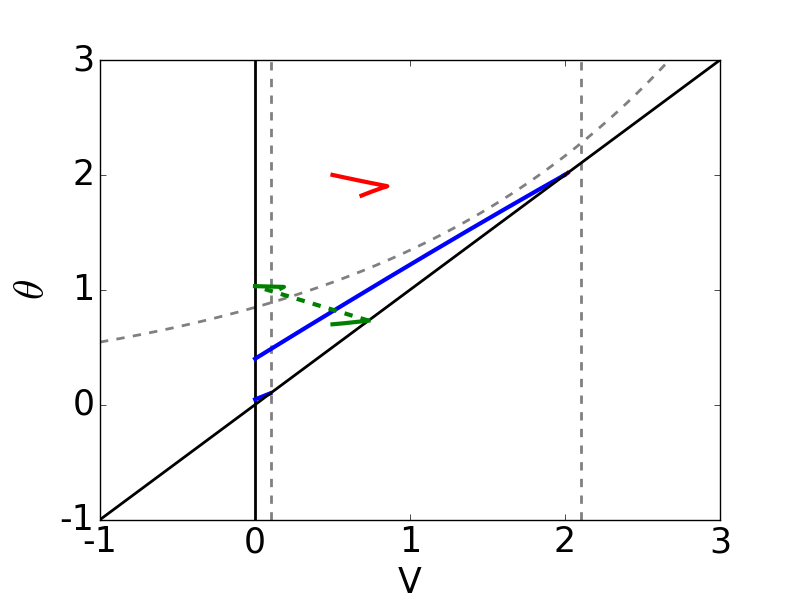}&
\includegraphics[width=0.36\textwidth]{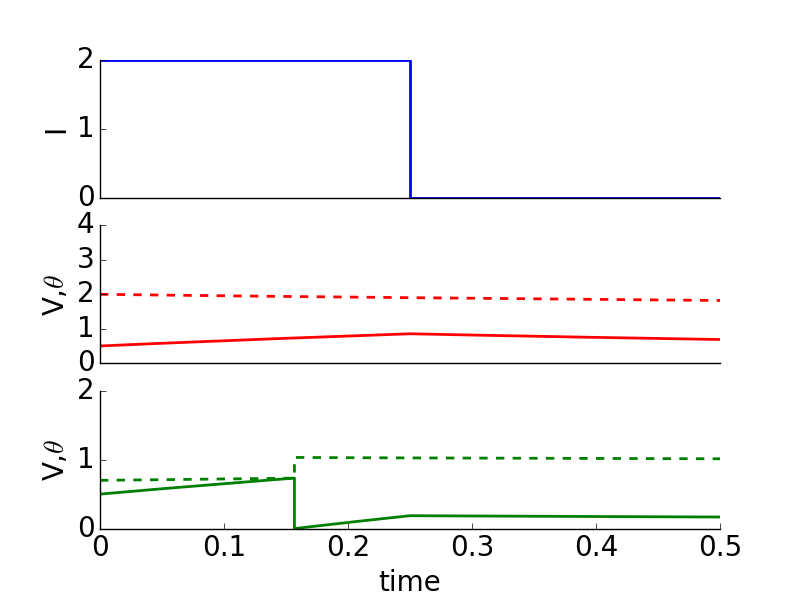} \\
\end{tabular}
\end{center}
\caption{(Left) Trajectories of
system~\eqref{eq:system_visual}-\eqref{eq:reset_ex} on the
$(V,\theta)$ phase space with initial conditions in sets $S_0$ (red)
and $S_1$ (green).  Blue curves correspond to the boundaries
$\Stran_n$, $n\geq1$ and the purple curve to the boundary $\Stan_1$.
Gray dashed curves correspond to $\theta$- and $V$-nullclines for
$I=0$ and $I=A$; and their intersection corresponds to the point $(V^*,\theta^*)$. Solid black curves correspond to the reset manifold
$\RM$ and the threshold manifold $\TM$. (Right) Times courses of
$I(t)$ and the variables $V$ (solid) and $\theta$ (dashed) for the
trajectories shown on the left with the same color.}
\label{fig:Sigmatan-Sigmatrans}
\end{figure}

One of the characteristics of system~\eqref{eq:system_visual} which
may influence having smooth or nonsmooth grazing bifurcations is the position
of the equilibrium point for $I=A$ (given by the intersection of the
nullclines, which occurs at $V^*=V_0+A$ and
$\theta^*=\theta_\infty(V_0+A)$). If this point happens to be far away
from the domain $\D$ ($V\gg \theta$), then the dynamics is fast,
trajectories spend little time between spikes and transversal grazing
is most likely to occur.
However, if the equilibrium point is close to the threshold manifold
(or even at $\D$), then the dynamics is slower and system may exhibit
tangencies with $\TM$. The nature of the function
$\theta_\infty(V)$ allows these two situations mainly by varying the
parameter $b$ between 0 and 1. For small $b$ (close to $0$) the nullcline $\theta=\theta_\infty(V)$
becomes almost flat in $\D$ and fixes the equilibrium point for $I=A$
outside $\D$ (see Figures~\ref{fig:Sigmatan-Sigmatrans} (a) and (b)).
However, for larger values of $b$ the function $\theta_\infty(V)$ may
be completely located in $\D$ (see
Figures~\ref{fig:Sigmatan-Sigmatrans} (c) and (d)).  Moreover, the latter case
has consequences from the neuron modeling point of view as these
systems are not capable to show repetitive firing for constant input
and they are referred as phasic neurons \cite{HugMenRin17}.

Apart from $A$ and $b$, other relevant parameters influencing these
different type of behaviours are $T$ and $d$, as they control the
integration time during the active part of the pulse. When the pulse is active for 
a short time, the regions $S_n$, $n>1$, occupy a small portion of the 
subthreshold domain (see Figures~\ref{fig:Sigmatan-Sigmatrans}(b) and (d)). 
In this work we keep $d=0.5$ fixed, and study bifurcations of the fixed points $\bz_0$
and $\bz_1$ when varying $b$, $A$ and $T$.

\begin{figure}
\begin{center}
\includegraphics[angle=-90,width=0.9\textwidth]{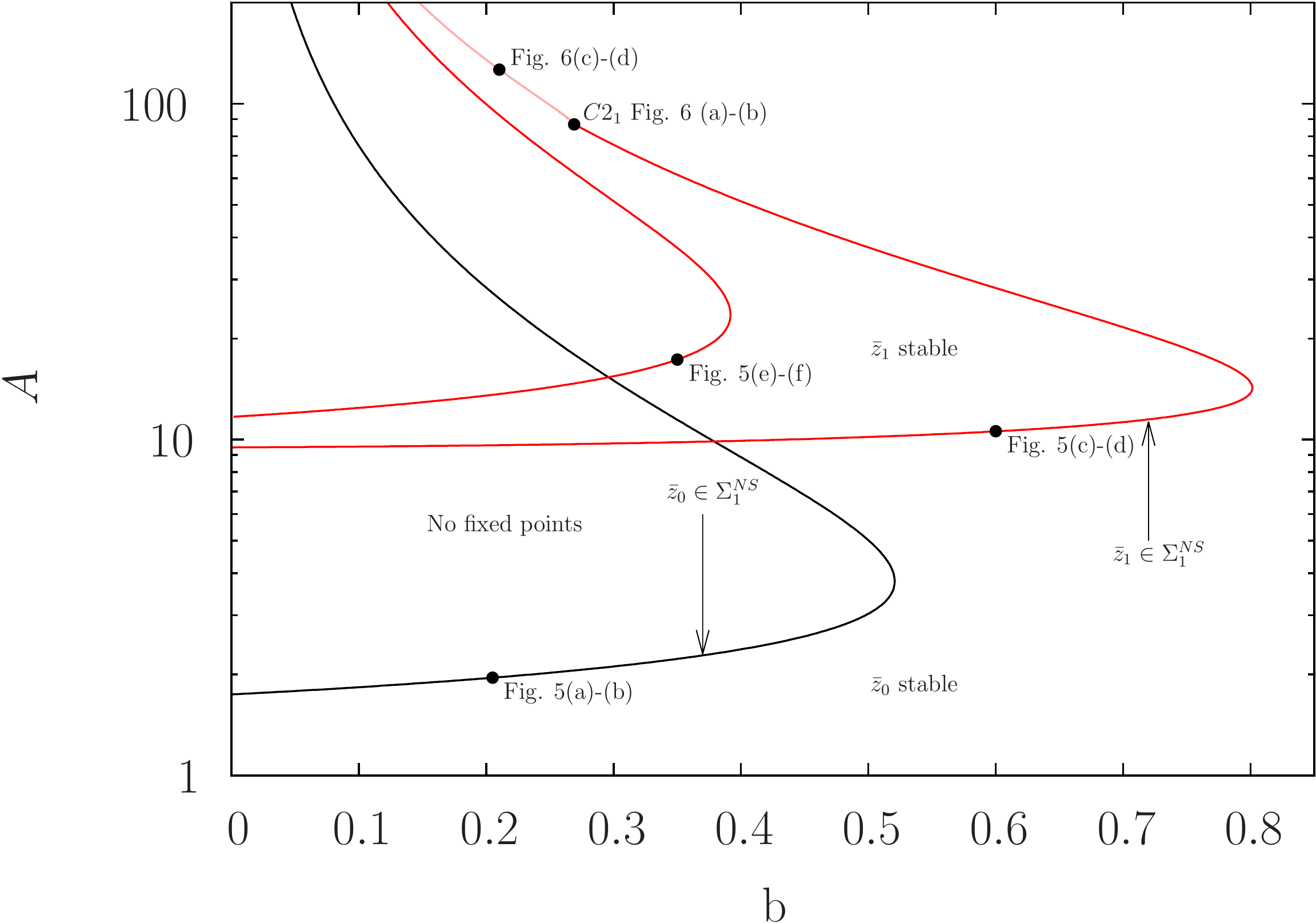}
\end{center}
\caption{Border collision bifurcation curves and regions of existence of $\bz_0$ and $\bz_1$ in the $(b,A)$
parameter space for $T=0.5$ (vertical axis is in logarithmic scale).
Curve in black corrsponds to a border collision of $\bz_0\in S_0$ with $\Stran_1$. Curve in dark red corresponds to a
border collision of $\bz_1\in S_1$ with $\Stran_1$ (outer curve) and $\Stran_2$ (inner curve).
Curve in light red corresponds to a border collision of $\bz_1\in S_1$ with $\Stan_1$. Dots on these curves indicate the
parameter values $(A,b)$ for which we show (in the Figure indicated nearby) 
the trajectory of the periodic orbit of the time continuous system that undergoes a grazing bifurcation.
}
\label{fig:2dbif_T0d5}
\end{figure}

We first fix $T=0.5$ and compute the bifurcation curves of the fixed
points $\bz_0$ and $\bz_1$ in the parameter space $(b,A)$ (see
Figure~\ref{fig:2dbif_T0d5}).  These curves have been computed
semi-analytically using a predictor-corrector method
detailed in Appendix~\ref{sec:compu_bif-curves}. As
we are computing periodic orbits close to bifurcations, the method may
predict or correct a point outside the feasible domain. However, as
explained in Section~\ref{sec:extension}, the system can be extended
to virtual domains and allow the Newton method to continue and
converge.\\
In black we show the border collision curve given by the collision of 
the fixed point $\bz_0\in S_0$ with $\Stran_1$ (nonsmooth grazing). Recall
that at this bifurcation a non-spiking $T$-periodic orbit
grazes the threshold precisely when the pulse is disabled (see
Figures~\ref{fig:Sigma1trans} and~\ref{fig:Sigma1trans_TS}). Hence, as
detailed in Appendix~\ref{sec:compu_bif-curves}, this curve has been
computed numerically solving Equations \biftranSn{} for $n=1$
which, in this particular case, become
\begin{equation}
\varphi_A(-dT;(V,V))=\varphi_0(T-dT;V,V),
\label{eq:Sigma1_trans}
\end{equation}
where $V$, $b$ and $A$ are the unknowns. Crossing this curve with
increasing $A$, the fixed point $\bz_0$ first disappears and reappears
again, as the curve exhibits a fold.  For the given value of $T$,
$\bz_0$ does not exhibit any border collision involving $\Stan_1$.\\
In red we show two border collision curves. The outer one corresponds
to the collision of the fixed point $\bz_1\in S_1$ with $\Stran_1$ (nonsmooth grazing), which 
corresponds to a spiking $T$-periodic orbit 
grazing the threshold precisely when the pulse is disabled
(see Figures~\ref{fig:Sigma1transz1} and~\ref{fig:Sigma1transz1_TS}).
This curve
has been computed numerically solving Equations \biftranSnn{} for
$n=1$ which, in this particular case, become
\begin{equation}
\varphi_A(-dT;(V,V))=\varphi_0(T-dT;V_r,V+\Delta)
\label{eq:Simga1_trans_z1},
\end{equation}
where $V$, $b$ and $A$ are the unknowns.\\
The outer red bifurcation curve stops at a point labeled as $C2_1$. At
this point, the grazing bifurcation occurring at $t=dT$ becomes smooth
(see Figures~\ref{fig:C21} and~\ref{fig:C21_TS}). This point is hence
a co-dimension-two bifurcation point as both smooth and nonsmooth
grazing bifurcation conditions
(Equations~\biftanSnn{} and~\biftranSnn{} for $n=1$) are simultaneously satisfied. From
this point on, the fixed point $\bz_1$ collides with $\Stan_1$, which
corresponds  to the light red bifurcation curve
in Figure~\ref{fig:2dbif_T0d5}. Recall that 
at this bifurcation a spiking $T$-periodic orbit tangentially grazes the
threshold $\TM$ (see Figures~\ref{fig:Sigma1tansz1}
and~\ref{fig:Sigma1tansz1_TS}). This curve has been computed
numerically solving Equations \biftanSnn{} for $n=1$ which, in
this particular case, become
\begin{align}
\varphi_A(-t_1;(V,V))&=\varphi_0(T-dT;\varphi_A(dT-t_1;V_r,V+\Delta)\label{eq:Simga1_tan_z1_1}\\
-V+V_0+A&=\frac{-V+\theta_\infty(V)}{\tau_\theta}
\label{eq:Simga1_tan_z1_2},
\end{align}
where $V$, $t_1$, $b$ and $A$ are the unknowns.\\
The inner red curve corresponds to the collision of the fixed point
$\bz_1$ with $\Stran_2$, which corresponds to a spiking $T$-periodic
orbit which attempts to exhibit a new spike by grazing the
threshold precisely when the pulse is disabled (see
Figures~\ref{fig:Sigma2transz1} and~\ref{fig:Sigma2transz1_TS}). 
As detailed in Appendix~\ref{sec:compu_bif-curves}, this curve has been
computed numerically solving
Equations \biftranSn{} for $n=2$
which, in this particular case, become
\begin{equation}
\begin{aligned}
\varphi_A(-t_1;(V_1,V_1))&=\varphi_0(T-dT;V_2,V_2)\\
(V_2,V_2)^T&=\varphi_A(dT-t_1;V_r,V_1+\Delta)
\label{eq:z1_Sigma2trans},
\end{aligned}
\end{equation}
where $V_1$, $V_2$, $t_1$, $b$ and $A$ are the unknowns.\\
For the given value of $T$, $\bz_1$ does not exhibit any border collision involving $\Stan_2$.\\
In the region limited by the inner and outer red curves defining the
bifurcations $\bz_1\in \Sigma_2$ and $\bz_1\in \Sigma_1$,
respectively, the fixed point $\bz_1\in S_1$ exists.
Note also that the curves defined by $\bz_0\in\Stran_1$ and
$\bz_1\in\Stran_1$ (black and outer dark red curves) cross transversally.
This implies the existence of a region where both fixed points $\bz_0$
and $\bz_1$ coexist and are stable, as well as the existence
of a region where, none of the fixed points $\bz_0$ and $\bz_1$ exist.
Instead, one finds higher periodic orbits organized by period
adding-like structures, which will be treated in more detail in
Section~\ref{sec:perorbits}.\\

\begin{figure}
\begin{center}
\includegraphics[angle=-90,width=0.9\textwidth]{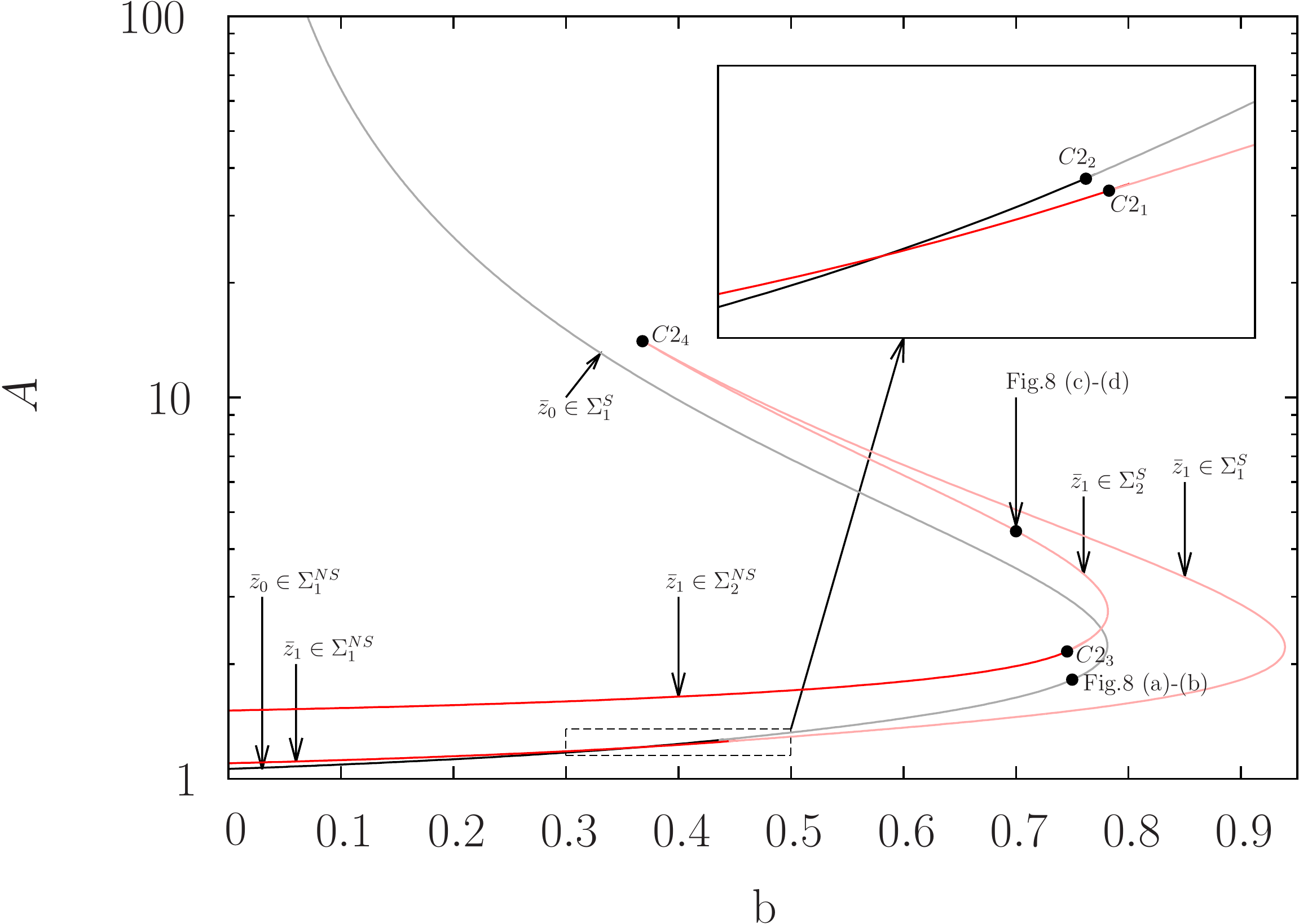}
\end{center}
\caption{Border collision bifurcation curves and regions of existence of the fixed
points $\bz_0$ and $\bz_1$ in the $(b,A)$ parameter space for $T=5$ (vertical axis in logarithmic scale).
Curves in black and gray correspond to border collision bifurcations of $\bz_0\in S_0$ with $\Stran_1$
and $\Stan_1$, respectively. Curves in red and light red correspond to border collision bifurcations of $\bz_1\in S_1$ with $\Stran$ and $\Stan$,
respectively. Co-dimension-two bifurcation points are labeled as $C2_i$, $i=2,3,4$ and are explained in the text.
Dots on these curves indicate the parameter values $(A,b)$ for which we show (in the Figure indicated nearby) 
the trajectory of the periodic orbit of the time continuous system that undergoes a grazing bifurcation.
}
\label{fig:2dbif_T5}
\end{figure}

In Figure~\ref{fig:2dbif_T5} we show the results of a similar
analysis for $T=5$. We observe the same nonsmooth grazing bifurcations as for $T=0.5$ but 
this case shows more smooth grazing bifurcations. Thus, we do not repeat the details for 
the nonsmooth grazing bifurcations that have already been discussed and we focus on the new ones.\\
We observe that for $T=5$ the fixed point $\bz_0\in S_0$ undergoes border collision
bifurcations through smooth grazing, $\bz_0\in\Stan_1$. That is, a
non-spiking $T$-periodic orbit tangencially grazes the threshold $\TM$ (see Figures~\ref{fig:z0Stan1}
and~\ref{fig:z0Stan1_TS}).
The corresponding bifurcation curve is shown in gray in
Figure~\ref{fig:2dbif_T5}, and it has been computed solving
Equations \biftanSn{}
for $n=1$ which, in this particular case,
become
\begin{align*}
\varphi_A(-t_1;(V,V))&=\varphi_0(T-dT;\varphi_A(dT-t_1;(V,V)))
\\
-V+V_0+A&=\frac{-\theta+\theta_\infty(V)}{\tau_\theta},
\end{align*}
where $V$, $t_1$, $b$ and $A$ are the unknowns. \\
We also observe that the fixed point $\bz_1\in S_1$ undergoes
border collision bifurcations when colliding with $\Stan_2$. We recall that, at this bifurcation
a $T$-periodic orbit exhibing one spike reaches the threshold
a second time by tangential grazing (see Figure~\ref{fig:z1Stan2} and~\ref{fig:z1Stan2_TS}). This
bifurcation curve, shown in light red in
Figure~\ref{fig:2dbif_T5},  has been computed solving
Equations \biftanSn{} for $n=2$ which, in this particular case,
become
\begin{equation}\label{eq:z1_Sigma2tan}
\begin{aligned}
\varphi_A(-t_1;V_1,V_1)&=\varphi_0(T-dT;\varphi_A(dT-t_2-t_1;V_2,V_2))\\
(V_2,V_2)^T&=\varphi_A(t_2;V_r,V_1+\Delta)\\
-V_2+V_0+A&=\frac{-V_2+\theta_\infty(V_2)}{\tau_\theta},
\end{aligned}
\end{equation}
where $V_1$, $V_2$, $t_1$, $t_2$, $b$ and $A$ are the unknowns. \\

For $T=5$ we observe $3$ new co-dimension-two bifurcation points apart
from the one reported in the case $T=0.5$, $C2_1$.
As for $C2_1$, two of these new points
are given by the transition from nonsmooth to smooth grazing
bifurcations. The curve defined by $\bz_0\in \Sigma_1$ transitions from 
$\bz_0\in\Stran_1$ (black curve) to $\bz_0\in\Stan_1$ (gray curve) at the point 
labeled as $C2_2$ (see zoomed box in Figure~\ref{fig:2dbif_T5}). At this
point, both Equations \biftanSn{} and \biftranSn{} are
simultaneously satisfied for $n=1$ and hence this is a
co-dimension-two bifurcation point. At these parameter values a non-spiking $T$-periodic orbit tangentially
grazes the threshold at $t=dT$ (see Figures~\ref{fig:z0Sigma1tantrans} and~\ref{fig:z0Sigma1tantrans_TS}).  Something similar occurs with the
fixed point $\bz_1\in S_1$: a grazing bifurcation transitions from 
nonsmooth (dark red) to smooth type (light red) at the point $C2_3$
($\bz_1\in\Stran_2\cap\Stan_2$). At this point, both
Equations \biftanSn{} and \biftranSn{} are simultaneously
satisfied for $n=2$ and hence this is a co-dimension-two bifurcation
point. At
these parameter values a $T$-periodic orbit exhibiting one spike grazes a
second time the threshold at $t=dT$, and does it tangencially (see Figures~\ref{fig:z1Stantrans2}
and~\ref{fig:z1Stantrans2_TS}).\\
The third co-dimension-two bifurcation point, $C2_4$, is of 
different type. Indeed, it is given by the intersection of the bifurcation
curves defined by $\bz_1\in\Stan_1$ and $\bz_1\in\Stan_2$. At these parameter values a $T$-periodic orbit
tangentially grazes the threshold twice (see Figures~\ref{fig:z1Sigma1tanSigma2tan}
and~\ref{fig:z1Sigma1tanSigma2tan_TS}).\\
As in the previous case, the curves defined by $\bz_0\in\Stran_1$ and
$\bz_1\in\Stran_1$ cross transversally defining four regions in the
parameter space regarding their existence. In two of them only
one fixed point exists (either $\bz_0$ or $\bz_1$), in another one both coexist and in the fourth
one none of them exist. In the latter region one finds higher periodic
orbits (see Section~\ref{sec:perorbits}). In the case where both
fixed points coexist one finds bi-stability, as both are
attracting.

\subsection{Periodic orbits of the stroboscopic map and bifurcations}\label{sec:perorbits}
In the previous section we have found the curves on the parameter
space $(b,A)$ where the fixed points $\bz_0$ and $\bz_1$ of the stroboscopic map
$\s$ undergo border collision 
bifurcations. As $\bz_0$ and $\bz_1$ collide with $\Sigma_1$ and disappear there might appear periodic orbits of the map visiting both $S_0$ and $S_1$. 
In order to explore the existence of such periodic orbits, we consider a set of initial conditions on the subthreshold regime (regions $S_0$ and $S_1$) 
and integrate them for several periods to identify the attracting periodic orbits of the system stepping on $S_0$ and $S_1$. See Appendix~\ref{ap:per_orb} 
for the numerical details. Of course, the same exploration can be done for orbits stepping on $S_n$, $n \geq 2$, 
but for the purposes of this paper we focus only on $S_0$ and $S_1$.

For $T=0.5$ we computed the number of attracting periodic
orbits of the stroboscopic map (see Figure~\ref{fig:2dbif_T05_numop_all})
and their periods (see Figure~\ref{fig:2dbif_T05_numop_period}).
Notice that several periodic orbits coexist for many parameter values. Hence, whenever there are several periodic orbits the colour in
Figure~\ref{fig:2dbif_T05_numop_period} has been modified 
in order to reproduce the effect of the intersection. As with the fixed point $\bz_1$, periodic orbits $(z_1,\ldots,z_n) \in (S_0 \cup S_1)^n$ 
of the stroboscopic map appear and disappear due to 
collisions of a point of the orbit $z_i$ with the border $\Sigma_1$ and $\Sigma_2$, 
bounding the region of existence of a given periodic orbit. Take for instance the orbit of period 2 $(z_1,z_2) \in (S_0 \cup S_1)^2$
(corresponding to the orbit with symbolic sequence 01, i.e $z_1 \in S_0$ and $z_2 \in S_1$, see Definition~\ref{def:symb_iti}), which 
can be found in the region colored in gray (plus intersections) in
Figure~\ref{fig:2dbif_T05_numop_period}. Notice that the shape of this region resembles that of the region 
of existence of $\bz_1$ in Figure~\ref{fig:2dbif_T0d5} (recall that this region is bounded by the red curves).
Clearly, the existence regions for different periodic orbits overlap as $b$ increases giving rise to regions with
multiple coexistence of periodic orbits.

For small values of $b$ (close to $0$), we observe that periodic orbits exist only in
the region where the stroboscopic does not have fixed points and,
moreover, these periodic orbits are unique.  Thus, for a fixed small $b$, as the amplitude increases the fixed point $\bz_0$ disappears
through a border collision bifurcation with $\Stran_1$ (black curve)
and a unique periodic orbit appears, undergoing most likely a period-adding
bifurcation (see also Figure~\ref{fig:periods_brut-force}) until the
fixed point $\bz_1$ appears through a border collision bifurcation
with $\Stran_1$ (see Section~\ref{sec:fixed_points_numeric}).  The
periodic orbits in this region are organized by bifurcation structures
that resemble the period-adding bifurcation structure of
$1$-dimensional circle or discontinuous maps (see~\cite{GraAlsKru17}).
More precisely, their ``rotation number'' (see Definition~\ref{def:rotn} and
Remark~\ref{rem:rotn})  resembles the devil's staircase, symbolic
sequences of periodic orbits are glued through gluing bifurcations and
their periods are added.  See Figure~\ref{fig:periods_brut-force} and Figure~\ref{fig:rotnums_b0d1},
where we show the periods and the ``rotation number'', respectively, of the 
periodic orbits along the 1-dimensional scan for $b=0.1$ (labeled in
Figure~\ref{fig:2dbif_T05_numop_all}). 

\begin{remark}
We emphasize that we refer to period-adding-like or
cascade of gluing bifurcations when we cannot assess that we have the infinite number of bifurcation 
curves that separate the regions of existence of periodic orbits or,
equivalently, a continuous curve of ``rotation numbers'', showing a devil's staircase through the complete Farey
tree.
\end{remark}

\begin{figure}
\begin{center}
\begin{tabular}{cc}
\subfigure[]{\includegraphics[angle=-90,width=0.5\textwidth]{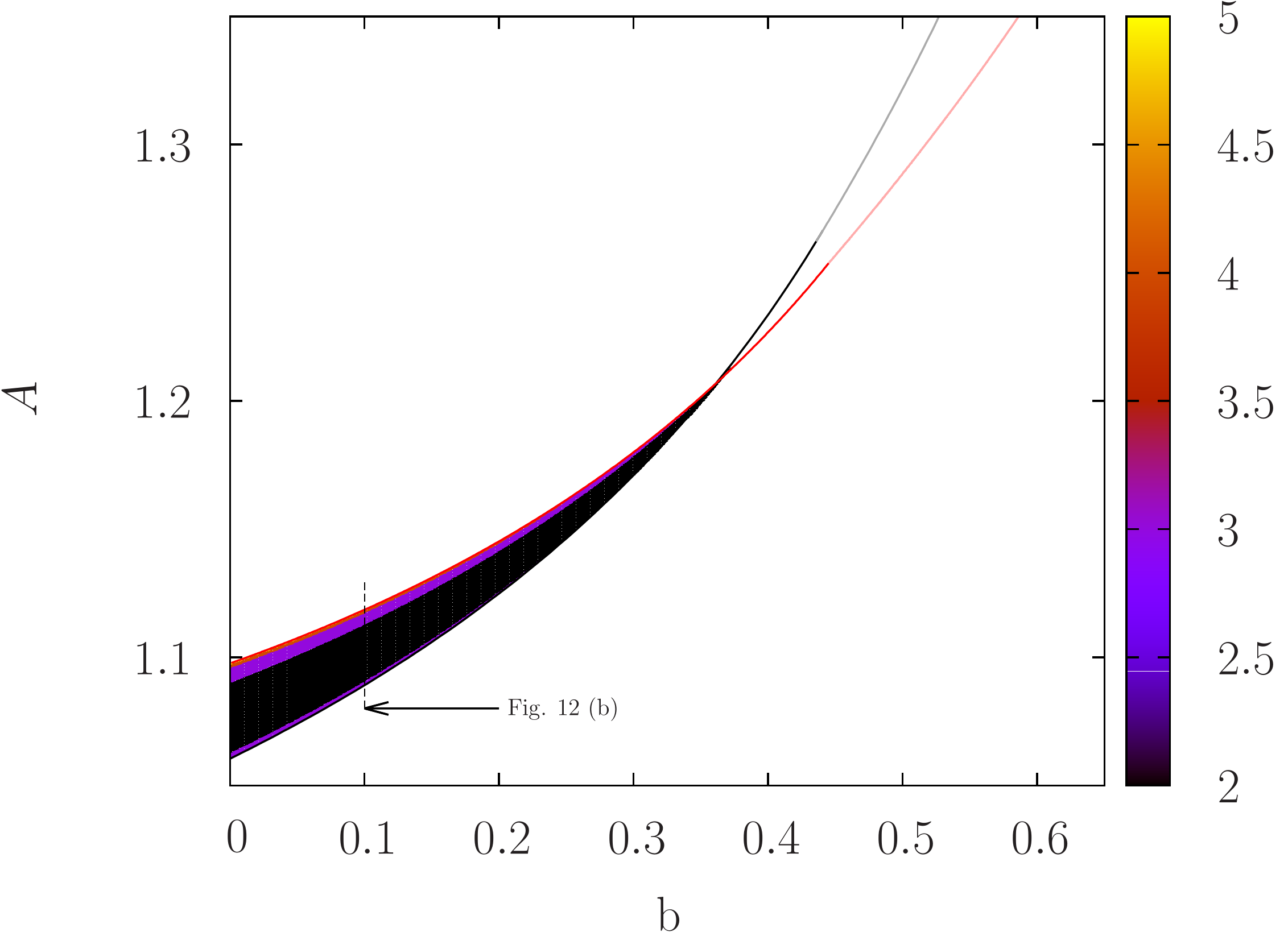}}&
\subfigure[\label{fig:periods_brut_force_T5}]{\includegraphics[angle=-90,width=0.5\textwidth]{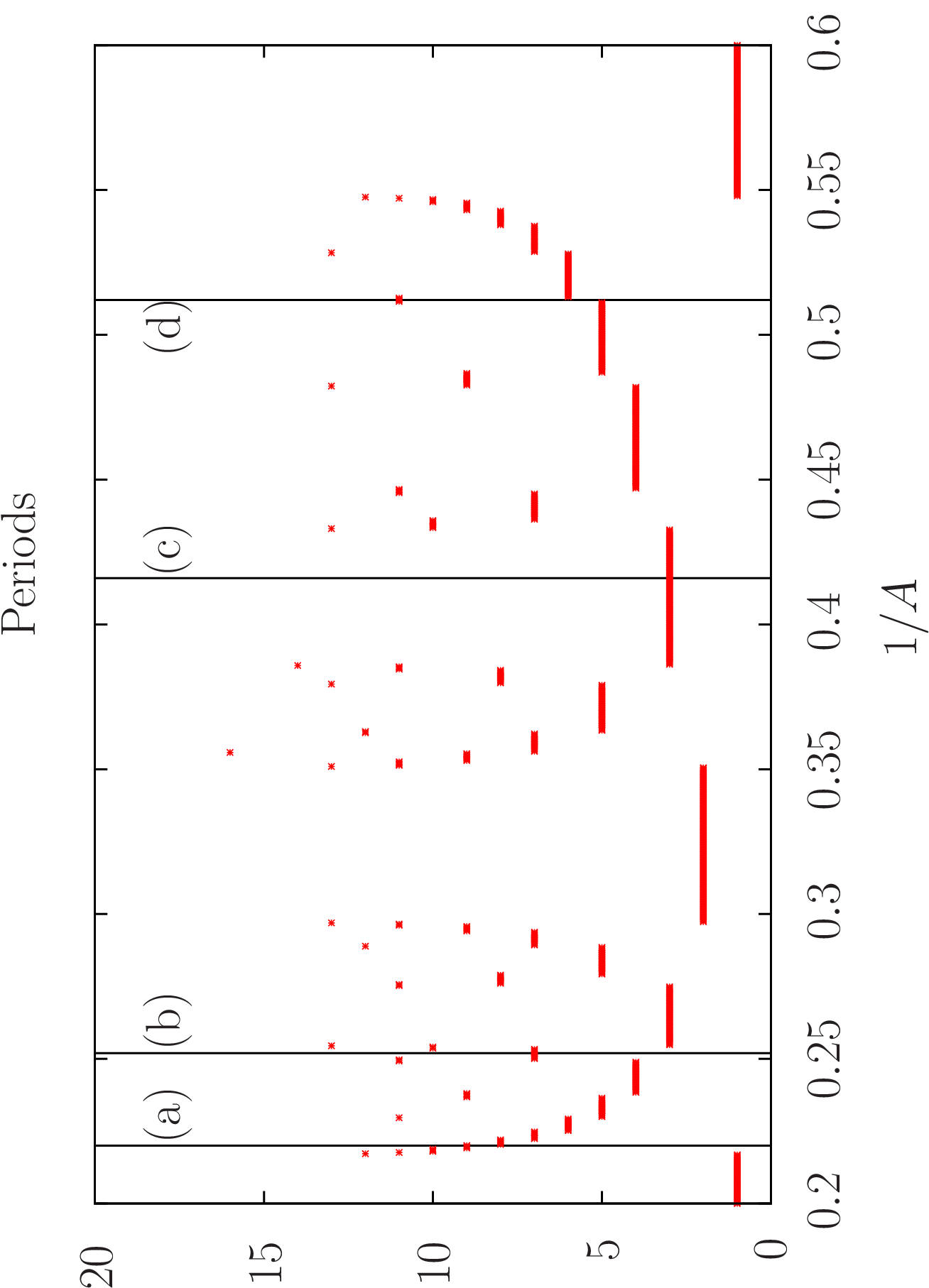}}\\
\end{tabular}
\end{center}
\caption{(a) Periods of unique maximin periodic orbits of the stroboscopic map for $T=5$ stepping on $S_0$ and $S_1$. 
We include also the bifurcation curves computed in Figure~\ref{fig:2dbif_T5}. 
(b) Periods of the periodic orbits found by varying the amplitude $A$ along the vertical line $b=0.1$ as indicated in panel (a).}
\label{fig:2dbif_T5_maximin_periods}
\end{figure}

For intermediate and large values of $b$ (approximately $0.5$ and above) there exist multiple periodic orbits of the stroboscopic map that coexist with fixed points. 
Indeed, the regions of existence of periodic orbits
expand towards the regions where $\bz_0$ or $\bz_1$ also exist, while at the same time 
start to intersect between them, showing multistability. Moreover,
many branches of ``rotation numbers'' are lost, leaving the Farey tree incomplete.
See for example a one-dimensional scan for $b=0.55$ in
Figure~\ref{fig:periods_brut_force_T5}, where many periodic orbits are
no longer found and one finds co-existence instead. Consequently,
the ``rotation number'' (shown in Figure~\ref{fig:rotnums_b0d55}) is
discontinuous, leading to the coexistence of periodic orbits and
overlapping of rotation numbers.

For $T=5$ we observe that for all values of $b$ between 0 and 1, periodic orbits only exist in the region confined between two nonsmooth border collision bifurcations, 
corresponding to the disappearance of the fixed point $\bz_0$ ($\bz_0$ impacts $\Stran_1$) and the appearance of the fixed point $\bz_1$ ($\bz_1$ impacts $\Stran_1$). 
See Figure~\ref{fig:2dbif_T5_maximin_periods}(a). 
In this case periodic orbits are all unique: we do not observe coexistence of several periodic orbits or coexistence of periodic orbits with fixed points. 
For a fixed value of $b$, as the amplitude $A$ increases these periodic
orbits undergo several gluing bifurcations (see also
Figure~\ref{fig:periods_brut_force_T5}). However, in this case our
numerical computations suggest that 
the Farey tree is incomplete and the curve or ``rotation numbers''
shows discontinuities without the overlapping observed in the
previous case.

Notice that, as opposed to the case $T=0.5$, periodic orbits are confined in a very small region of the bifurcation diagram. Alongside, for the 
case $T=0.5$ most of the border collision bifurcations correspond to
collisions with $\Stran_1$ and $\Stran_2$, i.e., nonsmooth grazing
bifurcations. In this case, it is expected that close to a border collision the dynamics of the map will map points of $S_0$ to $S_1$ and 
viceversa, and therefore there might appear periodic orbits whose iterates step on both regions $S_0$ and $S_1$. However, for $T=5$, for intermediate and large values of 
$b$ border collision bifurcations correspond to collisions with
$\Stan_1$ and $\Stan_2$, and in this case we do
not find periodic orbits stepping only on $S_0$ and $S_1$. Instead, it seems that fixed points $\bz_0$ and 
$\bz_1$ coexist, thus possibly preventing the existence of periodic orbits. 
The further exploration of the dynamics close to border collision
bifurcations with $\Stan_1$ and $\Stan_2$ lies ahead
(see Section~\ref{sec:discussion} for a discussion).

\subsection{Maximin itineraries}\label{sec:maximin_po}


In this section we explore the maximin properties of the computed periodic orbits for the stroboscopic map (see
equivalent Definitions~\ref{def:maximin}
and~\ref{def:well-ordered_symbolic-sequence}).
Using a simple algorithm (see Appendix~\ref{ap:per_orb} for details), we
find that all the itineraries are maximin. In
Figure~\ref{fig:periodic_orbits} we show the time series of the
periodic orbits obtained for the parameter values labeled in
Figure~\ref{fig:periods_brut-force} for $T=0.5$. Their symbolic
itineraries are $01^5$, $0101^2$, $001(01)^3$ and $0^71$, which are
all maximin. Recalling Definition~\ref{def:symb_iti}, symbol $0$ is
used when no spike is produced ($z_i\in S_0$), while $1$ means that
one spike is produced ($z_i\in S_1$). To study the existence of
maximin itineraries more rigorously we wonder if the conditions of
Theorem~\ref{theo:quasi-contraction} are satisfied.  Notice though
that conditions {\em i)--iii)} are difficult to check explicitly and
therefore we designed an algorithm to check them numerically. Next, we describe the numerical
algorithm and discuss the domain of application of the theoretical
result.


For a given value of $A$, we first numerically find sets $E_i$ that satisfy 
hypothesis \emph{i)} of Theorem~\ref{theo:quasi-contraction}.
We construct these sets to be as small as possible and later we check
whether they satisfy conditions \emph{ii)} and \emph{iii)}. For
convenience, we choose these sets to be
quadrilaterals whose union is a convex polygon. More
complex geometries are of course possible, although they would
significantly complicate the algorithm without guarantee of better
results.

To construct the sets $E_0$ and $E_1$, we first take a small segment
$\gamma\subset\Sigma_1$ ``close'' to the periodic orbit found by
direct simulation. This segment is iterated by $\s_0$ and $\s_1$. We
then consider the two quadrilaterals formed by the
segments $\gamma$ and $\s_0(\gamma)$, and $\gamma$ and $\s_1(\gamma)$
(see Figure~\ref{fig:sets_Ei}). We then grow the segment $\gamma$
until the union of these two quadrilaterals is a
convex polygon. If this cannot be done, then we stop the algorithm
and assume we could not find the desired sets $E_0$ and $E_1$.
If we succeed, we check whether the images of the two
quadrilateral candidates are contained in their
union (the convex polygon).  The fact that their union is a
convex polygon makes it easier to check this inclusion (see
Remark~\ref{rem:polygons}). If any of the images is not contained in
this union, then we further grow the initial segment $\gamma$ in the
direction that failed and we check again. If at some point we succeed,
then we have found sets $E_i$ satisfying hypothesis \emph{i)} of
Theorem~\ref{theo:quasi-contraction}.  For the four values of $A$
indicated in Figure~\ref{fig:periods_brut-force} we have been able to
find sets $E_0$ and $E_1$ as described (see
Figure~\ref{fig:sets_Ei}).

We then check the contracting condition~{\em ii)}. As discussed in
Section~\ref{sec:extension}, the map $\s_0$ is always contracting, as
it is the composition of two contracting stroboscopic maps
(provided that for $A>0$ the system~\eqref{eq:system_visual} for $I=A$ always possesses an
equilibrium point, eventually virtual).  However, $\s_1$ does not
necessary contract, even if system~\eqref{eq:system_visual} possesses for both $I=A$ and $I=0$ attracting equilibrium points, due to the collision with the threshold and the
reset condition. We check its contractiveness by computing the
differential $D\s_1$ (as described in Appendix~\ref{sec:diff_strobo})
and its eigenvalues in a mesh of points in $E_1$.  If, for all the
points in the mesh, both eigenvalues have modulus less than $1$ we
then assume that condition {\em ii)} is also satisfied. If this
condition is not fulfilled, we then say that we have not been able to
check the conditions of Theorem~\ref{theo:quasi-contraction}. This is
the case in panels (c) and (d) of Figure~\ref{fig:sets_Ei}, for which
we have found points for which the matrix $D\s_1$ has eigenvalues
outside the unit circle.

Finally, it remains to check whether condition {\em iii)} holds. This
is done by checking whether all points in the segment $\gamma\subset
\Sigma_1$ visit $S_0$ and $S_1$ altogether or they split
after intersecting $\Sigma_1$( $\s^n(\gamma)\cap \Sigma_1$ for
some $n>0$). In other words, we numerically check
whether all points in $\gamma$, are attracted towards the same fixed point of
$\s^p$, where $p$ is the period of the periodic orbit found by
iteration.  In Figure~\ref{fig:sets_Ei} we show examples of values of
$A$ for which condition {\em iii)} is  satisfied (panels~(a), (b) and
(d)) and not satisfied (panel (c)).  Thus, we conclude that for
parameter values $(A,b)$ corresponding to cases (a) and (b) we have
proven semi-rigorously the existence of a (locally) unique periodic orbit whose
symbolic itinerary is maximin, by means of numeric validation of the
hypothesis of Theorem~\ref{theo:quasi-contraction}. We want to
acknowledge that this numerical validation does not follow a computer
assisted proof procedure.

We have applied the numerical algorithm described above to all values
of $A$ in Figure~\ref{fig:periods_brut-force}.
In Figure~\ref{fig:satisfying_Gam} we show the values of
$A$ for which we have been able to validate conditions
{\em i)--iii)} of Theorem~\ref{theo:quasi-contraction}.

Notice that hypothesis of Theorem~\ref{theo:quasi-contraction} are difficult to check and very 
restrictive. Indeed, maximin periodic orbits exist far beyond the regions where 
the hypothesis can be checked using the numerical algorithm described above. 
Future work will be devoted to study the 
viability of the application of the techniques developed in~\cite{GleKowNor12} 
to cases with multistability.

\begin{remark}
This algorithm has more chances to suceed if the segment
$\gamma\subset\Stran_1$, which is the case in this example, as the
computed periodic orbits are located between the curves defining
$\bz_0\in\Stran_1$ and $\bz_1\in\Stran_1$.
\end{remark}
\begin{remark}\label{rem:polygons}
This algorithm takes advantage of the fact that the sets $E_i$
are quadrilaterals and that $E_0\cup E_1$ is a convex polygon in order
to easily check the inclusions $\s_0(E_0)\subset E_0\cup E_1$ and
$\s_1(E_1)\subset E_0\cup E_1$.  This is done by checking that the
images of the vertices are contained in $E_0\cup E_1$.  However, this
guarantees the inlusions $\s_i(E_i)\subset E_0\cup E_1$ only if 
these images are also quadrilaterals. This is in principle only true
if the flow is linear. Although in our case this is not true, the
non-linearity of the maps $\s_i$ is not significant and images of lines become
almost lines. The committed error by only checking the vertices is
drastically reduced by also checking the image of a middle point in
the segments of the quadrilaterals.
\end{remark}

\section{Discussion}\label{sec:discussion}
In this paper, we have studied the dynamics of hybrid systems submitted to a periodic forcing consisting of 
a square-wave pulse. 
In particular, we have studied a model
consisting of a leaky integrate-and-fire model with a dynamic threshold, combined with a reset rule 
that is applied whenever the trajectory crosses the threshold manifold. In our analysis, we have considered the 
stroboscopic map, which is a piecewise-smooth discontinuous two-dimensional map, for which we can apply existing theoretical results
by Gambaudo et al. \cite{GamTre88,GamGleTre88} for ``quasi-contractions'' (see Theorem~\ref{theo:quasi-contraction}). Thus, we have been able to ``prove'', in combination with numerics, 
the existence of periodic orbits of maximin type for the stroboscopic map for certain parameter values. 
Moreover, we have explored numerically nonsmooth bifurcations (border
collisions of fixed points and gluing bifurcations of $T$-periodic orbits), 
that provide a wider description of the dynamics.

Using particular geometries for the sets $E_i$
(quadrilaterals with $E_0\cup E_1$ being a convex polygon), we have shown
that the existing theory has very restrictive hypothesis, since we have
numerically observed the existence of maximin unique periodic orbits
even for cases that do not satisfy the hypothesis of Theorem~\ref{theo:quasi-contraction}.
Moreover, these hypothesis are difficult to check analytically.
Indeed, we have designed a semi-rigorous numerical procedure to check
them.  It is possible that using more complex geometries and
topologies for the sets $E_i$ one could obtain better results in checking these
hypothesis, although this may signifficantly complicate the numerical
algorithms. Thus, the existing theoretical results for
two-dimensional piecewise continuous maps are still very limited to
provide a complete description of the dynamics.

In our numerical analysis of nonsmooth bifurcations, 
we have varied the parameter $b$ in Equation~\eqref{eq:theta_infty} which sets the system in two different dynamical regimes. 
Thus, for small values of $b$ (close to $0$), in response to a constant input of sufficiently large amplitude, the system  
shows repetitive spiking (\emph{tonic regime}). However, as $b$ increases, a constant input cannot generate repetitive firing, only a few 
spikes before returning to resting potential (\emph{phasic regime}) \cite{MenHugRin12, RinHug13, HugMenRin17}. 
We have forced the system with a square-wave periodic pulse and explored fixed points and periodic solutions of the stroboscopic 
map in these two regimes as the amplitude $A$ increases.

%

In the tonic regime, we have observed a unique globally attractive maximin
periodic orbit with a period that, as the amplitude is increased, undergoes several
gluing bifurcations mimicking the so-called period adding bifurcations for
1-dimensional maps.  Thus, the
transition from a $T$-periodic orbit with no spikes (fixed
point $\bz_0$) and a $T$-periodic orbit with a spike (fixed
point $\bz_1$) occurs through a complex mechanism of concatenation of
sequences of periodic orbits whose ``rotation number''
is ordered as in
the Farey tree. When we compare the results for two different input
frequencies we observe different behaviours. For high frequency
($T=0.5$), their ``rotation number'' evolves, up to our numerical
accuracy, as in the $1$-dimensional case; that is, it evolves
continuously along the Farey tree showing a devil's staircase. This
at least the case for periodic orbits with periods up to 30. 
In contrast to the high frequency case, for lower input frequency
($T=5$), the ``rotation number'' becomes discontinuous. In other
words, the Farey tree of ``rotation numbers'' is not complete. In both
cases, the cascade of gluing bifurcations is confined in the region
between two border collision bifurcations (where $\bz_0$ and $\bz_1$ do
not exist).
 
The behavior drastically changes as $b$ increases and the
dynamical regime becomes phasic. The dynamics here is more complex and
very different for $T=5$ and $T=0.5$. Thus, for $T=0.5$ the
organization of periodic orbits in a Farey tree gets destroyed by the
overlap of neighboring periodic orbits (including periodic orbits of
period 1 corresponding to $\bz_0$ and $\bz_1$), leading first to
multistability and later to the disappearance of these periodic orbits
of period higher than 1. In the case $T=5$, periodic orbits of period
higher than 1 (that step on $S_0$ and $S_1$) cannot be found beyond
the regions where $\bz_0$ and $\bz_1$ do not exist.  Notice that in this
setting, for small values of $A$ the system displays a
$T$-periodic orbit with no spikes (fixed
point $\bz_0$) and as $A$ increases there appears a $T$-periodic
orbit with a spike (fixed point $\bz_1$), and both $T$-periodic orbits coexist.
This behavior agrees with the observation in \cite{MenHugRin12} (for a 
periodic rectified sinusoidal input), where the transition 
from firing patterns consisting of a $T$-periodic orbit
with no spikes to a $T$-periodic orbit with 1 spike, is more abrupt for
low input frequencies.
 
Notice that $b$ small sets the system in a dynamical regime (tonic) which can be 
modeled with a 1-dimensional integrate-and-fire
system with a fixed threshold. However, as $b$ increases the dynamical regime
becomes phasic, and a one-dimensional integrate-and-fire model
cannot reproduce these dynamics, thus suggesting that the results for
larger values of $b$ show characteristics of 2-dimensional systems.
Alongside, conditions for the existence of a period-adding bifurcation (for
which the rotation number is continuous and evolves showing a
devil's staircase along the Farey tree) have been
established theoretically for 1-dimensional piecewise maps \cite{GraKruCle14}. However, 
it still remains an open question whether it can occur in 2-dimensional maps
and under which conditions. We believe that our results might serve 
as a motivation and starting point to prove the conditions that guarantee the existence 
of period-adding structures in 2-dimensional hybrid systems.

We recall that, in the nonsmooth literature, the 
concept of grazing bifurcation refers to a
periodic orbit of a Filippov system that grazes a switching manifold (where the
vector field has a discontinuity). As the vector field is smooth
outside of the switching manifold, this grazing happens to be
tangent. In some sense, this classical grazing bifurcation can be
related to our case. After identifying the manifolds $\TM\sim\RM$
such that
\begin{equation*}
\TM\ni z_\TM\sim z_\RM\in\RM\text{ if }z_\RM=R(z_\TM),
\end{equation*}
and adding time as a new variable one gets a Filippov
system with three swtching manifolds ($\TM\sim\RM$, $\{t=dT\mod T\}$
and $\{t=0\mod T\}$). Then, what
we called smooth grazing bifurcation corresponds to grazing the
manifold $\TM\sim\RM$ only, while the nonsmooth grazing bifurcation corresponds to
simultaneously grazing the manifold $\TM\sim\RM$ and crossing
$\{t=dT\mod T\}$. The first case becomes similar to the classical
grazing with the addition of a second discontinuity surface which is
not influencing the bifurcation~\cite{DonHis04}. In our work,
however, we are not only interested in the conditions for which a
$T$-periodic orbit bifurcates through grazing, but also on how are the
dynamics after this bifurcation. We have shown that rich dynamics
involving gluing bifurcations are generated when the grazing and
crossing simultaneously occur.

Also regarding the nonsmooth literature, we emphasize that
classical
results~\cite{BerBudCha01a,DanNor00,DanZha05,Nordmark97,Nordmark01,NorKow06,Sim14,ZhaDan06,ThoDan06,Kowalczyk05}
cannot be applied to our case mainly due to the periodic forcing,
which makes the Poincar\'e map onto the grazed manifold 2-dimensional.
Indeed, this is precisely the main reason why one finds so rich
dynamics in this case. However, regarding multistability close to the
grazing bifurcation, a similar phenomenon is also found for autonomous
systems leading to 1-dimensional Poincar\'e\linebreak
maps~\cite{GleKowNor12,GleKowNor16} although one requires sliding
motion along the switching manifold, which cannot be found in our
case. We recall that, in this work we have focused on periodic orbits
stepping on $S_0$ and $S_1$ only. A further analysis of the
bifurcation structure close to the boundaries $\Stan_1$ and $\Stan_2$
may reveal that other type of periodic orbits also stepping on the
sets $S_i$, $i>1$, may also appear and may be related with previous
results on grazing bifurcations leading to multistability. 

We finally point out that, also due to the periodic forcing, our
analysis is different than others in the neuroscience
literature~\cite{TouboulB09, Rubinetal17a,JimMihBroNieRub13}. Instead
of considering the so-called firing or adaptation map (which is indeed
the Poincar\'e map onto the switching manifold), our analysis is based
on the stroboscopic map, which becomes a piecewise-smooth
discontinuous map. The impact map has the advantage that is
continuous but, as the system is non-autonomous, one would need to
keep track of the time variable at the threshold manifold, which
implies solving transcendental equations.  Moreover, unlike the
Poincar\'e map, the stroboscopic map is defined everywhere and hence
captures subthreshold dynamics (dynamics not interacting with the
threshold). Hence, this latter map is the most suitable one to provide
a complete description of the dynamics for the systems considered in this paper.

\newcommand{\etalchar}[1]{$^{#1}$}
\def\zh{Zh}\def\yu{Yu}\def\ya{Ya}

\section*{Appendix}
\appendix 

This appendix contains the details of the numerical computations performed along the paper. In particular, we describe the 
computation of the switching manifolds for the stroboscopic map, its differential, the curves corresponding to 
border collision bifurcations of the fixed points of the stroboscopic map and the computation of the periodic orbits of the 
stroboscopic map.

\section{Computation of the switching manifold}\label{ap:sigma}
In this section we explain how to compute the switching manifolds $\Sigma_n$, and therefore 
the areas $S_n$ in the subthreshold domain $\D$.
In the case of $\Stran_n$, we take points on $\TM$ and
integrate them backwards in time for $t=dT$. If the trajectory hits
the reset manifold $\RM$, we apply $R^{-1}$ (in the example considered, subtract $\Delta$ to the
$\theta$-variable and set $V=\theta$), count one hit and continue integrating the
trajectory backwards for $t=dT$. If the trajectory hits $\RM$ again for $t<dT$, we
count another hit and we apply the same procedure as before. We keep the last point. It 
will belong to $\Stran_{n+1}$, where $n$ counts the number of hits with
the manifold $\RM$. Of course, the trajectories always lie inside the
subthreshold region. On the contrary, if one of these trajectories
lies on the superthreshold region when integrating backwards, we may
suspect that $\Stan_n \neq \varnothing$. Then, we compute the
point on the threshold manifold $\TM$ where the trajectory is tangent. In the example considered, it must
satisfy $V=\theta$ and $(-1,1) \cdot (f(V,\theta)+(A,0))=0$, or, equivalently, $(V-V_0-A)\tau_\theta - V +
\theta_{\infty}(V)=0$. Once we find the point $(V,V)$ that solves the equations, we integrate this point backwards in time until
$t=dT$. A point on this trajectory belongs to $\Stan_{n+1}$, where $n$ corresponds to the number of times the trajectory has hit the reset manifold
before reaching the given point when integrated backwards. 

\section{The differential of the stroboscopic map}\label{sec:diff_strobo}

In this section we give details of the computation of the differential of 
the stroboscopic map $\s$ in $S_0$ and $S_1$, that is, $D\s_0$
and $D\s_1$, for the general setting described in Section~\ref{sec:general_setting}. The differential will be used to compute fixed points through a
Newton method as well as their stability. Recall that $\s_0$ is the
composition of two stroboscopic maps (see Equation~\eqref{eq:strobo_s0}). Since $I(t)$ is a
piecewise-constant function, the Jacobian matrix of the stroboscopic
map can be computed by solving the variational equations of system~\eqref{eq:general_system} for $I$ constant. Calling
$z=(x,y)$, these variational equations become
\begin{equation}
\frac{d}{dt}\delta_A(t;z)=
\left( 
\begin{array}[]{cc}
\frac{\partial}{\partial x}f^x(\varphi_A(t;z))&\frac{\partial}{\partial
y}f^x(\varphi_A(t;z))\\
\frac{\partial}{\partial x}f^y(\varphi_A(t;z))&\frac{\partial}{\partial
y}f^y(\varphi_A(t;z))\\
\end{array}
 \right)\delta_A (t;z)
\label{eq:variational}
\end{equation}
and
\begin{equation*}
\delta_A(t;z)=\left( 
\begin{array}{cc}
\frac{\partial}{\partial x}\varphi^x_A(t;z) &
\frac{\partial}{\partial y}\varphi^x_A(t;z))\\
\frac{\partial}{\partial x}\varphi^y_A(t;z))&
\frac{\partial}{\partial y}\varphi^y_A(t;z))
\end{array}
 \right),
\end{equation*}
where $\varphi_A^{x}$, $\varphi_A^{y}$, $f^x$ and $f^y$ denote the
$x$ and $y$ component of the flow $\varphi_A$ and the vector field $f$, respectively.
Equation~\eqref{eq:variational} must be integrated along with the flow
$\varphi_A(t;z)$ using the identity matrix as initial condition for
$\delta_A$, i.e $\delta_A(0;z)=Id$.\\
Hence, for $z\in S_0$, we have
\begin{equation}
D\s_0(z)=\delta_0\left(T-dT;\varphi_A\left(dT;z\right)\right)\cdot\delta_A(dT;z).
\label{eq:jacobian_S0}
\end{equation}
In order to compute the differential $D\s_1$ (see definition of $\s_1$ in Equation~\eqref{eq:s_in_S1}), we need to
take into account that the firing time, $t_*$, depends on the initial
conditions, i.e. $t_*=t(z)$. Thus, we compute $Dt$ by means of
applying the Implicit Function Theorem to the equation
\begin{equation}
F(t,z):=h(\varphi_A(t;z))=0.
\label{eq:threshold_condition}
\end{equation}
Let $(t_*,z)$ be a solution of Eq.~\eqref{eq:threshold_condition};
that is, $\varphi_A(t_*;z)\in\TM$. Then, assuming that
\begin{equation}
 \frac{\partial F}{\partial t} (t_*,z)= \nabla h
(\varphi_A(t_*; z)) \frac{d \varphi_A}{dt}(t_{*};z) \ne 0,
\label{eq:non-degeneracy}
\end{equation}
we obtain $t^*=t(z)$ and 
\begin{equation*}
Dt(z)=-\frac{\partial F/ \partial z}{\partial F / \partial t}=-\frac{\nabla h(\varphi_A(t_*;z))\cdot\delta_A(t_*;z)}{
\nabla h (\varphi_A(t_*;z))\cdot (f(z) + vA)}.
\end{equation*}
Note that condition~\eqref{eq:non-degeneracy} is equivalent to
requiring that the flow is not tangent to $\TM$ at
$\varphi_A(t_*;z)$, and hence this is only valid when the spike is
given by a transversal crossing between the flow and $\TM$.\\
Then, the differentials of the maps~\maps{} are
\begin{equation*}
DP_1(z)=\left( 
\begin{array}{c}
\delta_A(t_*;z) + (f(z)+vA) Dt(z) \\ Dt(z)
\end{array}
 \right)
\end{equation*}
\begin{equation*}
D\tRmap(z,t)=\left( 
\begin{array}{ccc}
\multicolumn{2}{c}{\multirow{2}{*}{$D\Rmap(z)$}}
&0\\
&&0\\
0&0&1
\end{array}
 \right)
\end{equation*}
\begin{equation*}
D\tP_2(z,t)=\left( 
\begin{array}{ccc}
\delta_A(dT-t;z)& -(f(z)+vA)
\end{array}
 \right)
\end{equation*}
\begin{equation*}
DP_3(z)=\delta_0(T-dT;z).
\end{equation*}
Hence, we get
\begin{equation*}
D\s_1(z)=DP_3(z_2)D\tP_2(z_r,t_1)D\tRmap(z_1,t_1)DP_1(z),
\end{equation*}
where
\begin{equation*}
(z_1,t_1)=P_1(z)\quad (z_r,t_1)=\tRmap(z_1,t_1)\quad
z_2=\tP_2(z_r,t_1).
\end{equation*}

\begin{remark}
Notice that in the computation of $DP_1$ we take into account that the original and perturbed 
orbit do not cross the switching manifold at the same time. 
The method described herein to compute $DP_1$ is equivalent to the one
developed in~\cite{AizGan58}, later adopted by the nonsmooth and
mechanical-engineering community and named
{\em saltation} (or {\em correction}) matrix~\cite{LeiNij04}.  
\end{remark}

\section{Numerical computation of bifurcation curves}\label{sec:compu_bif-curves}

The set of equations describing the bifurcation curves consist always of $n$ equations and $n+1$ unknowns, 
defining a curve on an $n$-dimensional space. In this section we will present a
predictor-corrector method to compute such curves. The corrector step
is based on a Newton method and Lagrange multipliers. We refer
to~\cite{AllGeo03} for further details. In
\begin{center}
\url{https://github.com/a-granados/ContLagMult}
\end{center}
we provide the code of a general purpose implementation of this method
in $C^{++}$. In this repository we also provide, to serve as example, the
details to compute all bifurcation curves discussed in
Section~\ref{sec:fixed_points_numeric}.\\
Let $\Gc(w)=0$ be the set of $n$ equations with $w \in \mathbb{R}^{n+1}$. Assume that there exists $w^{*} \in \mathbb{R}^{n+1}$ such that
$\Vert \Gc(w^*) \Vert$ is small. We look for an improved solution by means of a Newton method. Thus, we look for $\Delta w$ such that 
$w^*+\Delta w$ solves equation $\Gc(w)=0$ up to square order, i.e
we impose that $\Gc (w^*) + D \Gc (w^*) \Delta w =0$, with the extra condition that $\Vert \Delta w \Vert_2$ is a minimum.
It becomes a poblem of finding the local minima subject to an equality constraint, which we will solve using the method of Lagrange multipliers. 
We introduce the Lagrange multiplier $\mu \in \mathbb{R}^n$ and the Lagrange function:
\[ \mathcal{L}= \Delta w^T \Delta w + \mu^T (D \Gc (w^*) \Delta w + \Gc (w^*) ). \]
Now imposing that $\partial \mathcal{L}/\partial (\Delta w) =0$ and using the equality constraint, we are left with the 
following system of equations for $\mu$ and $\Delta w$:
\[
\begin{array}{rcl}
2 \Delta w^T + \mu^T D \Gc(w^*) & = & 0, \\
\Gc (w^*) + D \Gc (w^*) \Delta w &=&0, \\
\end{array}
\]
from where we obtain $\Delta w= - D \Gc(w^*)^T (D \Gc(w^*) D \Gc(w^*)^T) \Gc(w^*)$.
We repeat the process until we obtain a point $w_{new}$ that satisfies equation $\Gc$ up to the desired error. 
This provides a point on the curve. The next step is to compute another point on the curve.
To do so, we compute the tangent vector to the curve at the point $w_{new}$ and
take a point at a distance $\delta$ along this direction. Notice that the tangent vector is given by an 
element of the kernel of $D \Gc (w_{new})$. This will be the initial seed to repeat the 
Newton method procedure described above.

To illustrate the method, we provide the details for the computation
of two of the bifurcation curves reported in
Section~\ref{sec:fixed_points_numeric}, corresponding to $\bz_1\in\Stran_1$ and
$\bz_1\in\Stan_1$; the rest of the bifurcation curves can be obtained
proceeding similarly.\\
For the case $\bz_1\in\Stran_1$, the function $\Gc$ is given by
Equations~\eqref{eq:z1_Sigma2trans}, and becomes
\begin{equation*}
\Gc(\omega)=
\left(
\begin{array}{c}
\varphi_A(-t_1;V_1,V_1)-\varphi_0(T-dT;V_2,V_2)\\
(V_2,V_2)^T-\varphi_A(dT-t_1;V_r,V_1+\Delta)
\end{array}
\right),
\end{equation*}
with $\omega=(V_1,V_2,t_1,b,A)\in\RR^5$. Note that $\Gc(w) \in\RR^4$, as
the flow is 2-dimensional. The derivatives in $D\Gc$ involving the
variables $V_1$, $V_2$ are obtained integrating the variational
equations of system~\eqref{eq:system_visual}, and the ones involving
$t_1$ are given by the vector field. Regarding the derivatives with respect to
parameters $b$ and $A$, one needs to consider them as variables of the
system and increase its dimension from $2$ to $4$ by adding the
equations $\dot{b}=0$ and $\dot{A}=0$. Then, the derivatives with respect
to $b$ and $A$ are obtained from the variational equations of this new extended
system. Using the notation of Appendix~\ref{sec:diff_strobo}, the
variational equations for the extended system become
\begin{equation}
\frac{d}{dt}\td_A=
\left( 
\begin{array}[]{cccc}
\multicolumn{2}{c}{\multirow{2}{*}{$Df(\varphi_A(t;z))$}}&0&\frac{\partial}{\partial A}f^V(\varphi_A(t;z)) \\
&& \frac{\partial}{\partial b}f^\theta(\varphi_A(t;z)) &0\\
0&0&0&0\\
0&0&0&0
\end{array}
 \right)\td_A
\label{eq:variational_bA}
\end{equation}
where 
\begin{equation*}
\td_A(t;z,b,A)=\left( 
\begin{array}{cccc}
\multicolumn{2}{c}{\multirow{2}{*}{$\delta_A(t;z)$}}&
\frac{\partial}{\partial b}\varphi^V_A(t;z) &
\frac{\partial}{\partial A}\varphi^V_A(t;z) \\
&&
\frac{\partial}{\partial b}\varphi^{\theta}_A(t;z) &
\frac{\partial}{\partial A}\varphi^{\theta}_A(t;z) \\
0&0&1&0\\0&0&0&1
\end{array}
 \right),
\end{equation*}
and $\delta_A(t;z)$ is defined in Appendix~\ref{sec:diff_strobo} with
$z=(V,\theta)$.\\
With this notation we get
\begin{equation*}
\frac{\partial}{\partial V_1}\Gc(V_1,V_2,t_1,b,A)=\left( 
\begin{array}[]{c}
\delta_A^{1,1}(-t_1;V_1,V_1)+\delta_A^{1,2}(-t_1;V_1,V_1)\\
\delta_A^{2,1}(-t_1;V_1,V_1)+\delta_A^{2,2}(-t_1;V_1,V_1)\\
-\delta_A^{1,2}(dT-t_1;V_r,V_1+\Delta)\\
-\delta_A^{2,2}(dT-t_1;V_r,V_1+\Delta)
\end{array}
 \right),
\end{equation*}
\begin{equation*}
\frac{\partial}{\partial V_2}\Gc(V_1,V_2,t_1,b,A)=\left( 
\begin{array}[]{c}
-\delta_0^{1,1}(T-dT;V_2,V_2)-\delta_0^{1,2}(T-dT;V_2,V_2) \\
-\delta_0^{2,1}(T-dT;V_2,V_2)-\delta_0^{2,2}(T-dT;V_2,V_2) \\
1\\
1
\end{array}
 \right)
\end{equation*}
\begin{equation*}
\frac{\partial}{\partial t_1}\Gc(V_1,V_2,t_1,b,A)=\left( 
\begin{array}[]{c}
-(f^V(z_1)+A)\\-f^{\theta}(z_1)\\f^V(z_2)+A \\f^{\theta}(z_2)
\end{array}
 \right)
 \begin{array}{l}
  z_1=\varphi_A(-t_1;V_1,V_1)\\z_2=\varphi_A(dT-t_1;V_r,V_1+\Delta),
  \end{array}
 \end{equation*}
\begin{equation*}
\frac{\partial}{\partial b}\Gc(V_1,V_2,t_1,b,A)=\left( 
\begin{array}[]{c}
\td_A^{1,3}(-t_1;V_1,V_1,b,A)-\td_0^{1,3}(T-dT;V_2,V_2,b,A)\\
\td_A^{2,3}(-t_1;V_1,V_1,b,A)-\td_0^{2,3}(T-dT;V_2,V_2,b,A)\\
-\td_A^{1,3}(dT-t_1;V_r,V_1+\Delta,b,A)\\
-\td_A^{2,3}(dT-t_1;V_r,V_1+\Delta,b,A)
\end{array}
 \right),
\end{equation*}
and
\begin{equation*}
\frac{\partial}{\partial A}\Gc(V_1,V_2,t_1,b,A)=\left( 
\begin{array}[]{c}
\td_A^{1,4}(-t_1;V_1,V_1,b,A)\\
\td_A^{2,4}(-t_1;V_1,V_1,b,A)\\
-\td_A^{1,4}(dT-t_1;V_r,V_1+\Delta,b,A)\\
-\td_A^{2,4}(dT-t_1;V_r,V_1+\Delta,b,A)
\end{array}
 \right).
\end{equation*}
Similarly, for the smooth grazing bifurcation corresponding to $\bz_1\in\Stan_2$
the function $\mathcal{G}$ is given by Equations~\eqref{eq:z1_Sigma2tan} and has the form
\begin{equation*}
\Gc(\omega)=
\left(
\begin{array}{c}
\varphi_A(-t_1;(V_1,V_1))-\varphi_0(T-dT;\varphi_A(dT-t_1-t_2;V_2,V_2))\\
(V_2,V_2)^T-\varphi_A(t_2;V_r,V_1+\Delta)\\
-V_2+V_0+A-1/\tau_\theta\left(-V_2+\theta_\infty(V_2) \right)
\end{array}
\right),
\end{equation*}
with $\omega=(V_1,V_2,t_1,t_2,b,A)\in\RR^6$. Note that now
$\Gc\in\RR^5$, and the last equation imposes a tangency when reaching
the threshold at $(V_2,V_2)$. The differential $D\Gc$ is obtained
proceeding as before taking into account that we need to apply the
chain rule due to the composition $\varphi_0\circ\varphi_A$. 

\section{Numerical computation of periodic orbits}\label{ap:per_orb}
In this section we describe the algorithm to compute periodic orbits of the stroboscopic map.  In
\begin{center}
\url{https://github.com/ghuguet/PerOrbVtheta}
\end{center}
we provide the code in $C^{++}$ used to compute the periodic orbits for the stroboscopic map of 
system~\eqref{eq:system_visual}-\eqref{eq:reset_ex}, considered in this paper. 

Periodic orbits of the system are computed by integrating forward a set of initial conditions
on the region $S_0$. For each initial condition we compute $N=100$ iterates of the stroboscopic map and 
we keep the last point $x_f$ in a list of points $L_f$ as well as the previous one $x_a$ ($x_f=\s(x_a)$) in a list of points $L_a$. We 
also keep track of the number of spikes that occur between $x_a$ and $x_f$.
When several initial conditions yield the same points 
$x_a$ and $x_f$, we keep them only once. Once we have computed the iterations for all points in the subthreshold domain, we check whether the points in the list $L_a$
coincides with the points in the list $L_f$. If there is a point $\bar{x}_a$ in $L_a$ which is not 
present in $L_f$, we add $\bar{x}_a$ in the list $L_f$ and $\s^{-1}(\bar{x}_a)$ in the list $L_a$.
Analogously, if a point $\bar{x}_f$ in $L_f$ is not 
present in $L_a$, we add $\bar{x}_f$ in the list $L_a$ and $\s(\bar{x}_f)$ in the list $L_f$.
We repeat this procedure until both lists contain the same points. If after several iterations both lists 
are still different, the program returns the message that the computation of periodic orbits has failed.
Otherwise, once both lists are equal, we order them while constructing a map $\nu$ that maps the position $i$ of the point $x_a$ in the ordered list $L_a$ to the 
position $j$ of its image $x_f=\s(x_a)$ in the ordered list $L_f$, i.e. $\nu(i)=j$. For each pair $(i,j)$ we also keep track of the number of spikes that occur. 
Let us call this number $s_{i}$
A fixed point of the stroboscopic map corresponds to $\nu(i)=i$, while a periodic orbit of period $q$ corresponds to 
\[\nu^q(i)=i, \, \textrm{where} \, \nu^k(i) \neq i, \, \textrm{for} \, k=1,\ldots,q-1.\]
Hence, we use the map $\nu$ to compute the number of fixed points and periodic orbits. Moreover, for each periodic orbit we compute its period $q$ and for 
periodic orbits stepping only on $S_0$ and $S_1$, we keep its itinerary given by $(s_{i}, s_{\nu(i)},\ldots, s_{\nu^{q-1}(i)}) \in \{0,1\}^q$.\\

Using this itinerary we compute whether the periodic orbit is maximin. From a
computational point of view, Definition~\ref{def:maximin} is not
practical to decide whether a symbolic sequence is maximin or not.
Instead, we use the following equivalent definition of maximin using
the notion of $p,q$-ordered sequences
(see~\cite{GraAlsKru17} for details):
\begin{definition}\label{def:well-ordered_symbolic-sequence}
Let $\x\in W_{p,q}$ be a periodic symbolic sequence.  Consider the
(lexicographically) ordered sequence given by the iterates of $\x$ by
$\sigma$
\begin{equation}
\sigma^{i_0}(\x)<\sigma^{i_1}(\x)<\sigma^{i_2}(\x)<\dots<\sigma^{i_{q-1}}(\x).
\label{eq:sequence_of_sequences}
\end{equation}
We say that the sequence $\x$ is a $p,q$-ordered (symbolic) sequence
if
\begin{equation*}
i_j-i_{j-1}=\text{constant}.
\end{equation*}
\end{definition}
In other words, a sequence $\x\in W_{p,q}$ is $p,q$-ordered if
$\sigma$ acts on the
sequence~\eqref{eq:sequence_of_sequences} as a cyclic permutation, i.e.
there exists some $k\in\N$, $0< k<q$, such that
\begin{equation}
i_j=i_{j-1}+k\pmod{q}.
\label{eq:maximin_cyclie-permutation}
\end{equation}
Then, from~\cite{GamLanTre84}, one has that a symbolic sequence
$\x\in W_{p,q}$ is maximin if, and only if, it is $p,q$-ordered.\\
Using this definition, checking computationally whether a symbolic
sequence is maximin becomes simpler and faster, as one does not need
to compare, element by element, the sequence with the rest of sequences in the set
$W_{p,q}$. The algorithm is as follows. Given a symbolic sequence
$\x\in W_{p,q}$ we consider the $q$ sequences given by
$\sigma^i(\x)$ with $0\le i <q$. To order them lexicographically, we
simply consider them as integer numbers written in binary. Thus, for
a symbolic sequence $\x=(x_0,\ldots,x_{q-1})$ we compute
\begin{equation*}
a(\x)=\sum_{i=0}^{q-1}x_i2^i.
\end{equation*}
Let us define $a_i:=a(\sigma^i(\x))$. Then, note that $\sigma^i(\x)<\sigma^j(\x)$ if, and only if, $a_i<a_j$.
Hence we easily get a sequence
\begin{equation*}
a_{i_0}<a_{i_1}<\cdots<a_{i_{q-1}}.
\end{equation*}
Then, if condition~\eqref{eq:maximin_cyclie-permutation} is satisfied
for all indices $i_j$, $0\le j<q$ the sequence $\x$ is maximin.


\section*{Figures}

\unitlength=\textwidth

\unitlength=\textwidth
\begin{figure}
\begin{center}
\begin{picture}(1,1)
\put(0,1.12){
\subfigure[\label{fig:Sigma1trans}]{\includegraphics[angle=-90,width=0.5\textwidth]{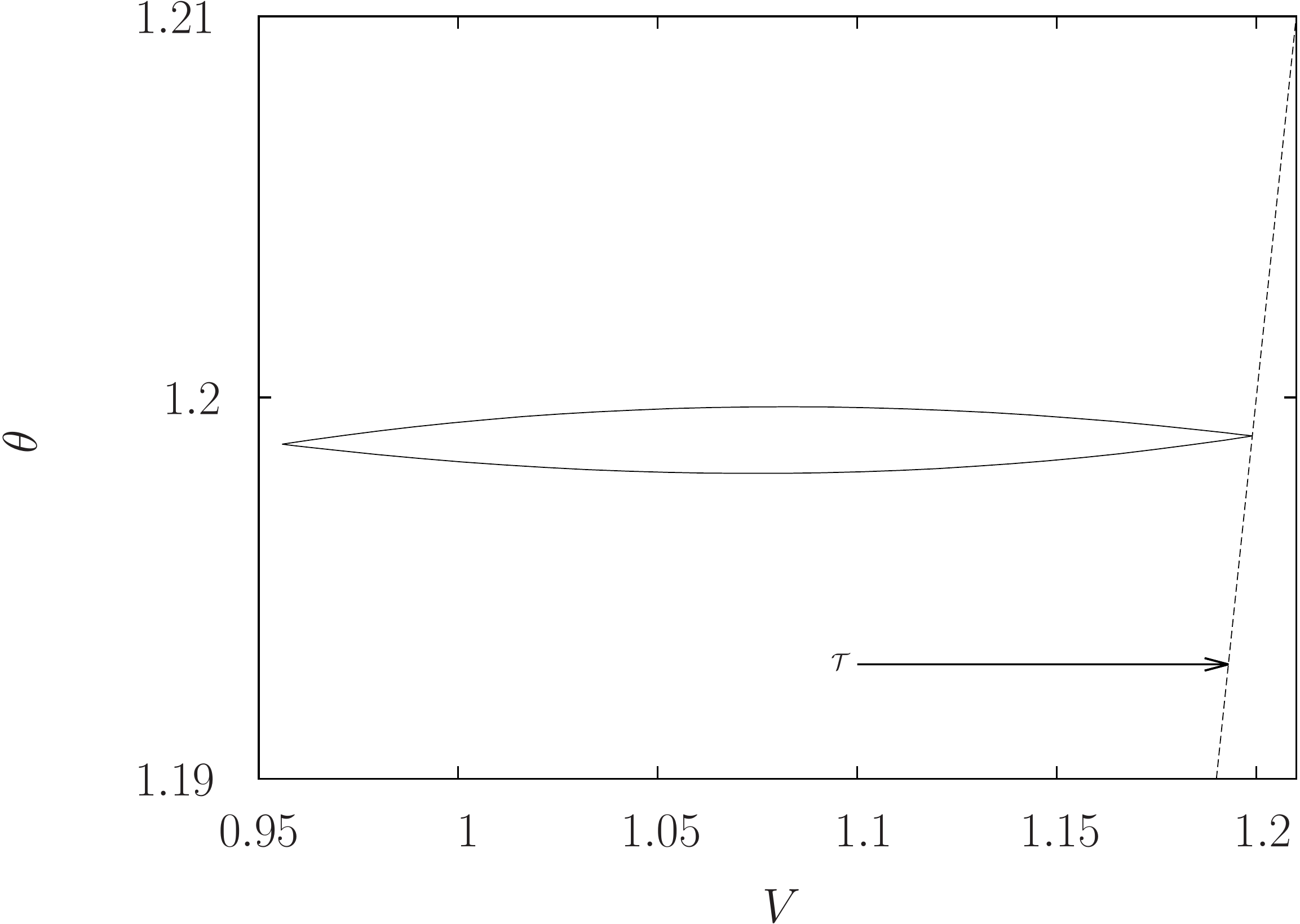}}
}
\put(0.5,1.12){
\subfigure[\label{fig:Sigma1trans_TS}]{\includegraphics[angle=-90,width=0.5\textwidth]{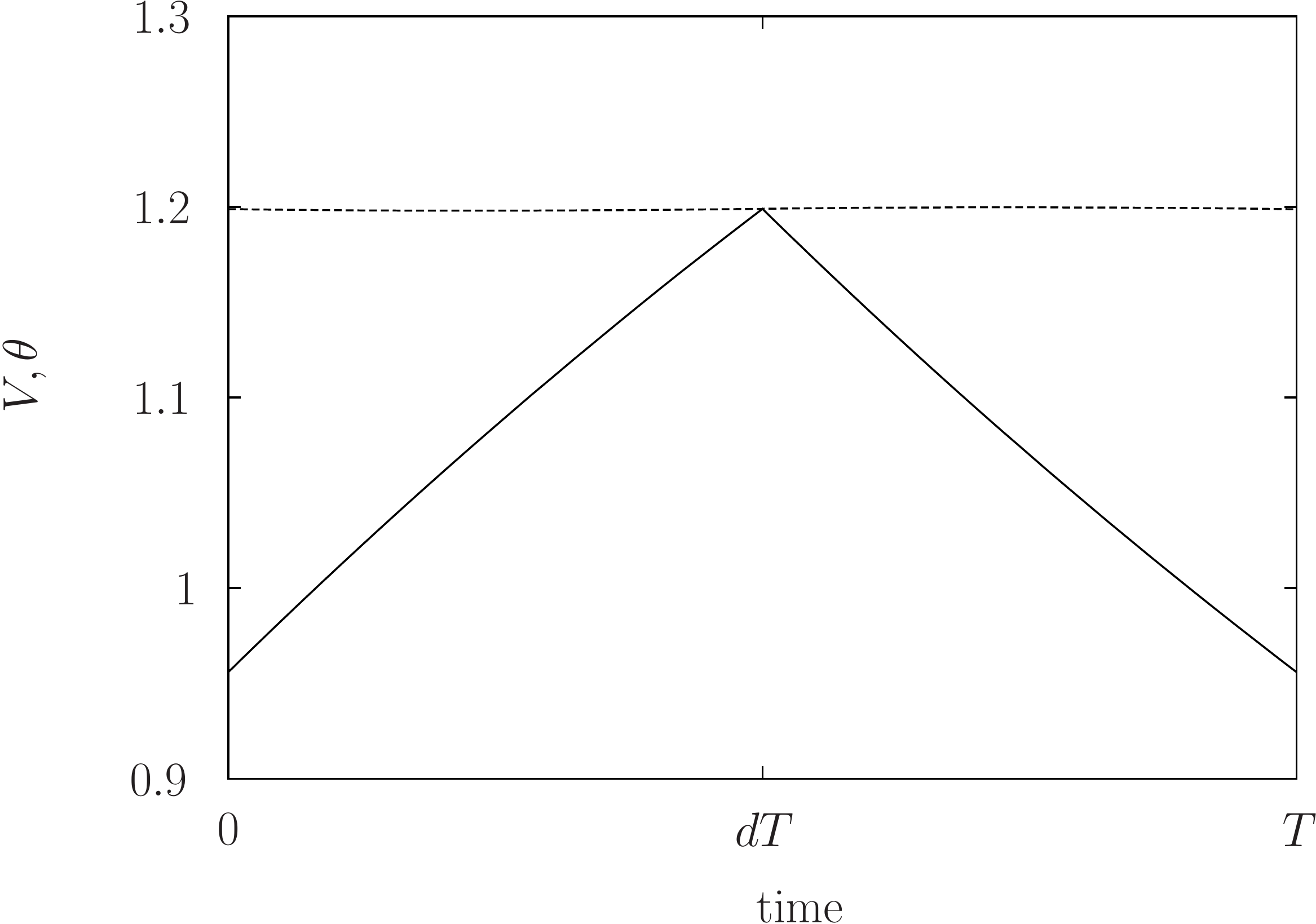}}
}
\put(0,0.72){
\subfigure[\label{fig:Sigma1transz1}]{\includegraphics[angle=-90,width=0.5\textwidth]{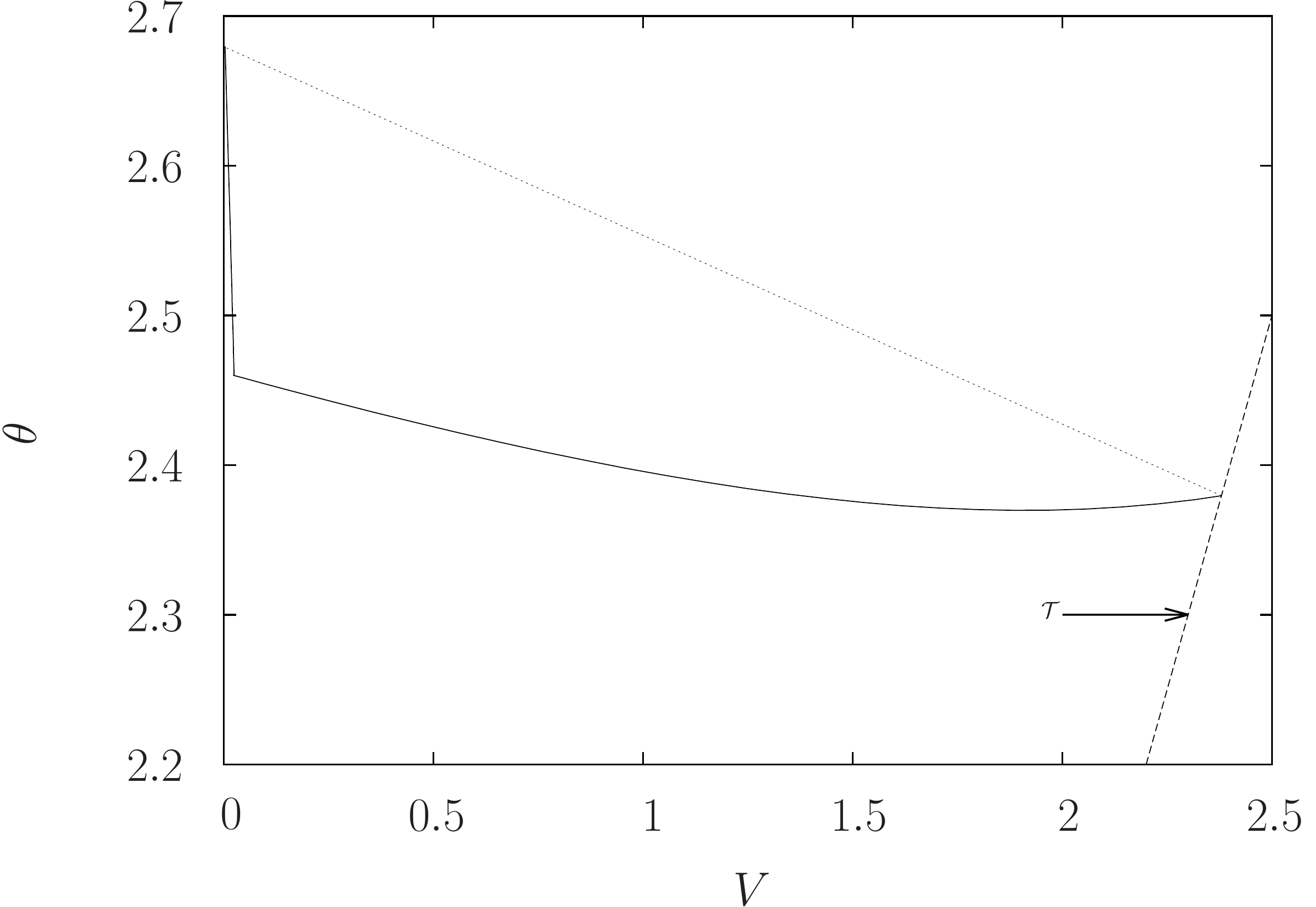}}
}
\put(0.5,0.72){
\subfigure[\label{fig:Sigma1transz1_TS}]{\includegraphics[angle=-90,width=0.5\textwidth]{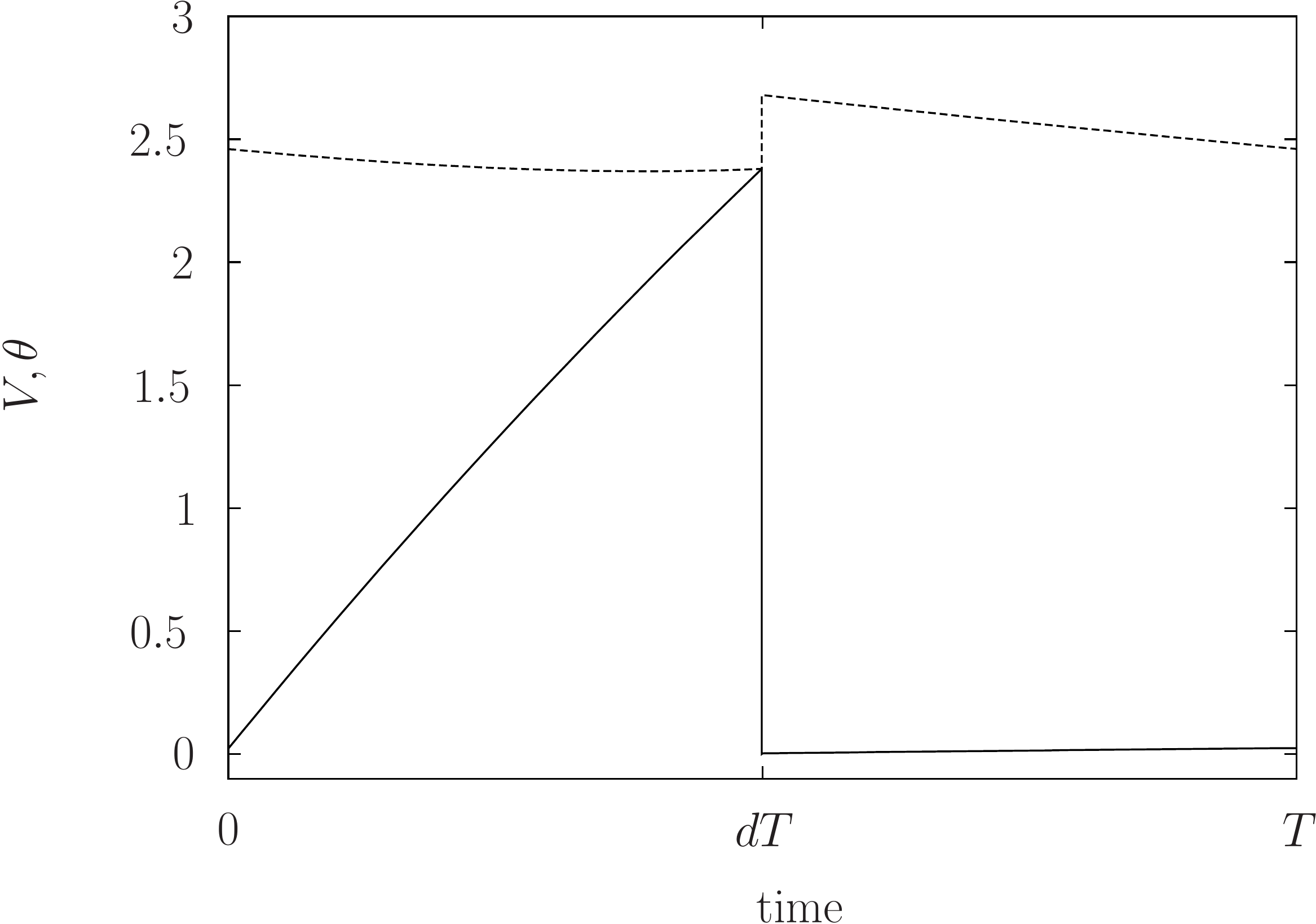}}
}
\put(0,0.32){
\subfigure[\label{fig:Sigma2transz1}]{\includegraphics[angle=-90,width=0.5\textwidth]{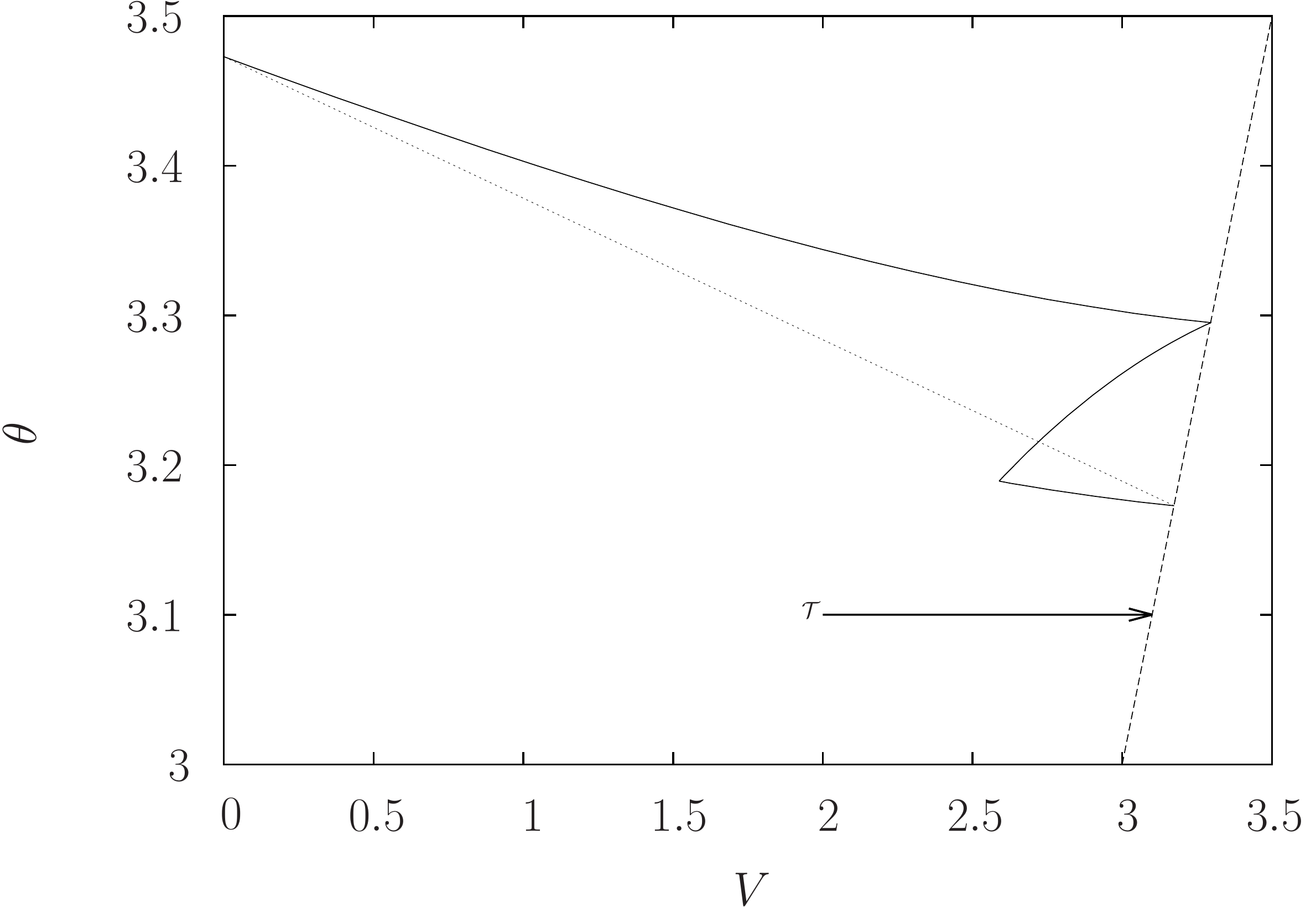}}
}
\put(0.5,0.32){
\subfigure[\label{fig:Sigma2transz1_TS}]{\includegraphics[angle=-90,width=0.5\textwidth]{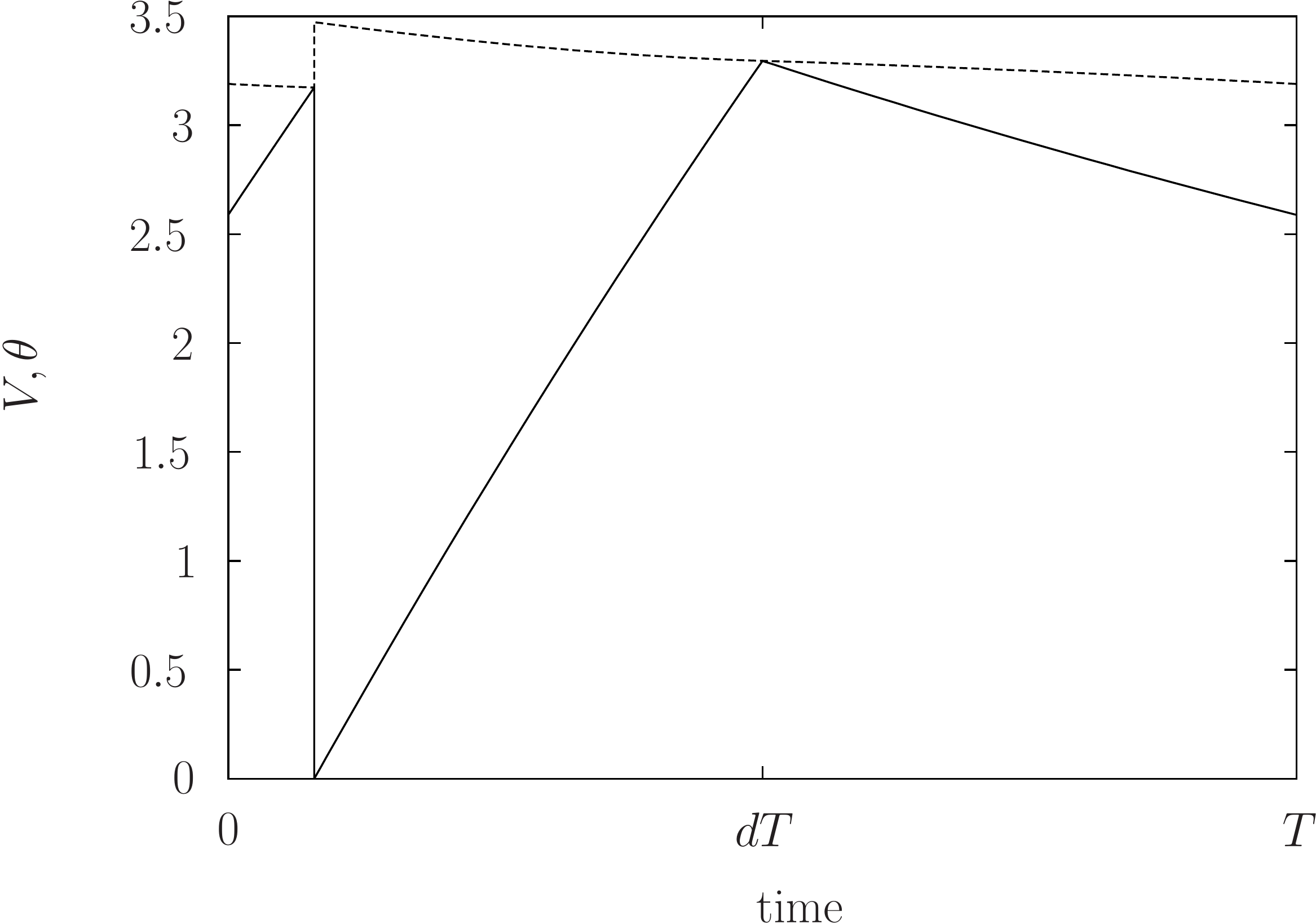}}
}
\end{picture}
\end{center}
\caption{Trajectories of $T$-periodic orbits undergoing the grazing bifurcations labeled in Figure~\ref{fig:2dbif_T0d5}.  
(a)-(b) Nonsmooth grazing bifurcation of a non-spiking periodic orbit
(border collision of $\bz_0\in S_0$ with $\Stran_1$). (c)-(d) Nonsmooth
grazing bifurcation of a 1-spiking periodic orbit (border collision of
$\bz_1\in S_1$ with $\Stran_1$). (e)-(f) Nonsmooth grazing
bifurcation of a 1-spiking periodic orbit (border collision of
$\bz_1\in S_1$ with $\Stran_2$).
Left panels show trajectories on the $(V,\theta)$-phase space while right panels show the corresponding time courses of the variables 
$V$ (solid line) and $\theta$ (dashed line) over 1 period.
}
\label{fig:grazing_T0d5}
\end{figure}

\begin{figure}
\begin{center}
\begin{picture}(1,0.7)
\put(0,0.7){
\subfigure[\label{fig:C21}]{\includegraphics[angle=-90,width=0.5\textwidth]{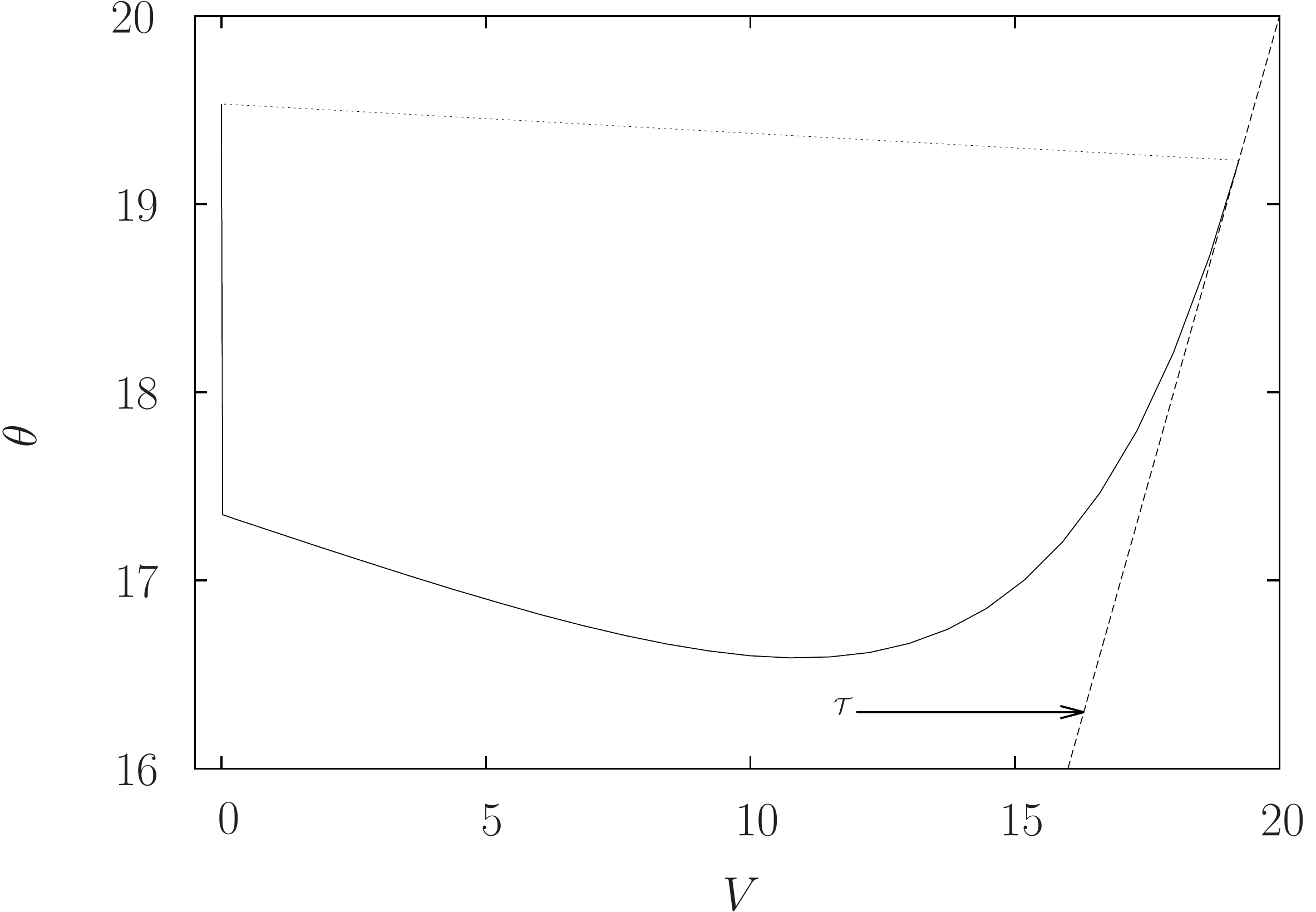}}
}
\put(0.5,0.7){
\subfigure[\label{fig:C21_TS}]{\includegraphics[angle=-90,width=0.5\textwidth]{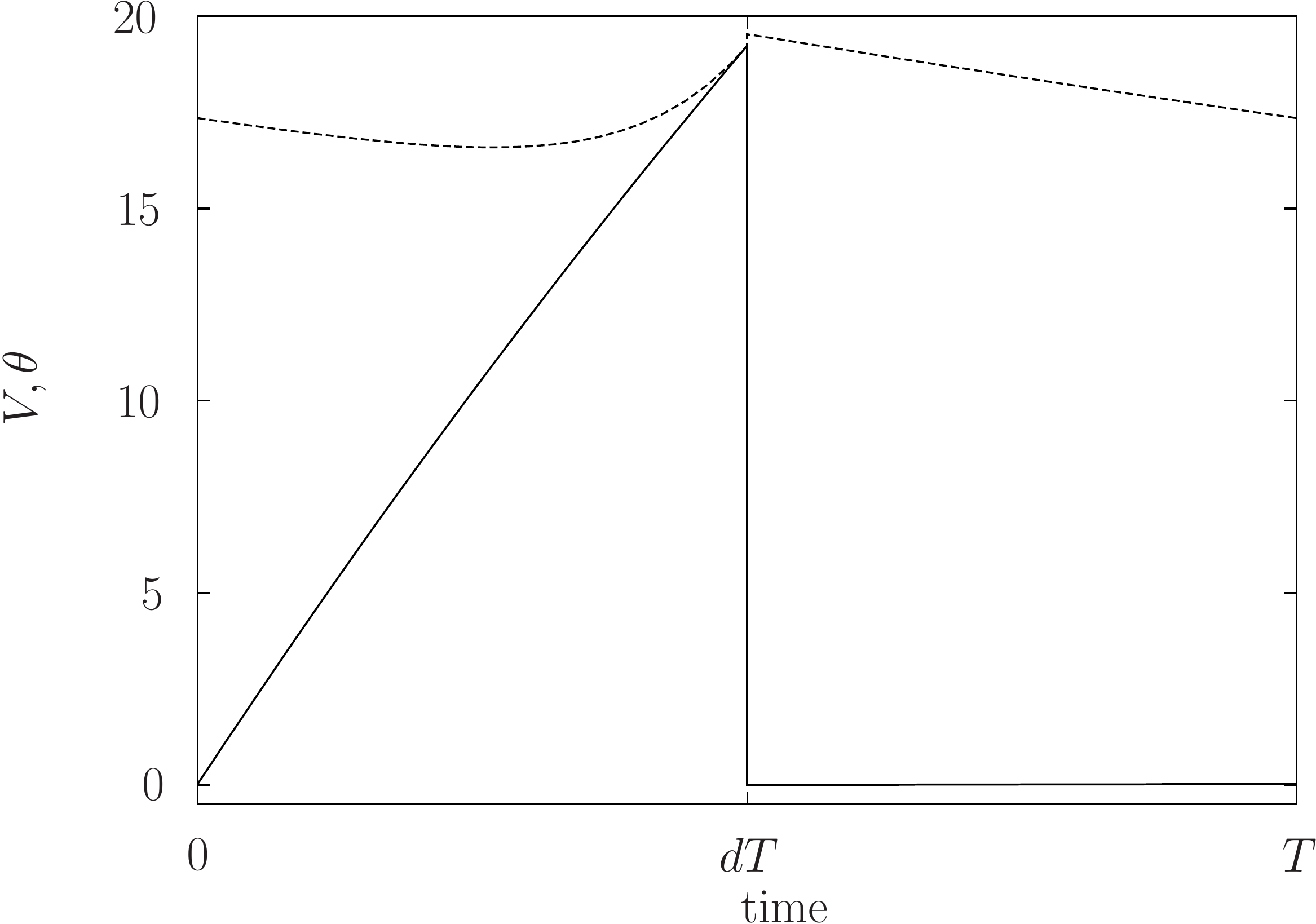}}
}
\put(0,0.32){
\subfigure[\label{fig:Sigma1tansz1}]{\includegraphics[angle=-90,width=0.5\textwidth]{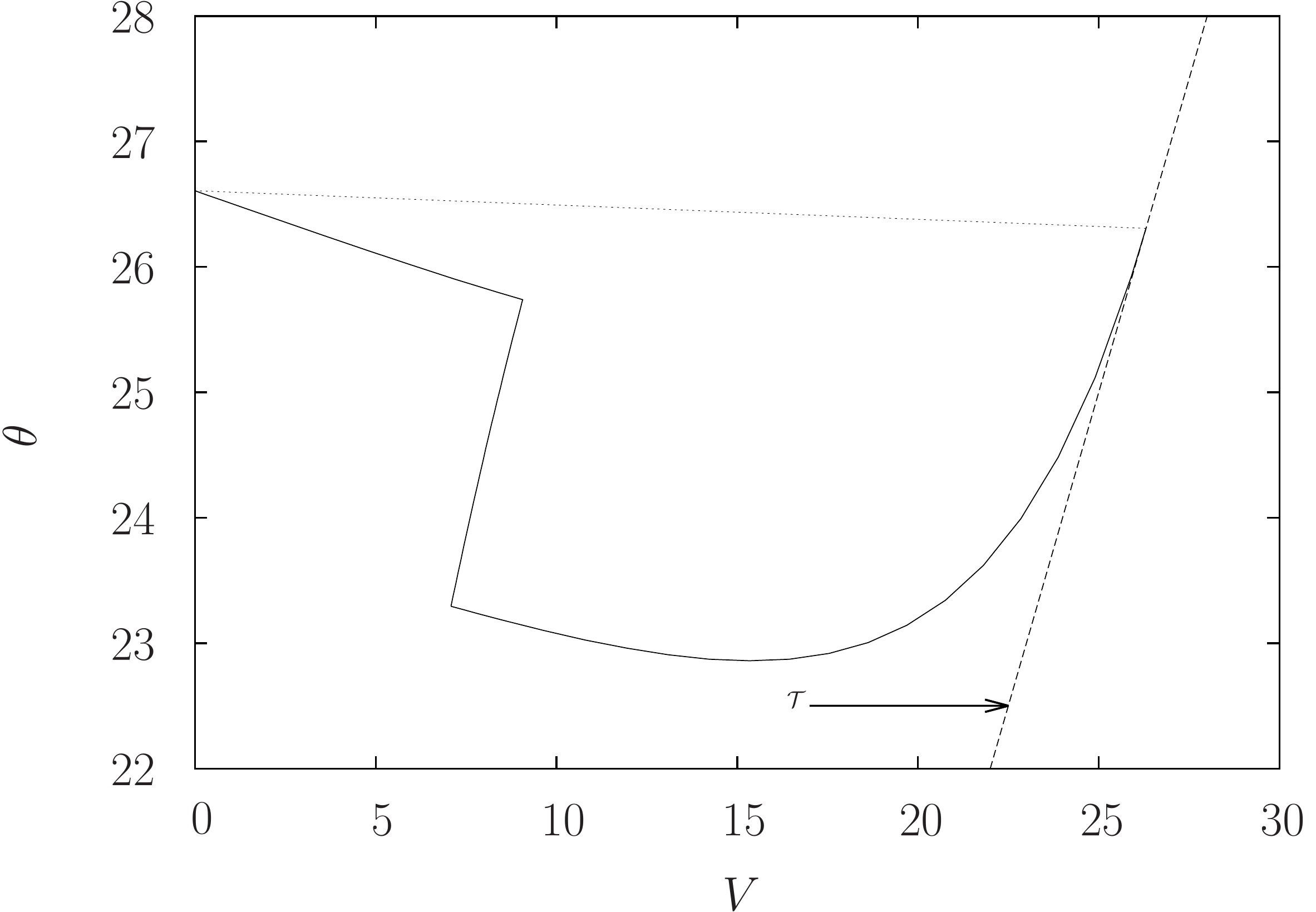}}
}
\put(0.5,0.32){
\subfigure[\label{fig:Sigma1tansz1_TS}]{\includegraphics[angle=-90,width=0.5\textwidth]{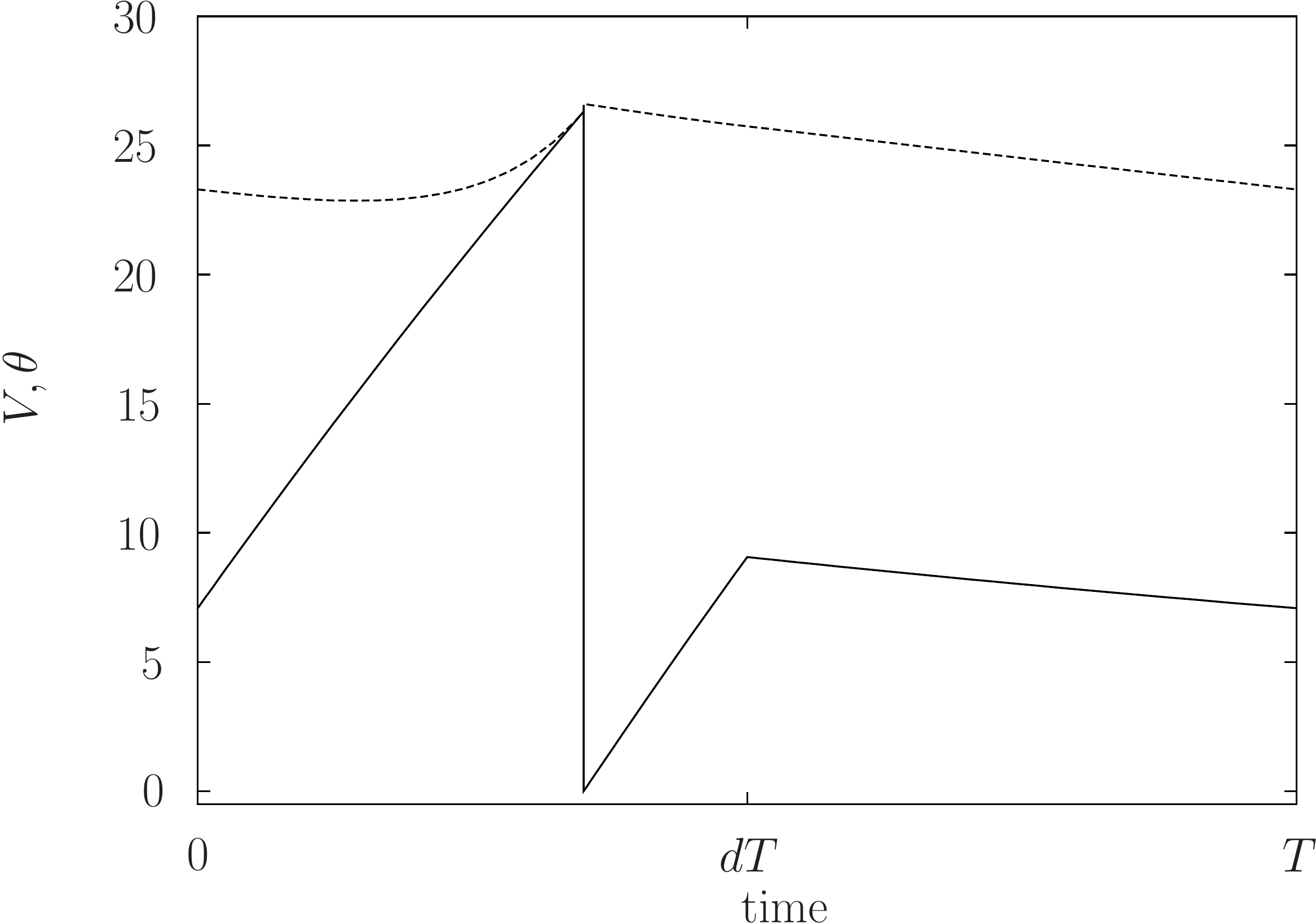}}
}
\end{picture}
\end{center}
\caption{Trajectories of $T$-periodic orbits undergoing the grazing bifurcations labeled in Figure~\ref{fig:2dbif_T0d5}.
(a)-(b) Co-dimension two bifurcation labeled as $C2_1$ in Figure~\ref{fig:2dbif_T0d5}; a periodic orbit undergoes a 
smooth bifurcation at $t=dT$ ($\bz_1 \in S_1$ collides with $\Stran_1$ and $\Stan_1$ simultaneously).  
(c)-(d) Smooth grazing bifurcation of a 1-spiking periodic orbit (border collision of $\bz_1\in
S_1$ with $\Stan_1$.)
Left panels show trajectories on the $(V,\theta)$-phase space while right panels show the corresponding time courses of the variables 
$V$ (solid line) and $\theta$ (dashed line) over 1 period.
}\label{fig:grzing2_T0d5}
\end{figure}

\begin{figure}
\begin{center}
\begin{picture}(1,0.6)
\put(0,0.72){
\subfigure[\label{fig:z0Stan1}]{
\includegraphics[angle=-90,width=0.5\textwidth]{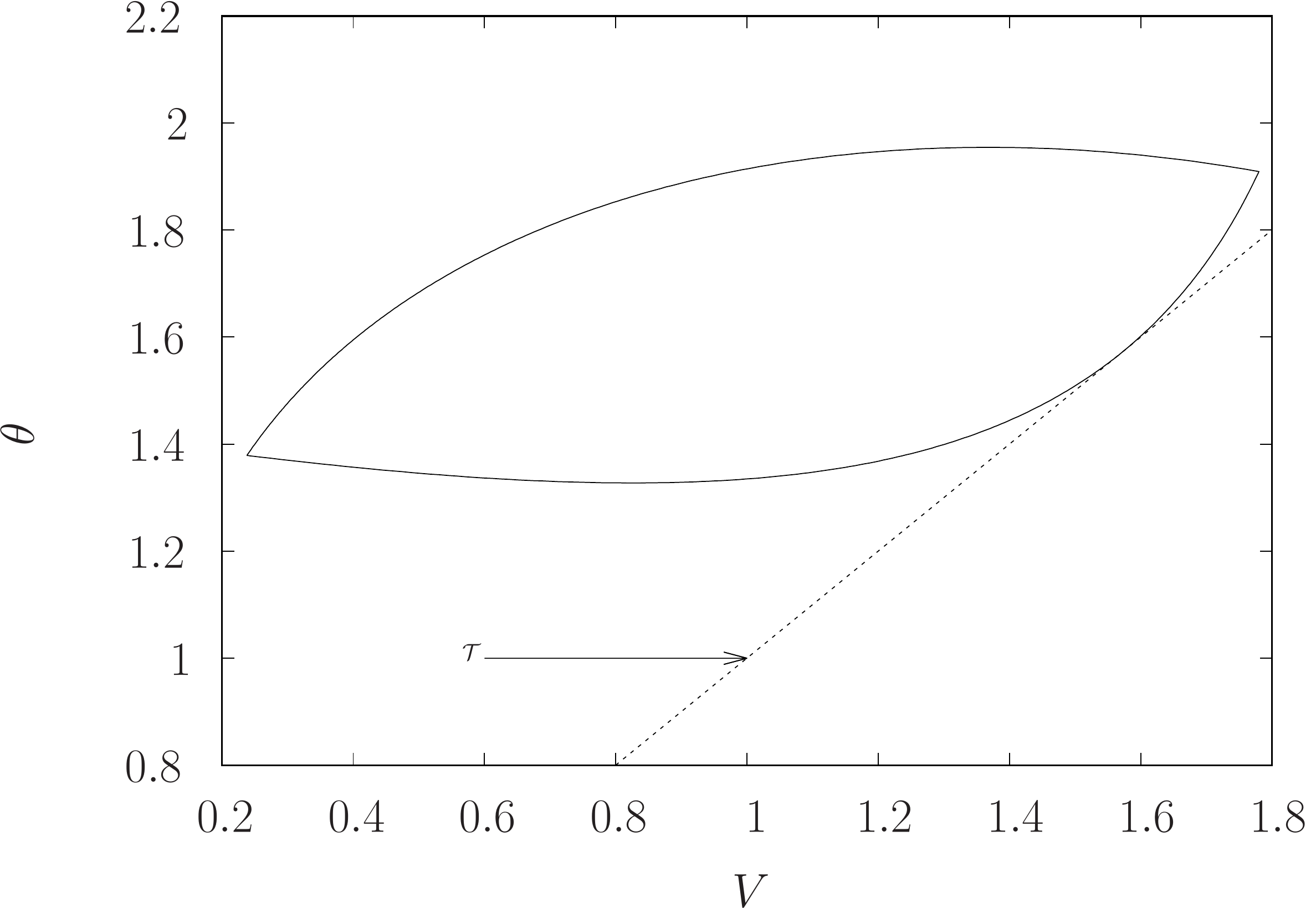}}
}
\put(0.5,0.72){
\subfigure[\label{fig:z0Stan1_TS}]{
\includegraphics[angle=-90,width=0.48\textwidth]{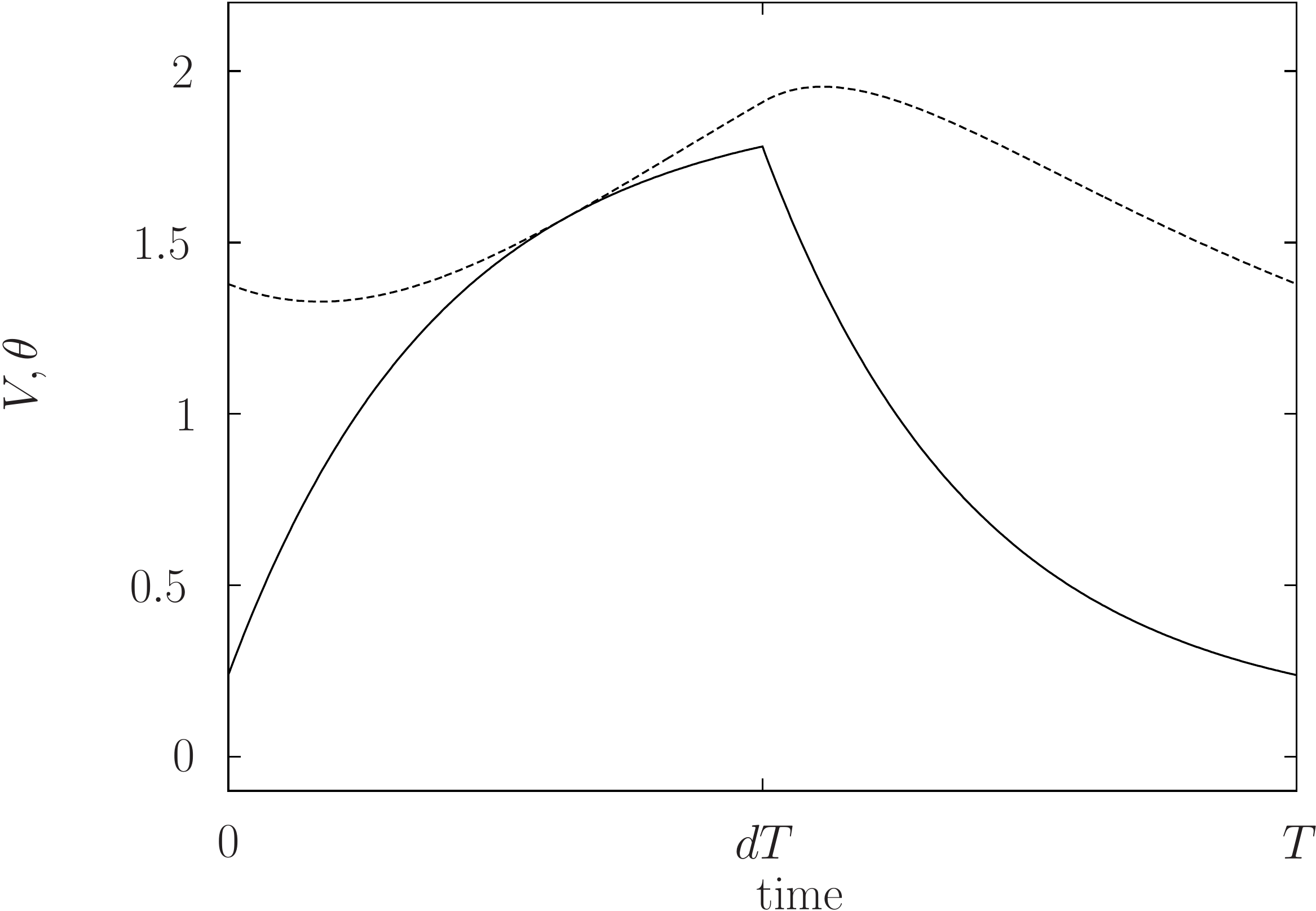}}
}
\put(0,0.32){
\subfigure[\label{fig:z1Stan2}]{
\includegraphics[angle=-90,width=0.5\textwidth]{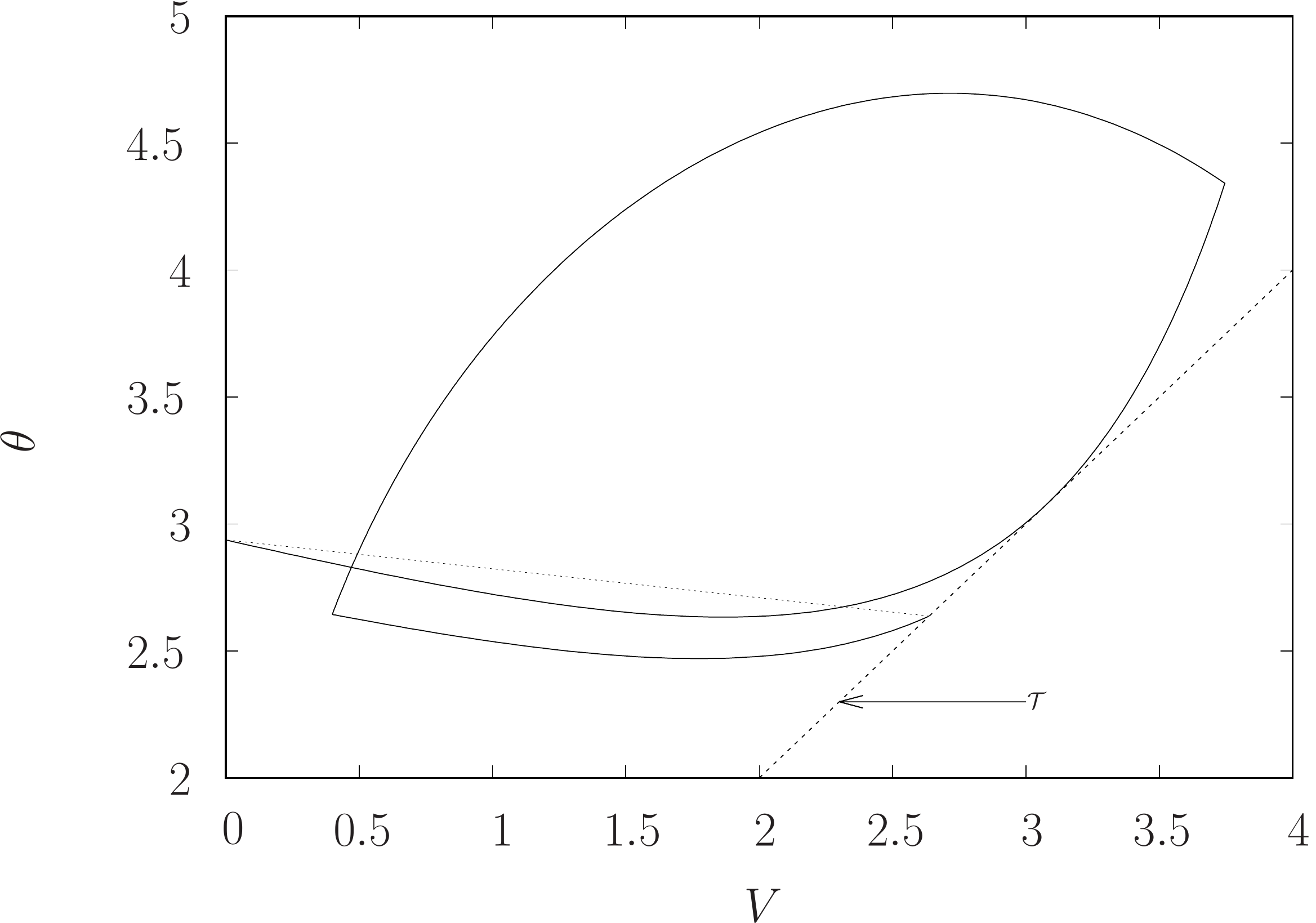}}
}
\put(0.5,0.32){
\subfigure[\label{fig:z1Stan2_TS}]{
\includegraphics[angle=-90,width=0.49\textwidth]{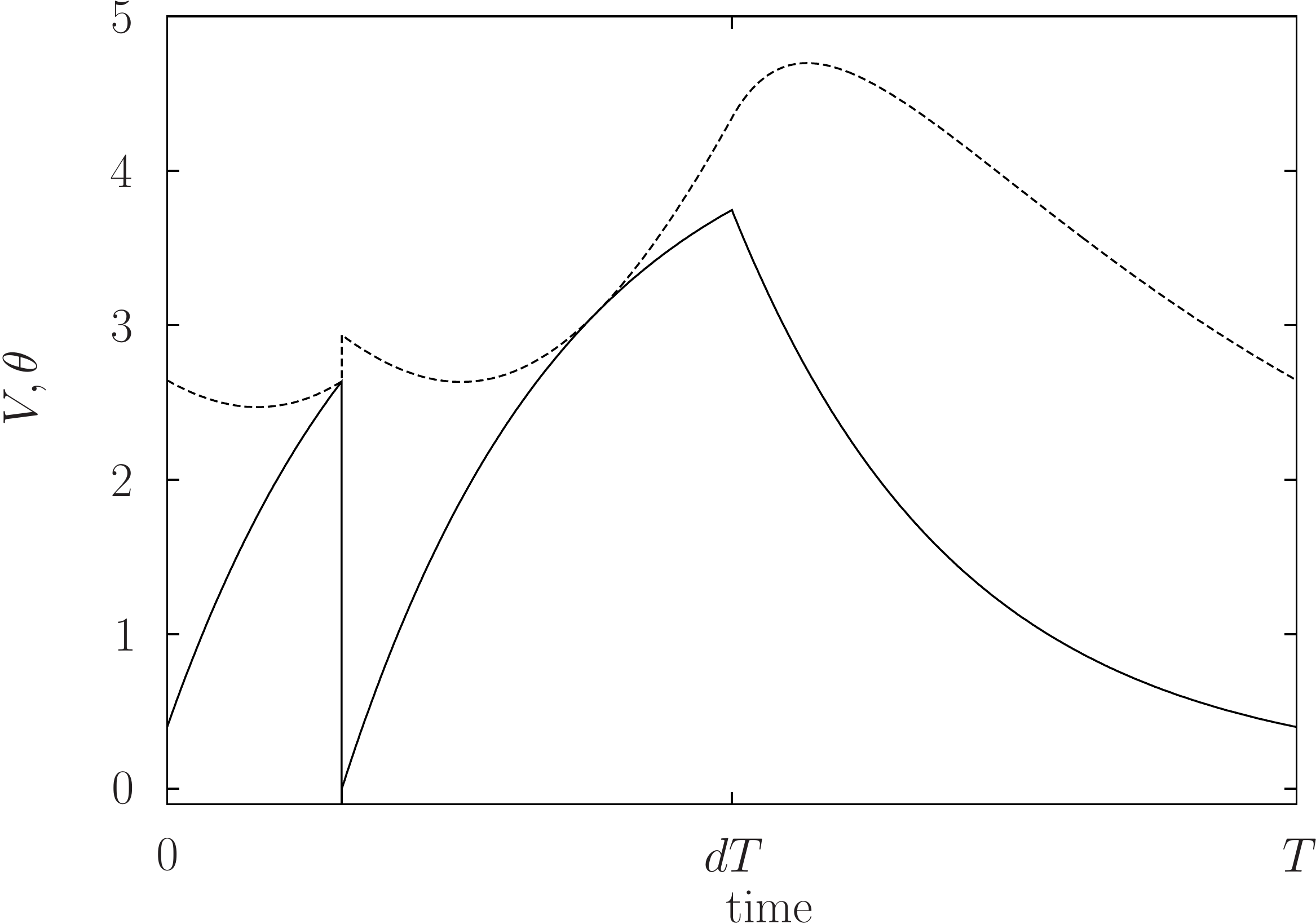}}
}
\end{picture}
\end{center}
\caption{Trajectories of $T$-periodic orbits undergoing the smooth grazing bifurcations labeled in
Figure~\ref{fig:2dbif_T5}. (a)-(b) Smooth grazing bifurcation of a
non-spiking periodic orbit (border collision of $\bz_0\in S_0$ with
$\Stan_1$). (c)-(d) Smooth grazing bifurcation of 1-spiking periodic orbit 
(border collision of $\bz_1\in S_1$ with $\Stan_2$).
}
\label{fig:tanbifs}
\end{figure}

\begin{figure}
\begin{center}
\begin{picture}(1,1)
\put(0,1.04){
\subfigure[\label{fig:z0Sigma1tantrans}]{
\includegraphics[angle=-90,width=0.45\textwidth]{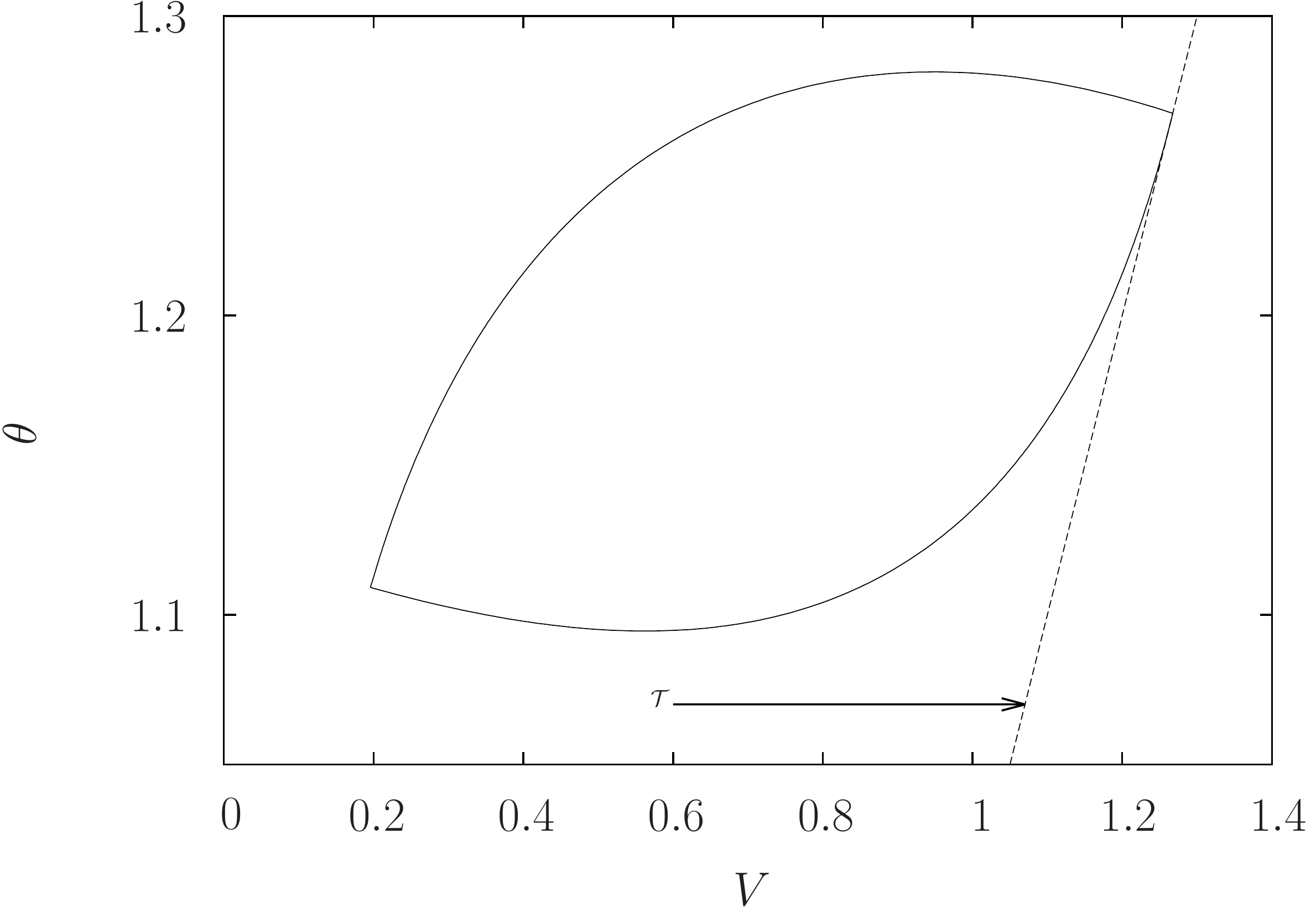}}
}
\put(0.5,1.04){
\subfigure[\label{fig:z0Sigma1tantrans_TS}]{
\includegraphics[angle=-90,width=0.45\textwidth]{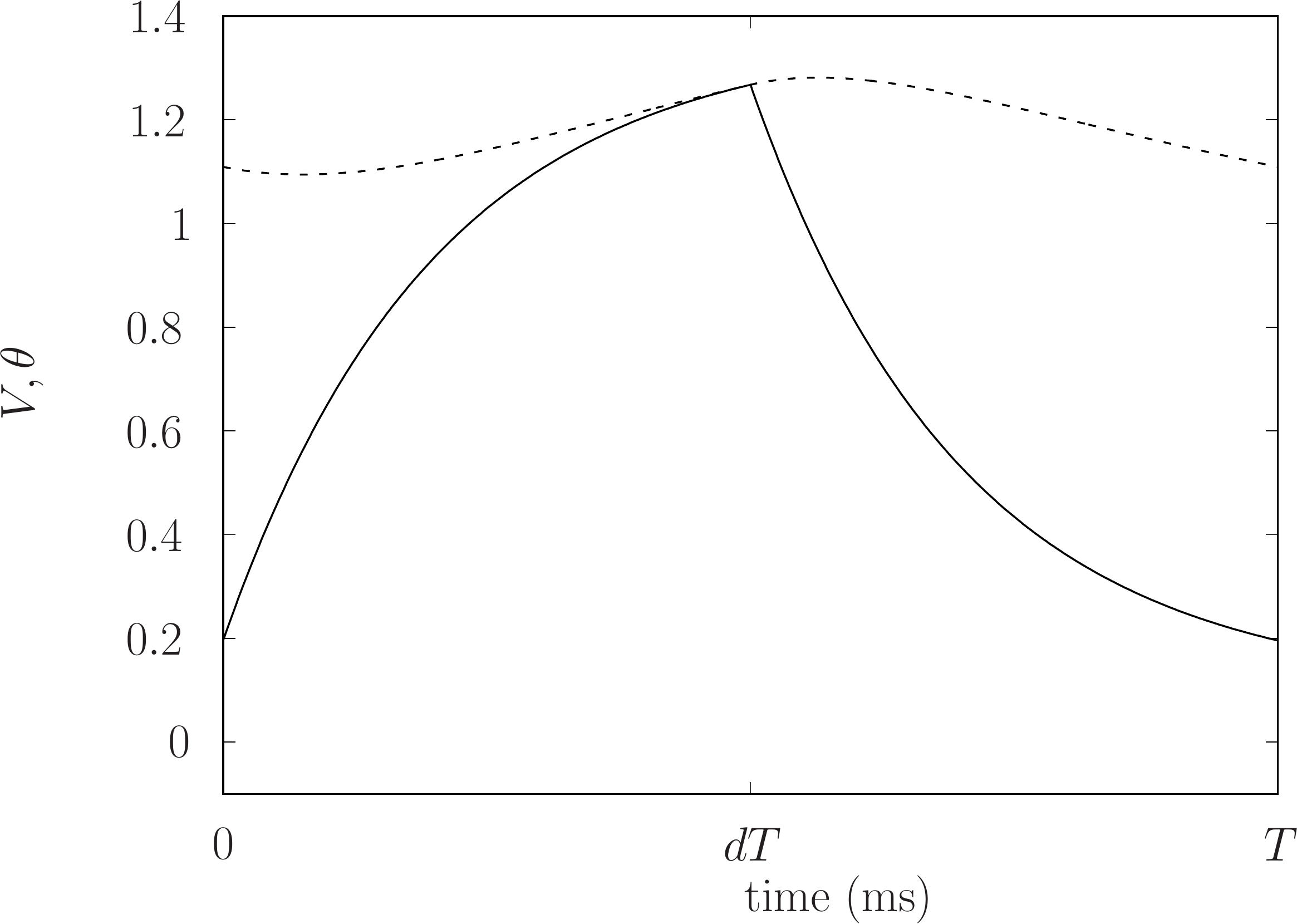}}
}
\put(0,0.68){
\subfigure[\label{fig:z1Stantrans2}]{\includegraphics[angle=-90,width=0.45\textwidth]{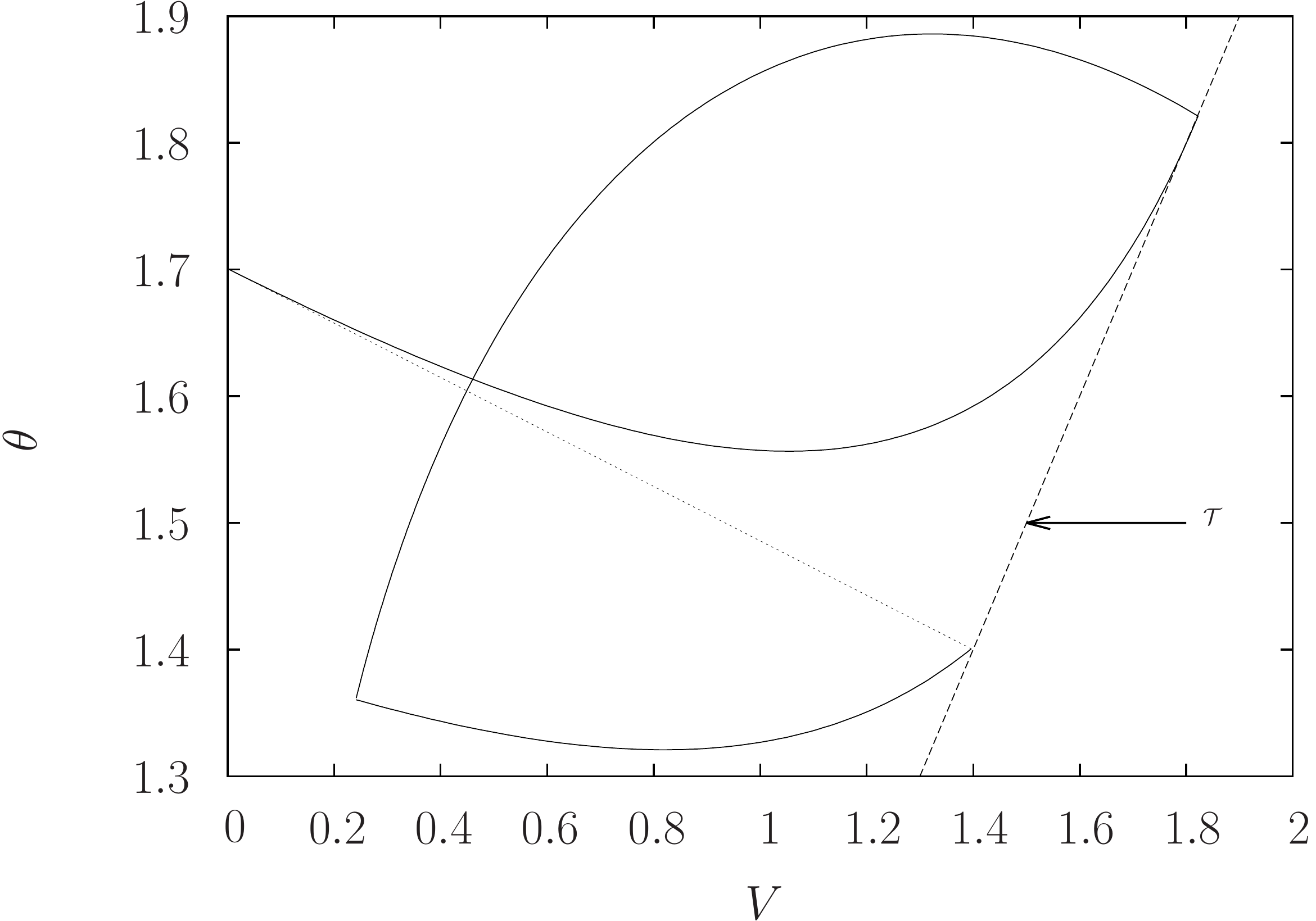}}
}
\put(0.5,0.68){
\subfigure[\label{fig:z1Stantrans2_TS}]{\includegraphics[angle=-90,width=0.45\textwidth]{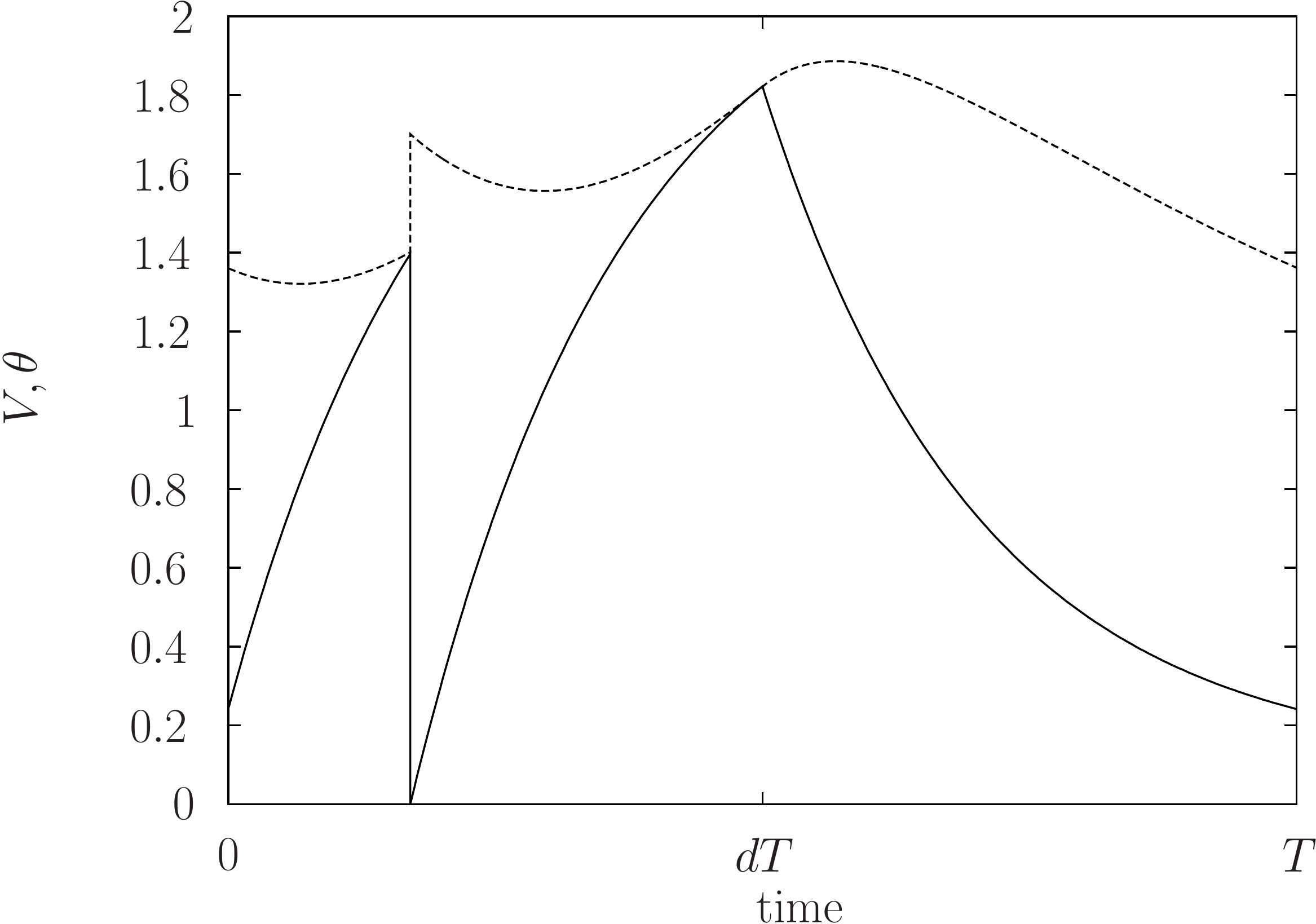}}
}
\put(0,0.32){
\subfigure[\label{fig:z1Sigma1tanSigma2tan}]{\includegraphics[angle=-90,width=0.45\textwidth]{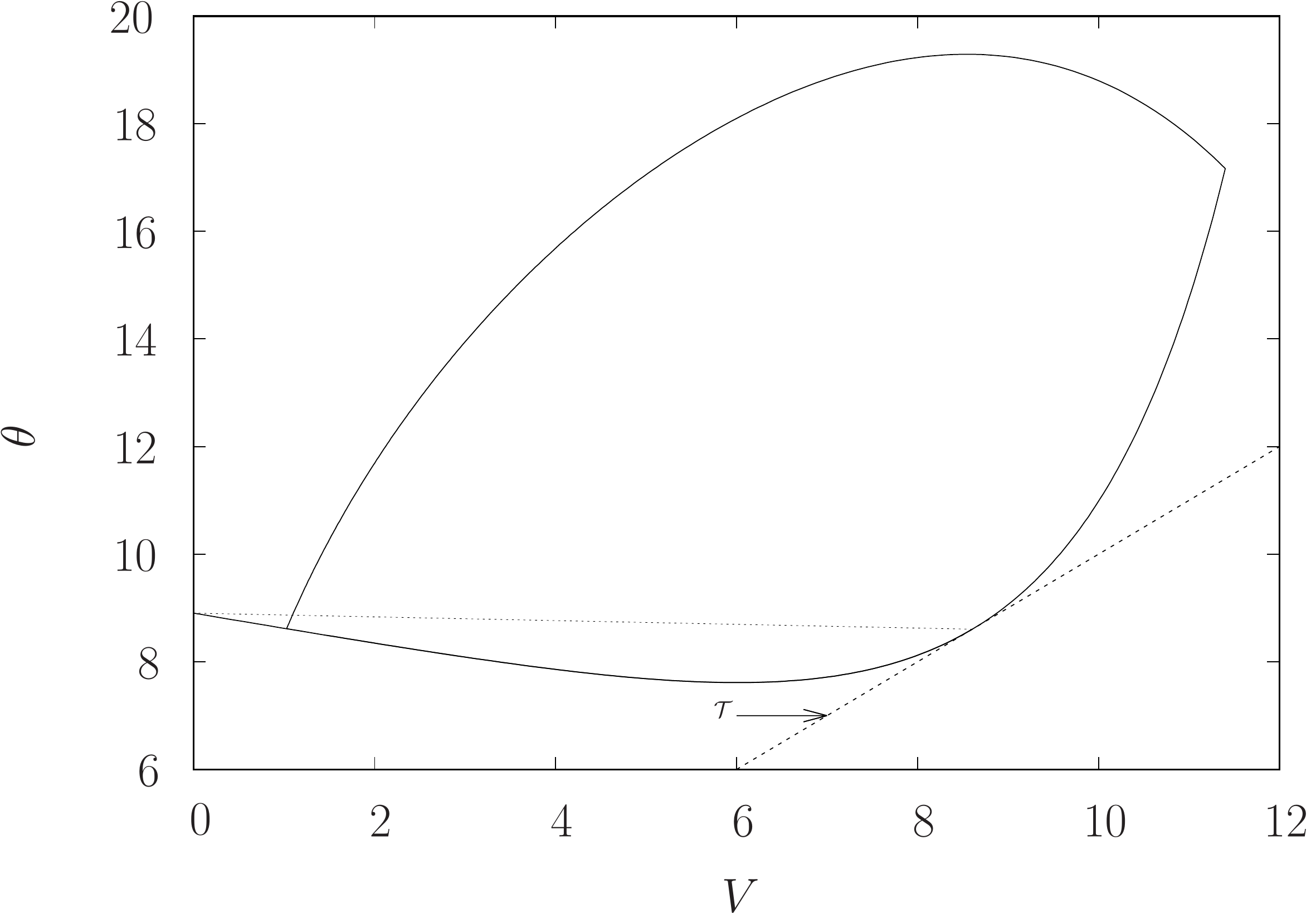}}
}
\put(0.5,0.32){
\subfigure[\label{fig:z1Sigma1tanSigma2tan_TS}]{\includegraphics[angle=-90,width=0.45\textwidth]{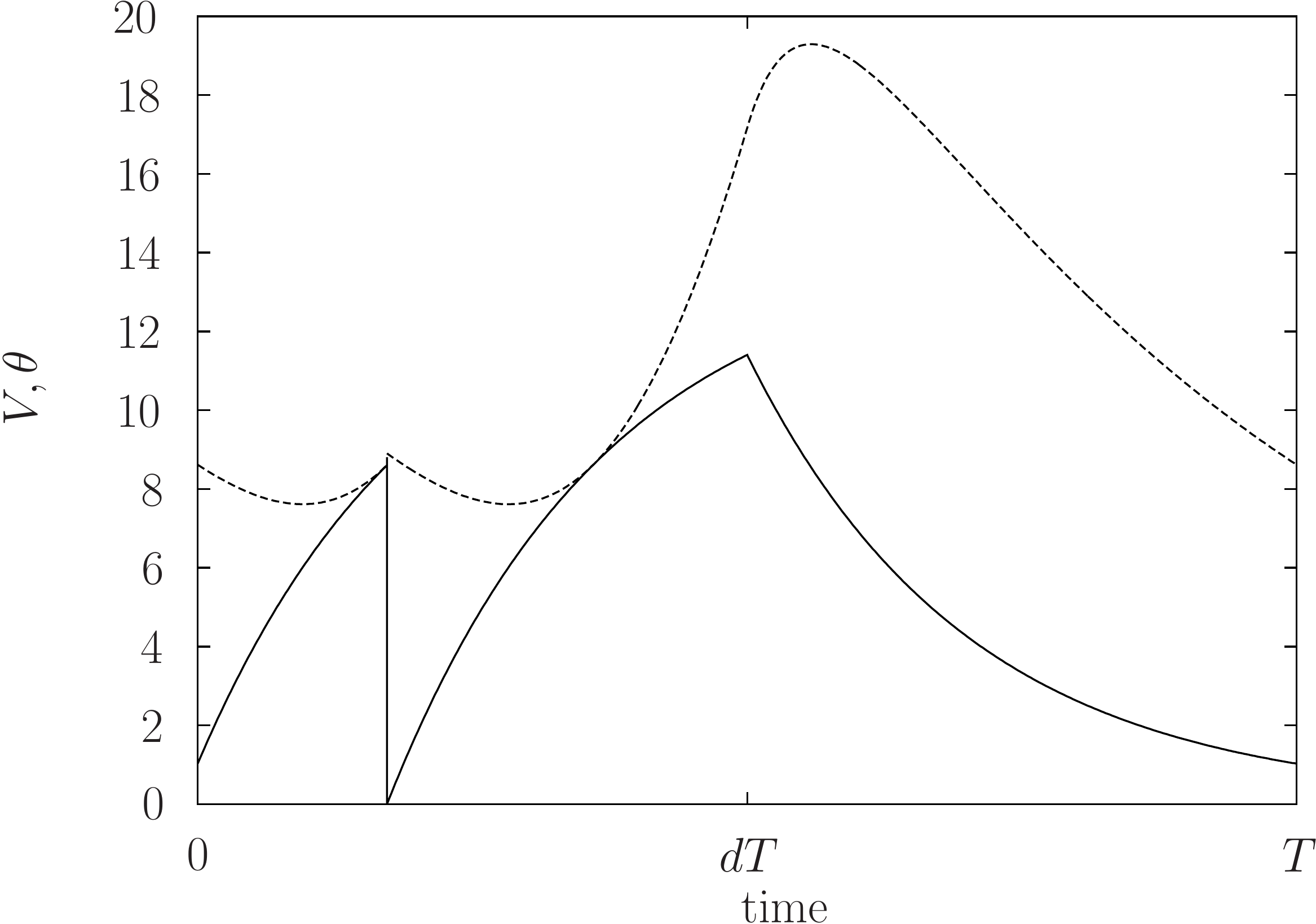}}
}
\end{picture}
\end{center}
\caption{Trajectories of the $T$-periodic orbits at the co-dimension two bifurcation
points labeled in Figure~\ref{fig:2dbif_T5}. (a)-(b) Point $C2_2$: a non-spiking periodic orbit
tangentially grazes the threshold at $t=dT$
($\bz_0\in\Stran_1\cap\Stan_1$). (c)-(d) Point $C2_3$: a 1-spiking periodic orbit
tangencially grazes the threshold precisely at $t=dT$, when the pulse is disabled
($\bz_1\in\Stran_2\cap\Stan_2$).  (e)-(f) Point
$C2_4$: a 1-spiking periodic orbit grazes the threshold twice, both
tangentially ($\bz_1\in\Stan_1\cap\Stan_2$).}
\label{fig:c23}
\end{figure}

\begin{figure}
\begin{center}
\begin{picture}(1,0.5)
\put(0,0.3){
\subfigure[\label{fig:2dbif_T05_numop_all}]{\includegraphics[angle=-90,width=0.45\textwidth]{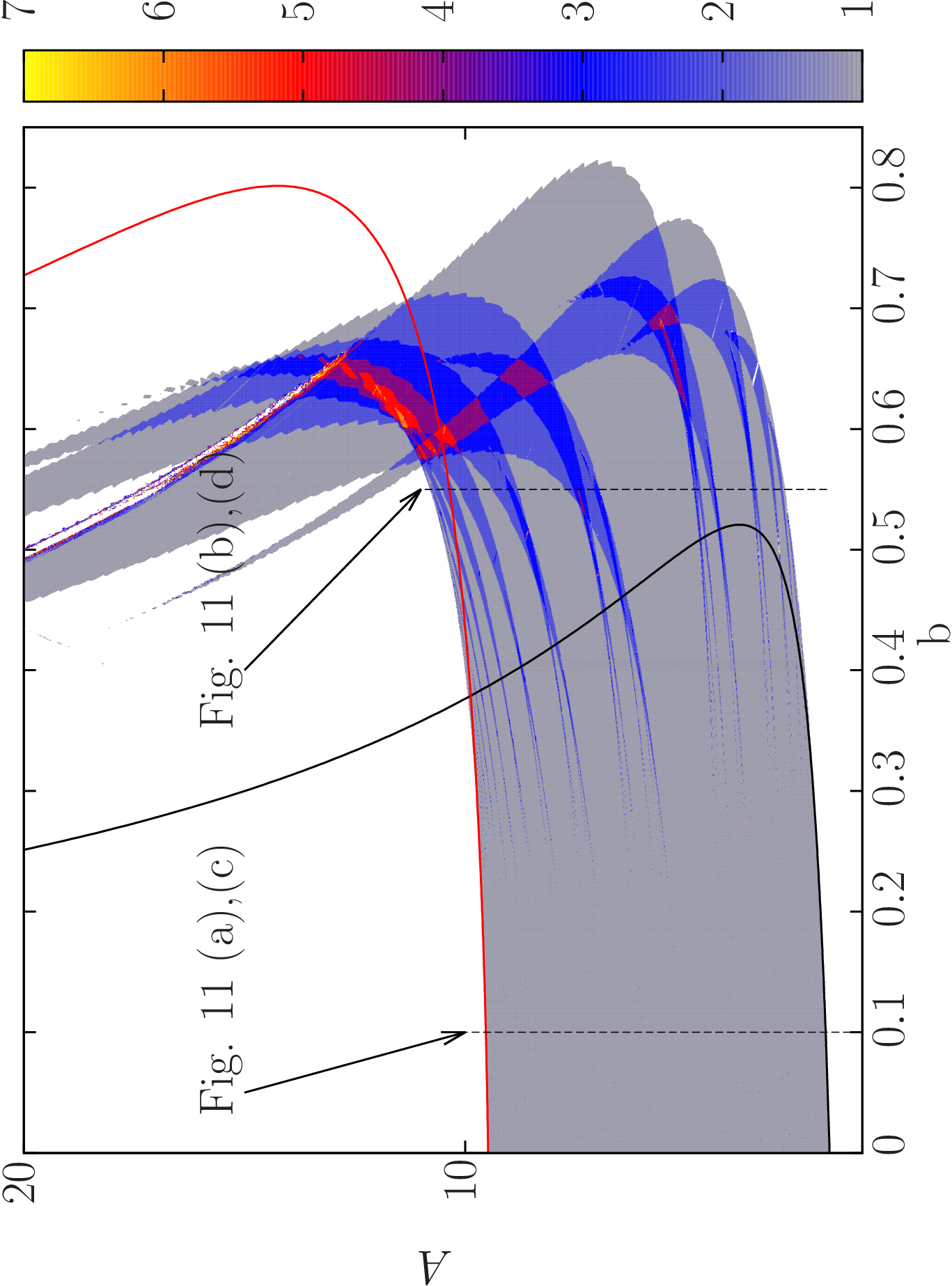}}
}
\put(0.5,0.3){
\subfigure[\label{fig:2dbif_T05_numop_period}]{\includegraphics[angle=-90,width=0.45\textwidth]{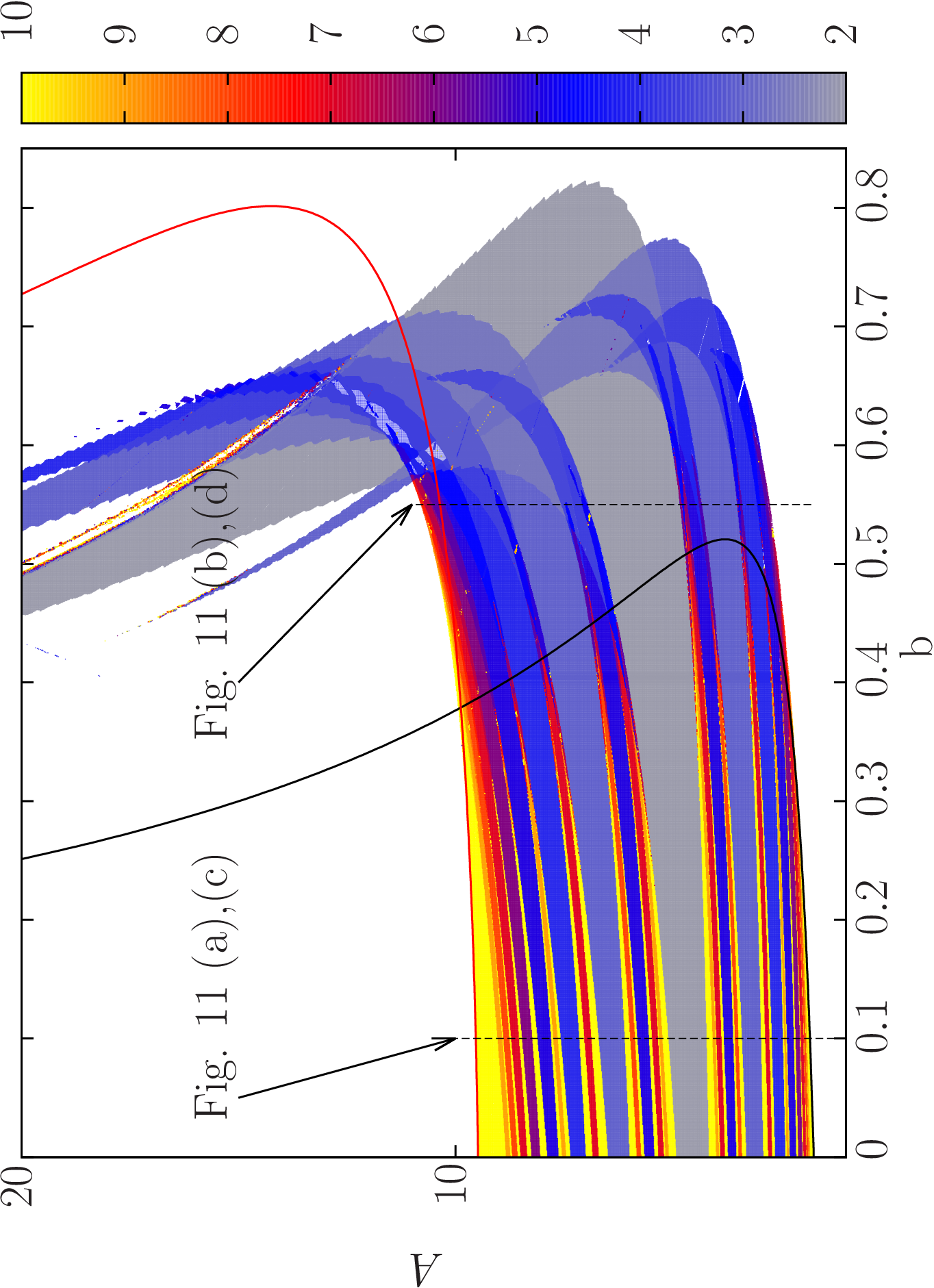}}
}
\end{picture}
\end{center}
\caption{(a) Number and (b) period of the periodic orbits of the stroboscopic map for $T=0.5$ stepping on $S_0$ and
$S_1$ found by means of the numerical algorithm described in Appendix~\ref{ap:per_orb}. 
We include also the bifurcation curves computed in Figure~\ref{fig:2dbif_T0d5}. 
All periodic orbits are maximin.
The areas where there is coexistence of two or more periodic orbits are colored according to an averaged period in order to 
reflect superposition of colors. Regions with orbits of period equal or higher than 10 have the same color.
Notice that the region where periodic orbits exist shows a jagged edge due to numerical issues related to grazing of the orbits.}
\label{fig:2dbif_po}
\end{figure}

\begin{figure}
\begin{center}
\begin{picture}(1,0.7)
\put(0,0.69){
\subfigure[\label{fig:periods_brut-force}]{\includegraphics[angle=-90,width=0.5\textwidth]{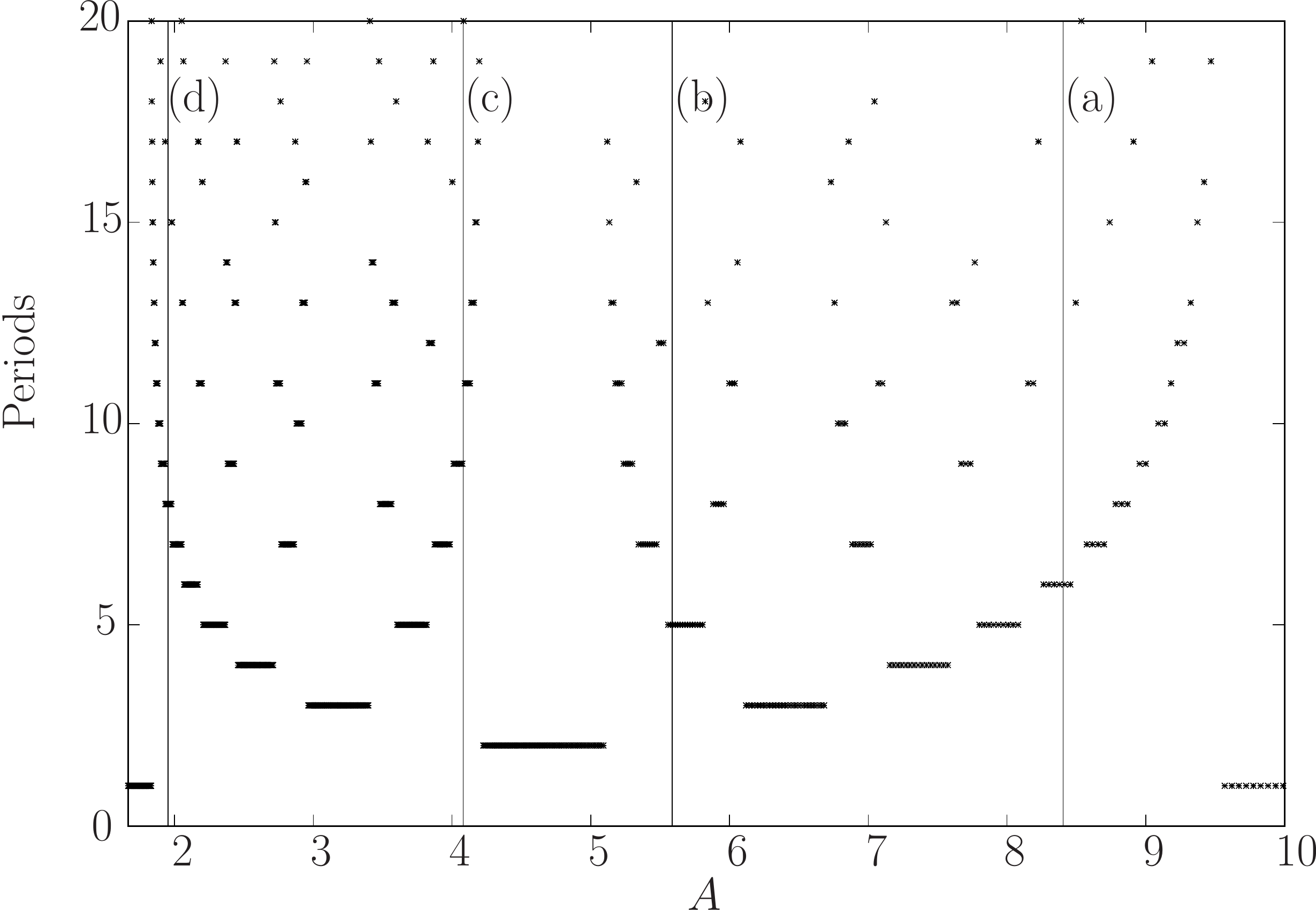}}
}
\put(0.5,0.69){
\subfigure[\label{fig:periods_brut-force2}]{\includegraphics[angle=-90,width=0.5\textwidth]{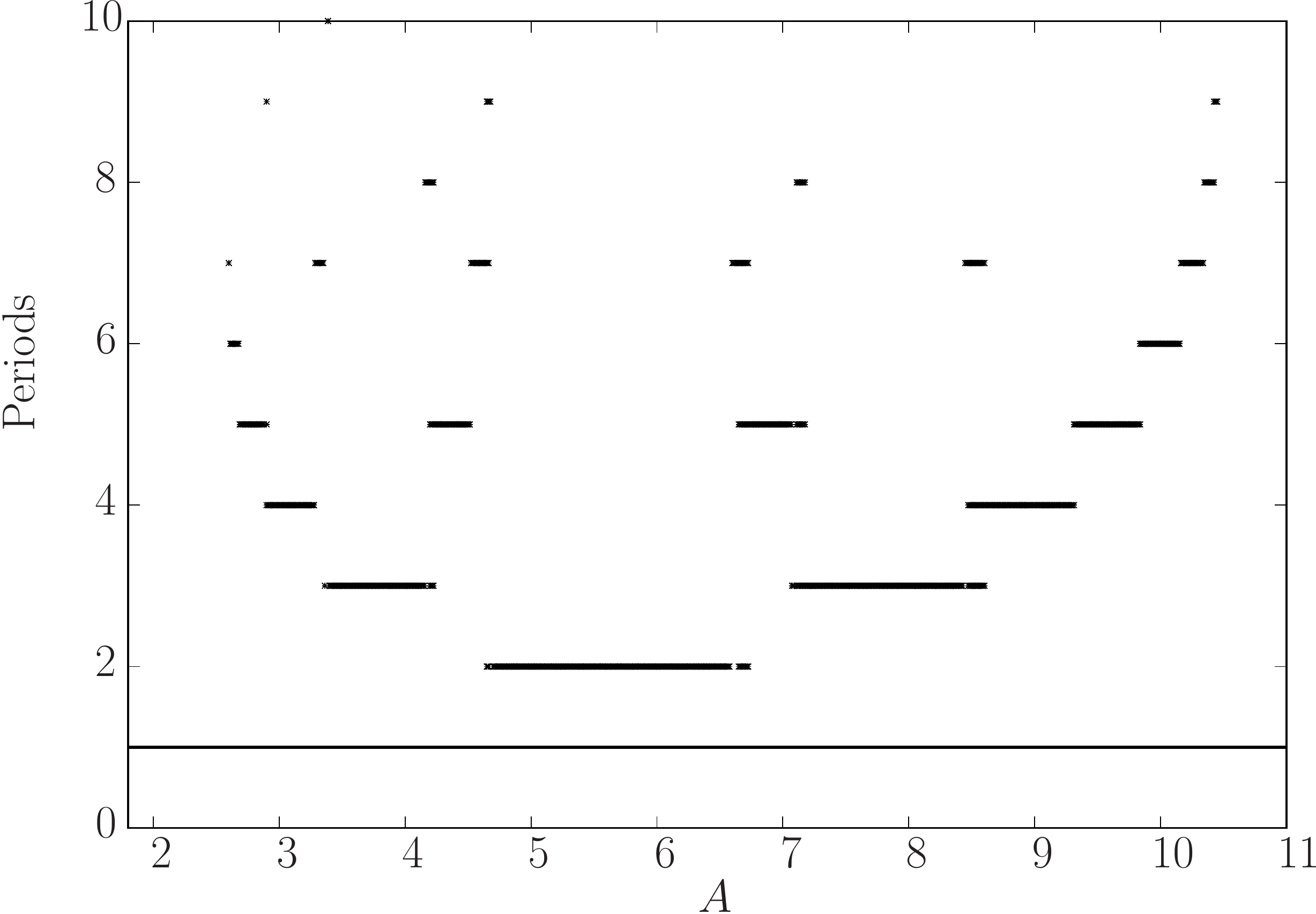}}
}
\put(0,0.31){
\subfigure[\label{fig:rotnums_b0d1}]{\includegraphics[angle=-90,width=0.5\textwidth]{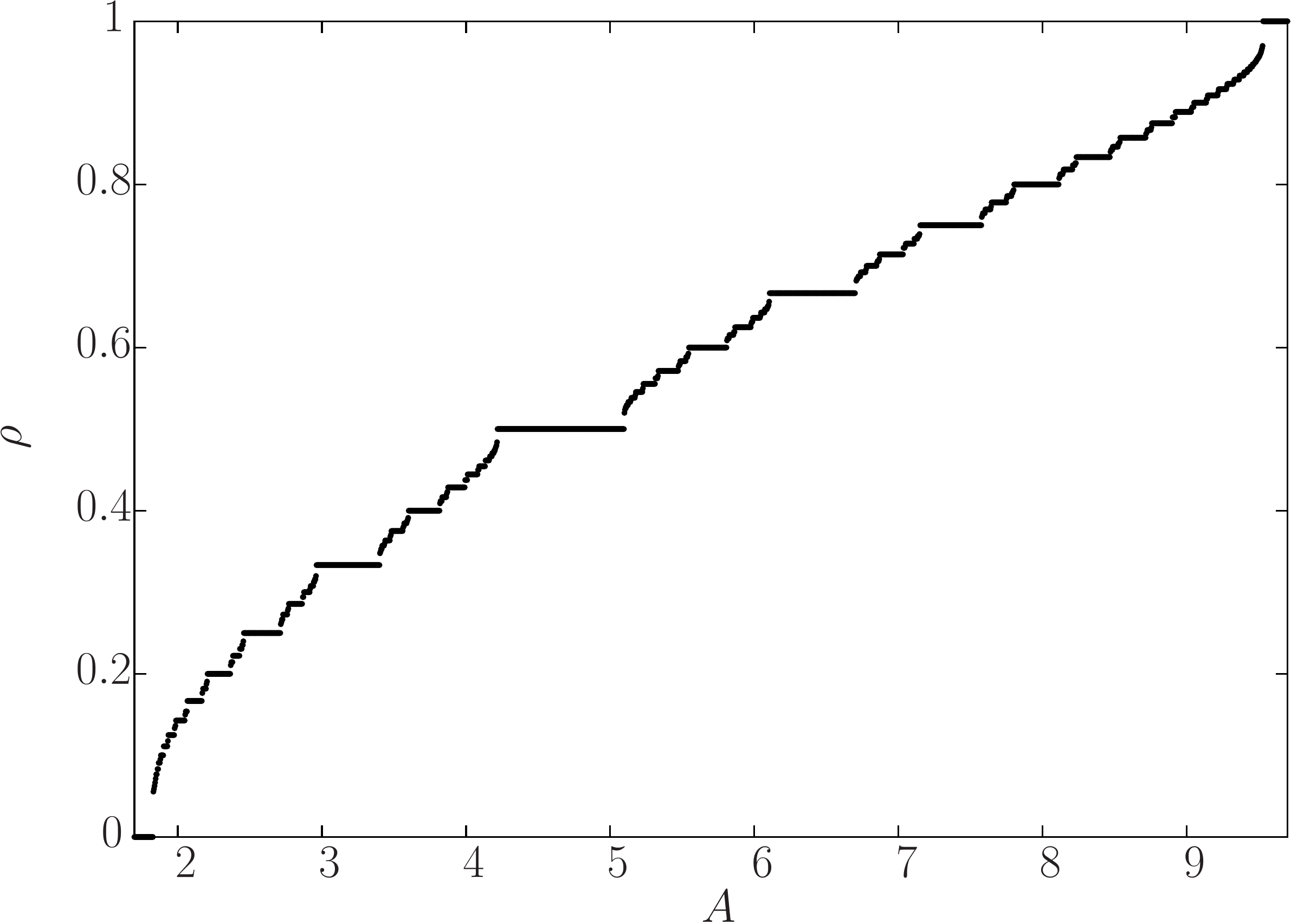}}
}
\put(0.5,0.31){
\subfigure[\label{fig:rotnums_b0d55}]{\includegraphics[angle=-90,width=0.5\textwidth]{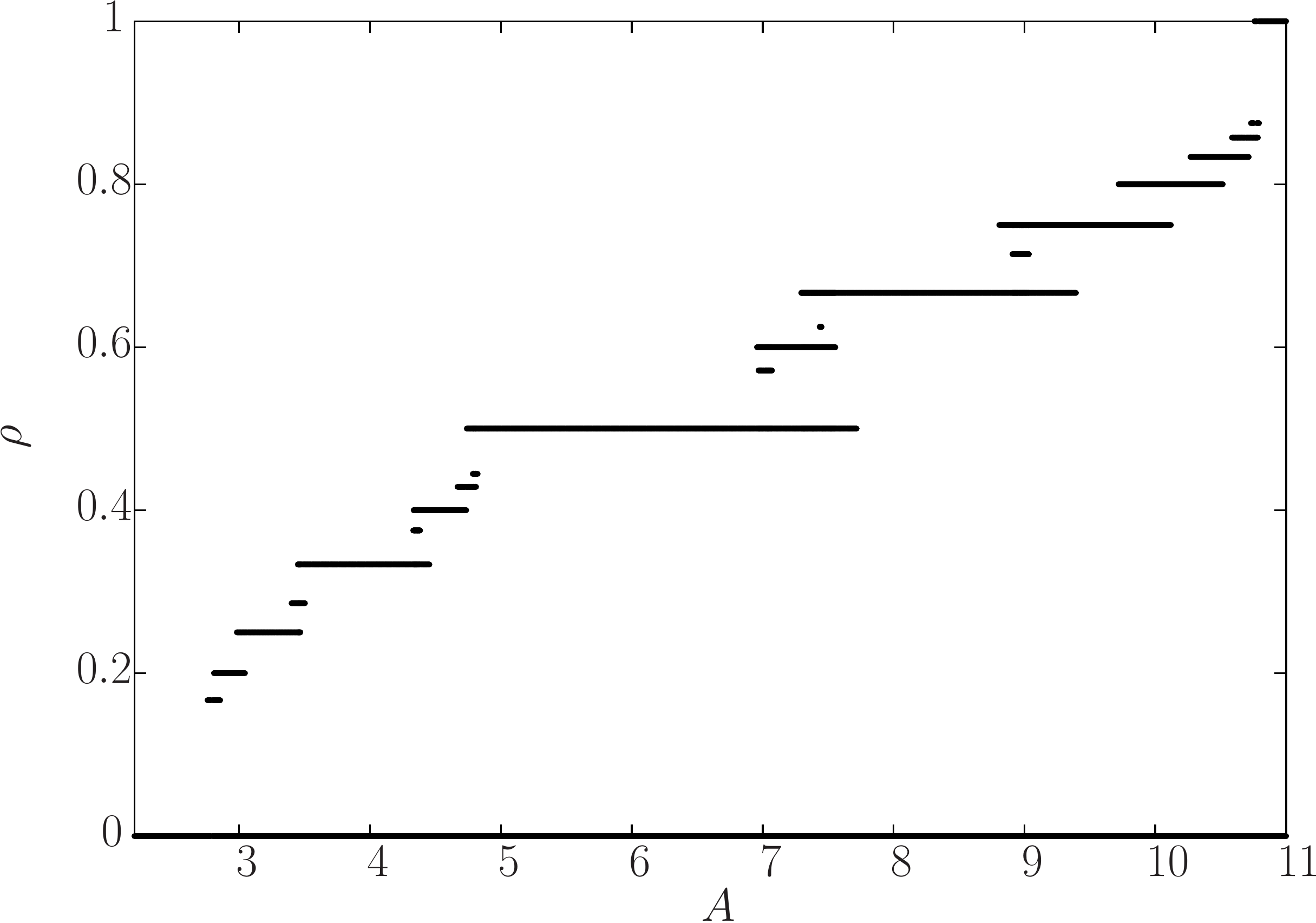}}
}
\end{picture}
\end{center}
\caption{Periods (top) and ``rotation numbers'' (bottom) of the
periodic orbits found by varying the amplitude $A$ along the vertical
lines $b=0.1$ (left) and $b=0.55$ (right) as indicated in Figure~\ref{fig:2dbif_T0d5}. 
Periodic orbits have been computed numerically using the algorithm described in
Appendix~\ref{ap:per_orb}.  For $b=0.1$, we repeated the computations
with a smaller stepsize along the $A$-axis and we have found higher
periods interleaved according to the Farey tree structure (results not
shown), suggesting that the ``rotation number'' shown in (c) might be continuous along the
Farey tree.}
\label{fig:1dscans}
\end{figure}

\begin{figure}
\begin{center}
\begin{picture}(1,0.7)
\put(0,0.69){
\subfigure[]{\includegraphics[angle=-90,width=0.5\textwidth]
{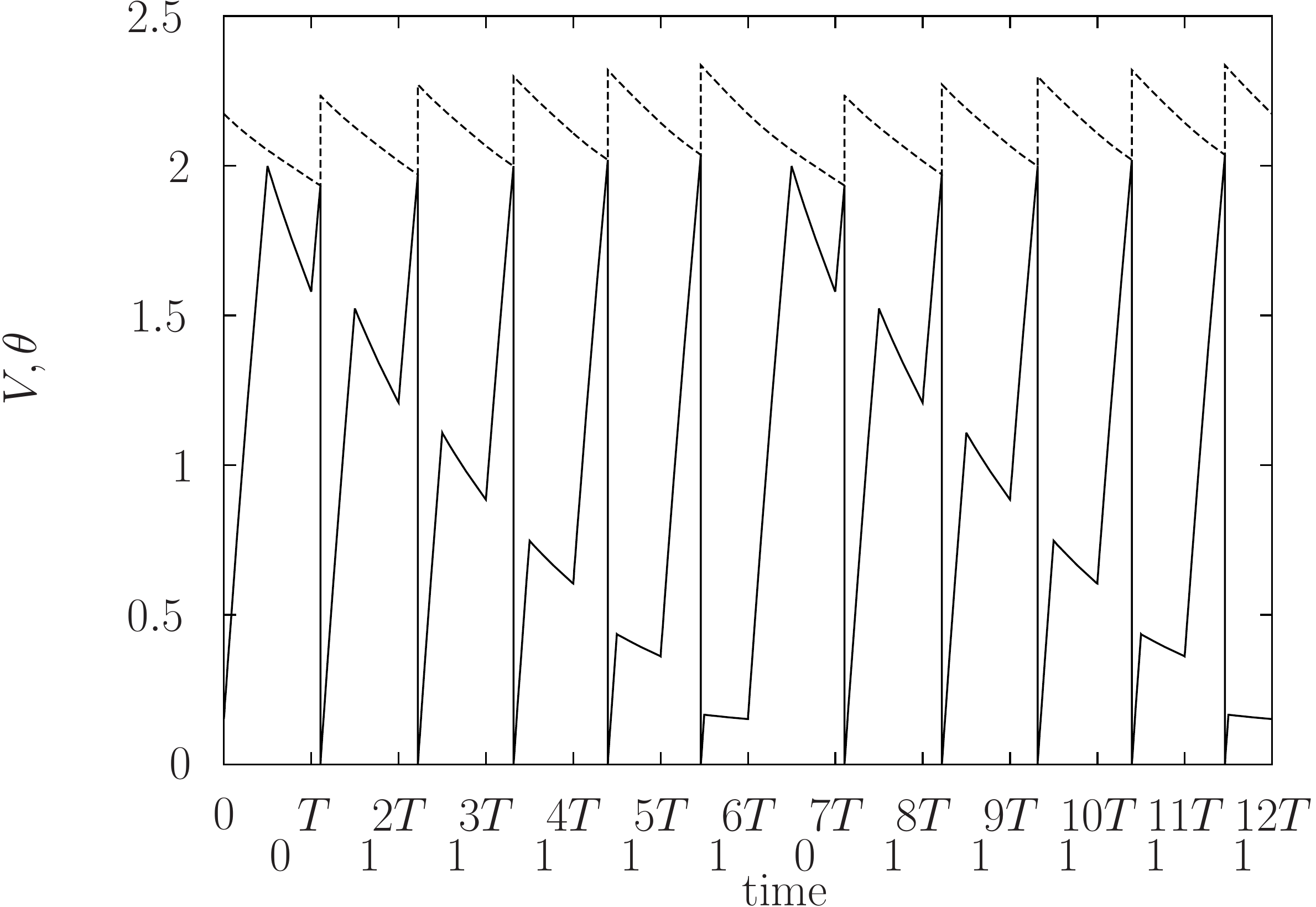}}
}
\put(0.5,0.69){
\subfigure[]{\includegraphics[angle=-90,width=0.5\textwidth]
{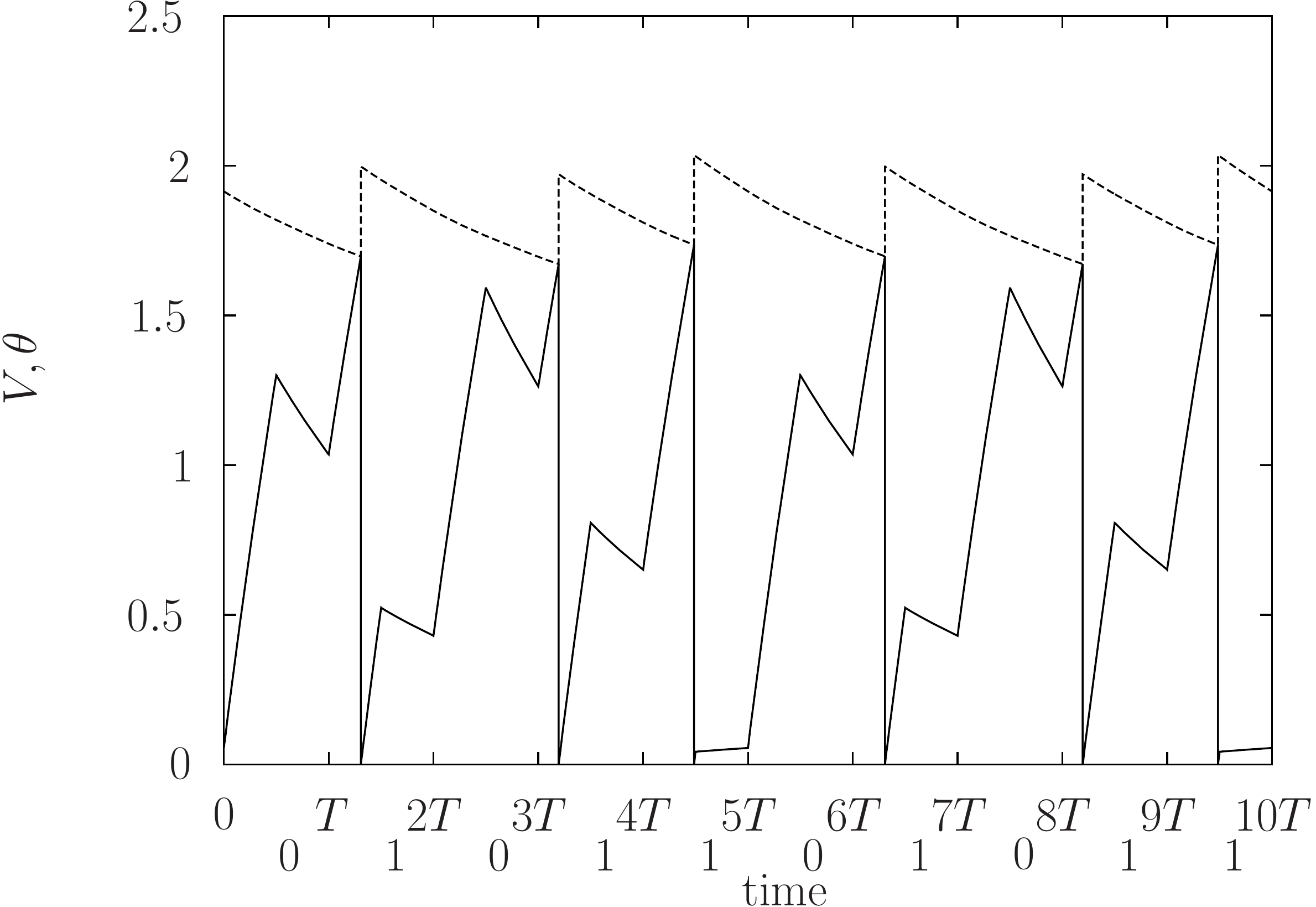}}
}
\put(0,0.31){
\subfigure[]{\includegraphics[angle=-90,width=0.5\textwidth]
{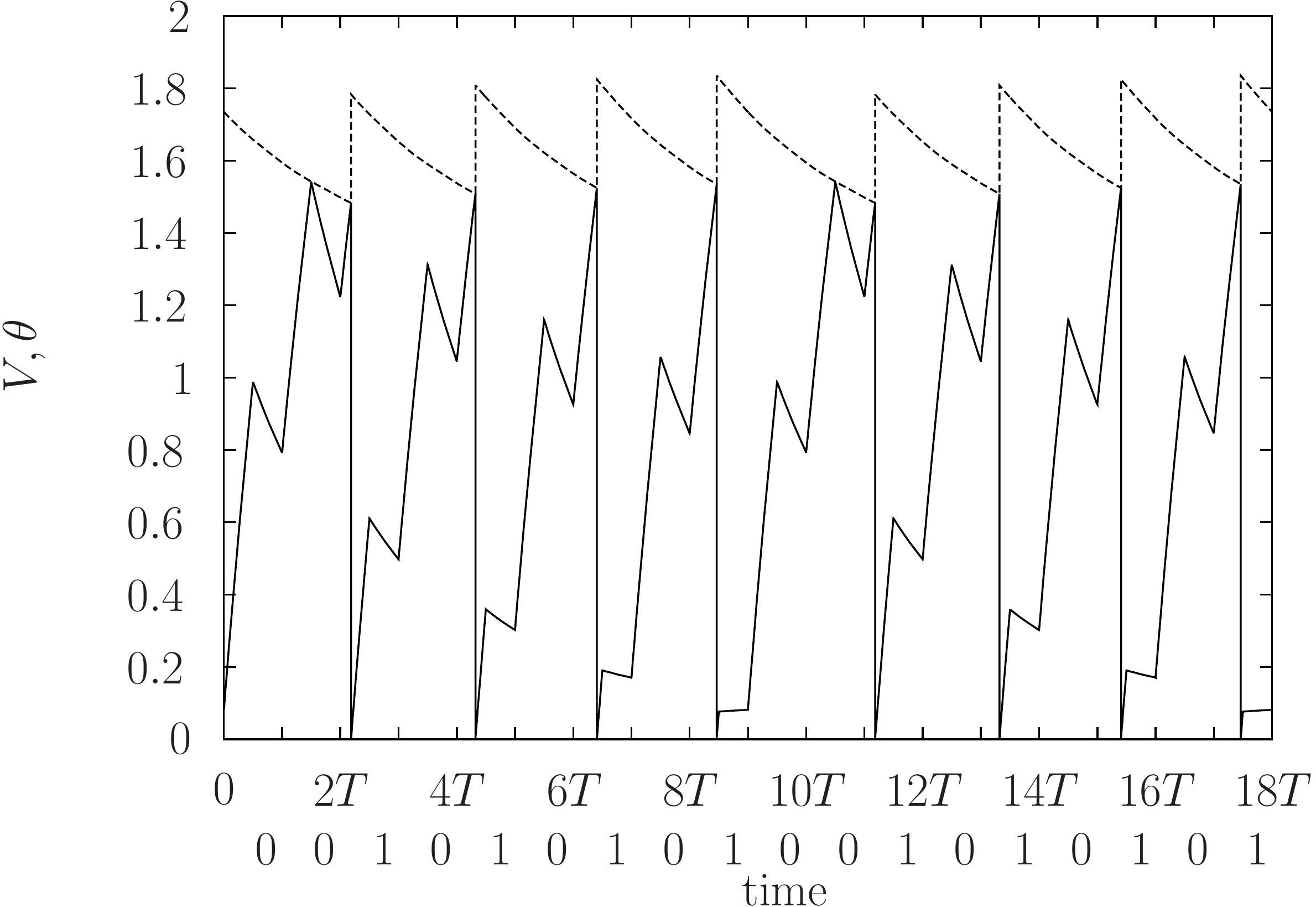}}
}
\put(0.5,0.31){
\subfigure[]{\includegraphics[angle=-90,width=0.5\textwidth]
{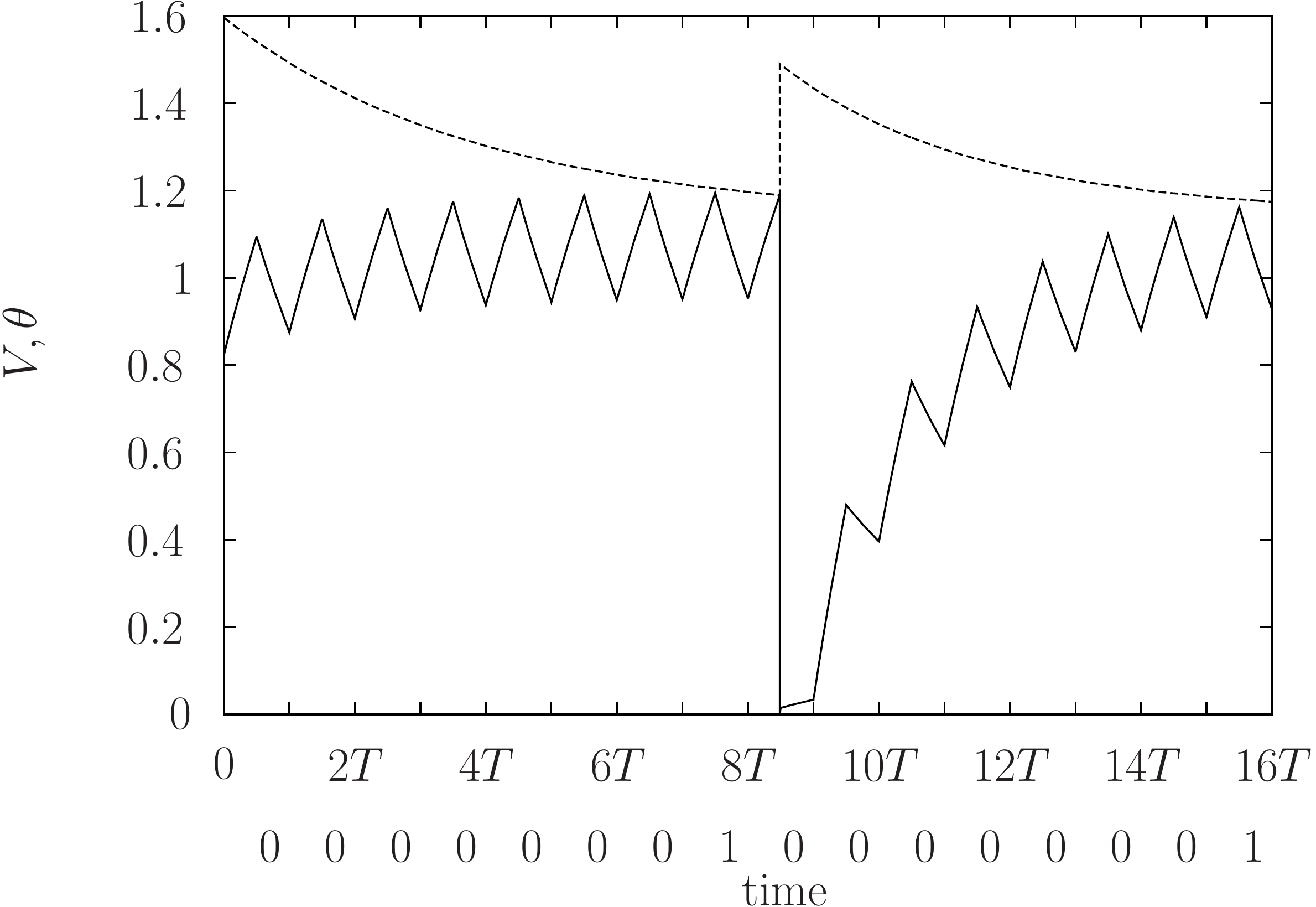}}
}
\end{picture}
\end{center}
\caption{Time course of the variables $V$ (solid) and
$\theta$ (dashed) corresponding to the four $nT$-periodic orbits
labeled in Figure~\ref{fig:periods_brut-force} (with $n=6$, $5$,
$9$ and $8$ for panels (a), (b), (c) and (d), respectively). Notice
that we show two periods of the $nT$-periodic orbits.
We indicate their symbolic itinerary at the bottom of the time axis.}
\label{fig:periodic_orbits}
\end{figure}

\begin{figure}
\begin{center}
\begin{picture}(1,0.7)
\put(0,0.67){
\includegraphics[angle=-90,width=0.5\textwidth]
{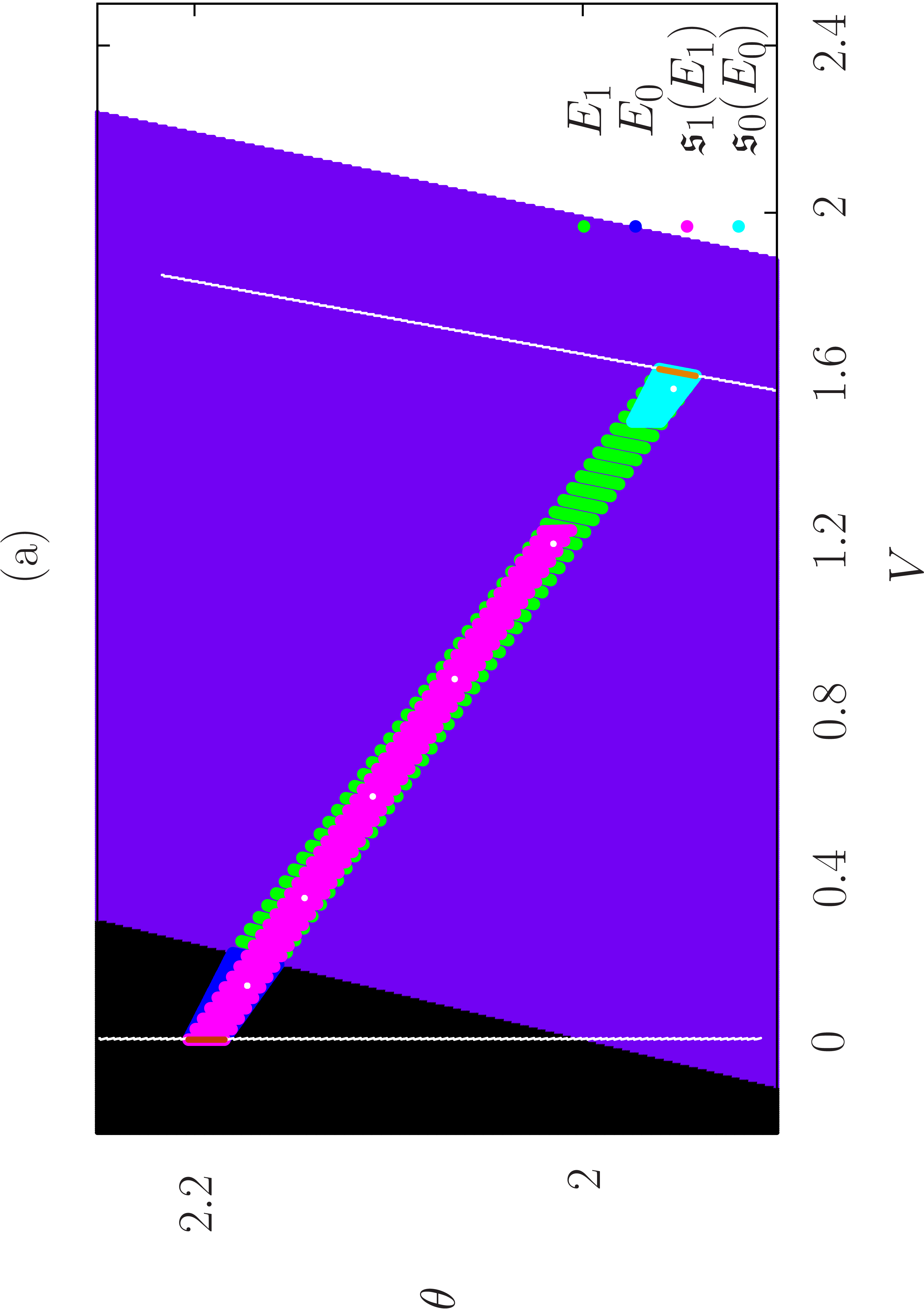}
}
\put(0.5,0.67){
\includegraphics[angle=-90,width=0.5\textwidth]
{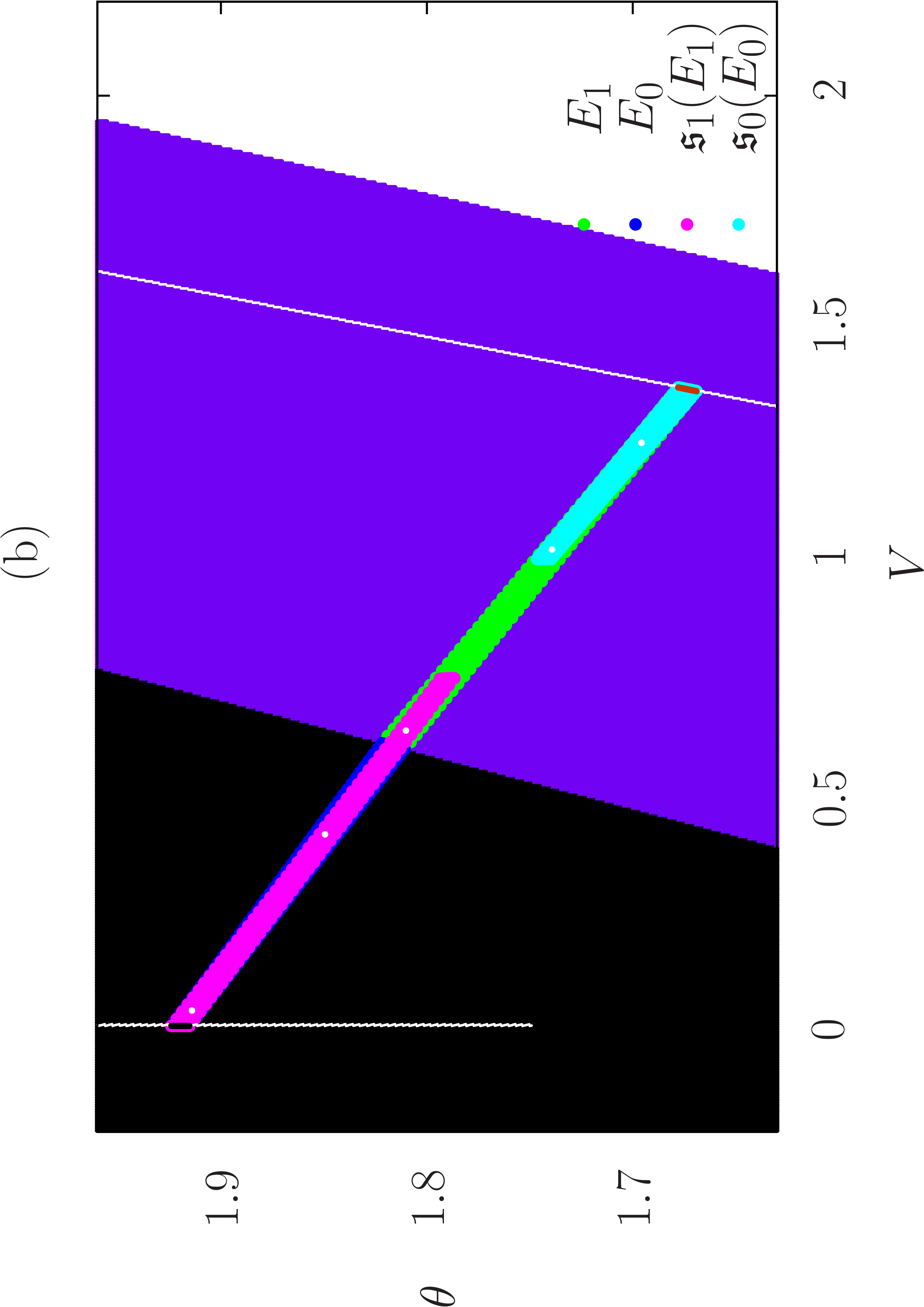}
}
\put(0,0.3){
\includegraphics[angle=-90,width=0.5\textwidth]
{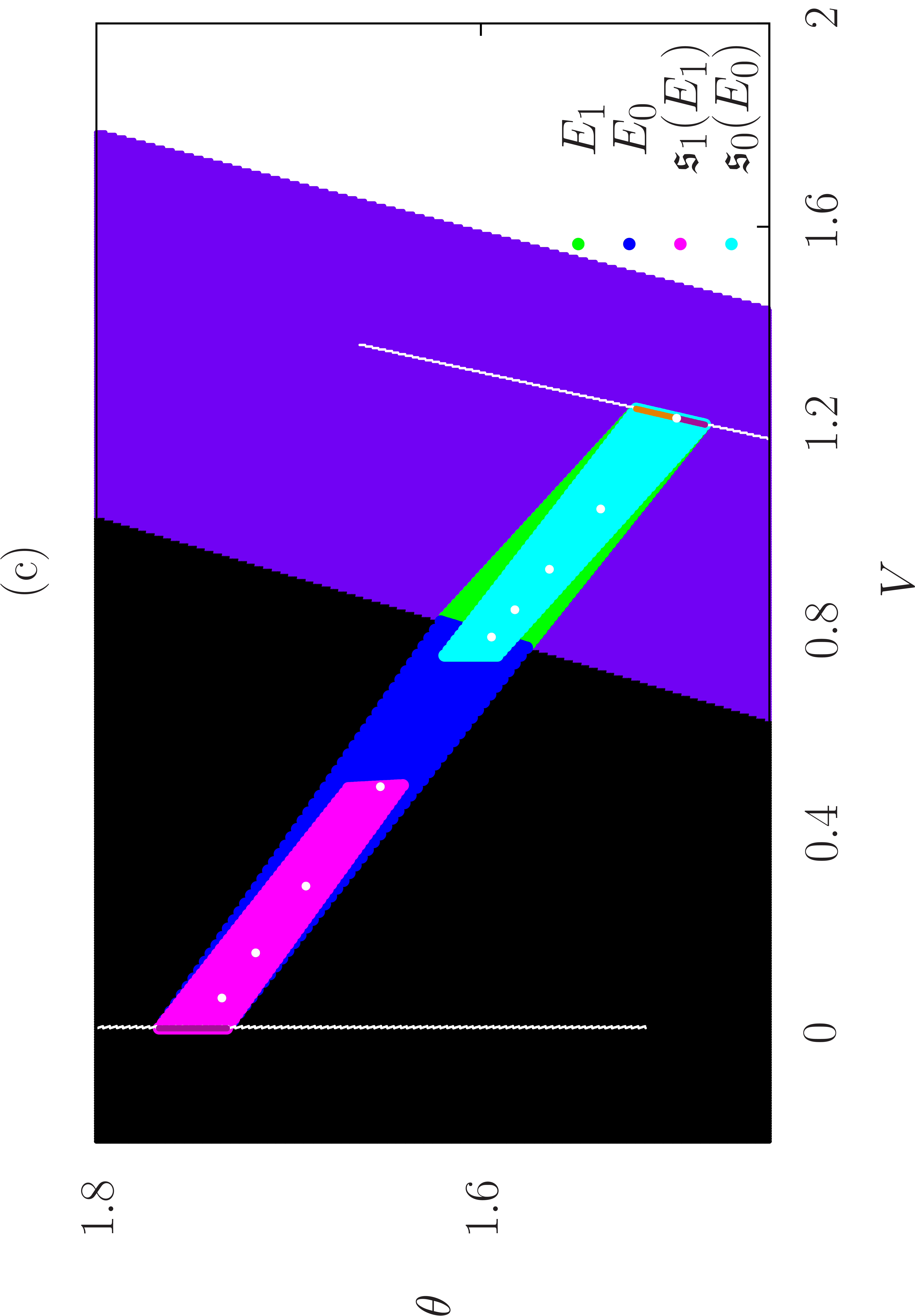}
}
\put(0.5,0.3){
\includegraphics[angle=-90,width=0.5\textwidth]
{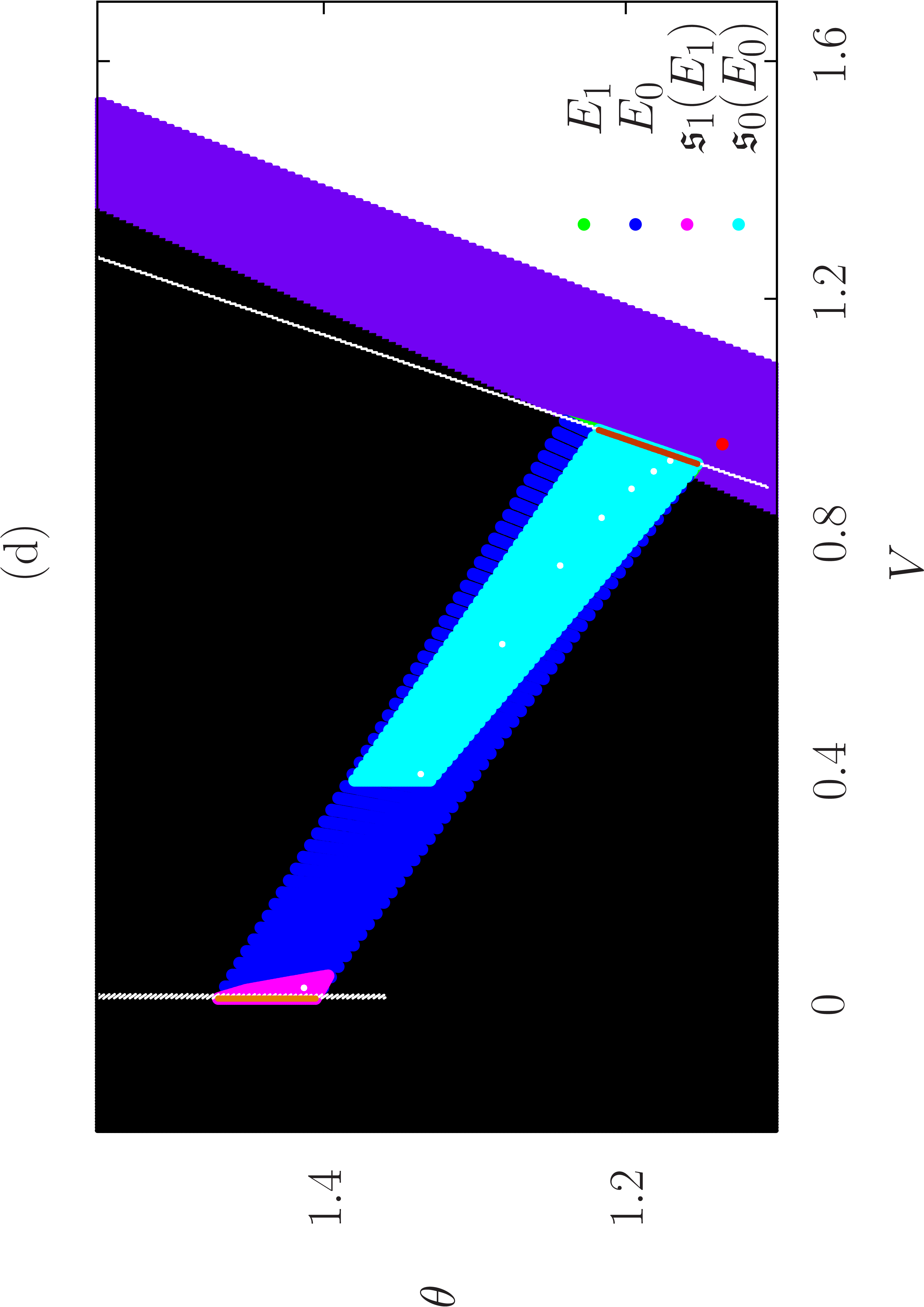}
}
\end{picture}
\end{center}
\caption{Sets $E_i$ satisfying condition \em{i)} of Theorem~\ref{theo:quasi-contraction} for the
parameter values labeled in Figure~\ref{fig:periods_brut-force}. The
found periodic orbits have symbolic itineraries: (a) $01^5$, (b)
$0101^2$, (c) $0^21\left( 01 \right)^3$ and (d) $0^71$. The
corresponding evolution of $V$ and $\theta$ with respect to time is
shown in Figure~\ref{fig:periodic_orbits}.  White points are the
periodic orbit found by direct simulation. Red point is a virtual
fixed point found by exteding the maps $\s_i$ in their virtual
domains, as explained in Section~\ref{sec:extension}.  Black
background is $S_0$, blue background is $S_1$. Blue and green sets are
$E_0$ and $E_1$, respectively. Pink and light blue are $\s_1(E_1)$ and
$\s_0(E_0)$, respectively. The white lines are the images
$\s_0(\Sigma_1)$ and $\s_1(\Sigma_1)$. Colored lines are the images of
$\gamma$ (see text) by $\s_i$. Those segments that stay connected for
all iterates of $\s$ are plotted with the same color.  For Figure~(c),
not all iterates of the $\Sigma_1$ stay connected and, hence, does not
satisfy condition~\emph{iii)}. For (c) and (d) we have found points in
$E_1$ whose differential $D\s_1$ has eigenvalues outside the unit
circle, and hence condition \em{ii)} is not satisfied. (a) and
(b) satisfy conditions {\em i)--iii)} of Theorem~\ref{theo:quasi-contraction}.}
\label{fig:sets_Ei}
\end{figure}

\begin{figure}
\begin{center}
\includegraphics[angle=-90,width=0.9\textwidth]{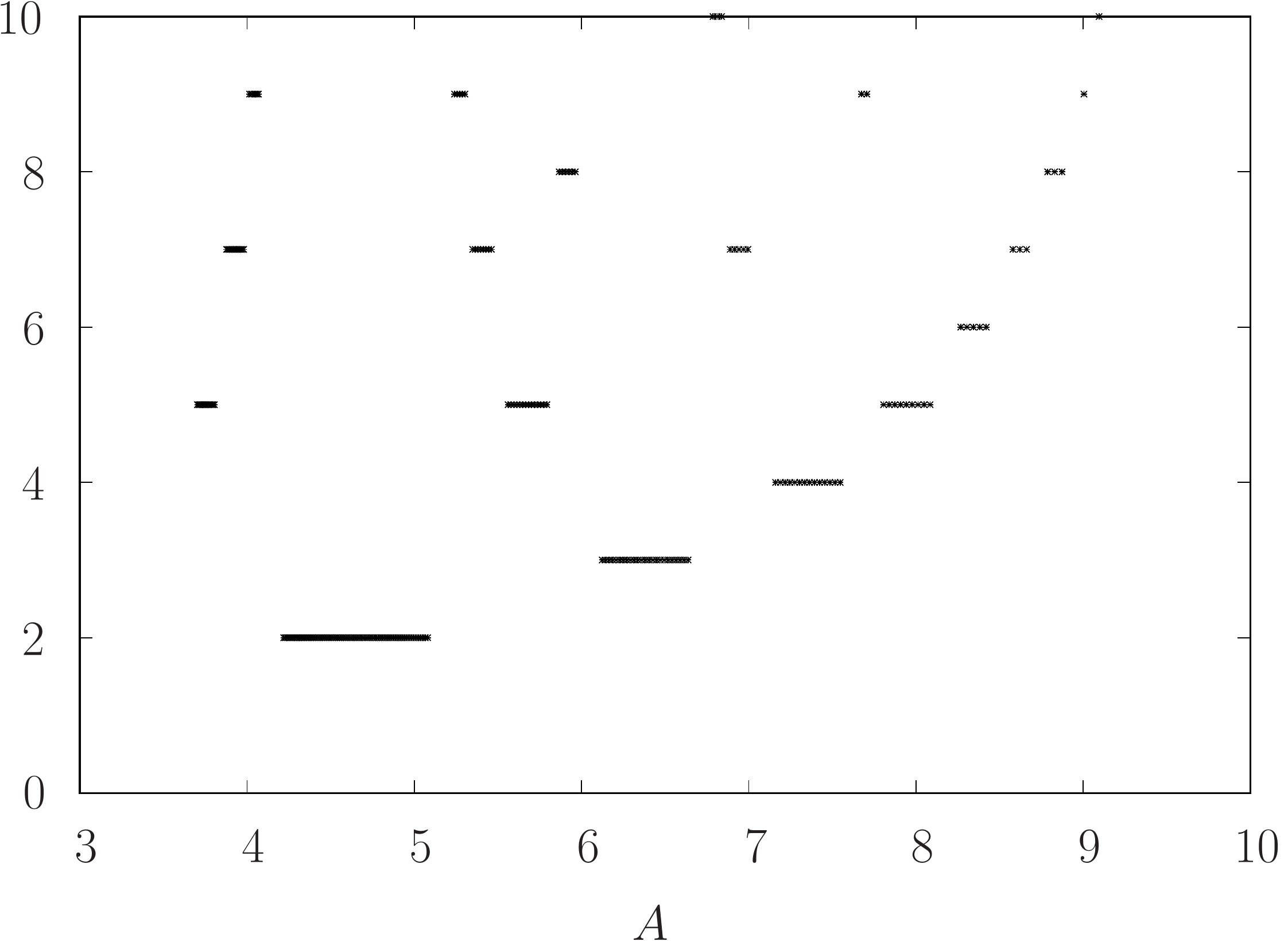}
\end{center}
\caption{Periods of the periodic orbits in Figure~\ref{fig:periods_brut-force} satisfying conditions of
Theorem~\ref{theo:quasi-contraction}.}
\label{fig:satisfying_Gam}
\end{figure}

\end{document}